\documentstyle[12pt]{article}
\begin{document}
\newtheorem{Def}{Definition}[section]
\newtheorem{conj}{Conjecture}[section]
\newtheorem{thm}{Theorem}[section]
\newtheorem{lem}{Lemma}[section]
\newtheorem{rem}{Remark}[section]
\newtheorem{prop}{Proposition}[section]
\newtheorem{cor}{Corollary}[section]
\newtheorem{clm}{Claim}[section]
\newtheorem{step}{Step}[section]
\newtheorem{sbsn}{Subsection}[section]
\renewcommand{\theequation}{\thesection.\arabic{equation}}
\newcommand{\add}{ \ \mbox{to}\ }
\newcommand{\qbar}{ \bar q }
\title
{The distance function to the boundary, Finsler geometry and
the singular set of viscosity solutions of
some Hamilton-Jacobi equations}
\author{YanYan Li\thanks{Partially
 supported by
  NSF grant DMS-0100819.
}\\
Department of Mathematics\\
Rutgers University\\
New Brunswick, NJ 08903\\
\\
Louis Nirenberg\\
Courant Institute\\
251 Mercer Street \\
New York, NY 10012\\
\\
Dedicated in memory to Jacques Louis Lions
}

\input { amssym.def}
\input epsf

\date{}
\maketitle
\setcounter{section}{0}

\section{Introduction  }
\noindent {\bf 1.1}\
This paper is concerned with viscosity solutions of
Hamilton-Jacobi equations of the form
\begin{equation}
H(x,u,\nabla u)=1\quad \mbox{in}\ \Omega,
\label{1.1}
\end{equation}
a $C^{2,1}$ bounded domain (connected open
set) in $\Bbb R^n$, and
\begin{equation}
H(x,t,p)\in C^\infty(\overline \Omega\times\Bbb R\times \Bbb R^n).
\label{1.2}
\end{equation}
We consider positive solutions $u$ satisfying
\begin{equation}
u|_{\partial \Omega}=0.
\label{1.3}
\end{equation}

For definitions and properties of viscosity solutions
we refer to \cite{L} and \cite{BC}.  Our main results are for special 
$H=H(x,p)$, i.e.,
\begin{equation}
H(x,\nabla u)=1\quad \mbox{in}\
\Omega
\label{1.1prime}
\end{equation}
under suitable conditions we show that the
$(n-1)$-dimensional Hausdorff measure of the
singular set of solution (the complement
of the open set where $u\in C^{1,1}$) is finite.

In addition, we prove the corresponding result for $H(x,t,p)$ but under
very special conditions.
See Theorem \ref{thm9.1} and, simple consequences, 
Proposition \ref{prop1.1}, \ref{prop1.2} and \ref{prop1.3}.

We were brought to the problem by first studying the singular set of
the distance function to the boundary of $\Omega$.  This set
is sometimes called the ridge of $\Omega$, or medial axes.  Our
interest in the set arises in connection with
nonlinear elliptic boundary value problems
(\cite{LN}). We first describe this set $\Sigma$.

Let $G$ be the largest open subset of $\Omega$ such that
every point $x$ in $G$ has a unique closest point on 
$\partial \Omega$.  The set $\Sigma$ is defined to be
$$
\Sigma=\Omega\setminus G.
$$
In $G$, the distance $u$ to the boundary is smooth
(i.e. of class $C^{1,1}$, or $C^\infty$ in case
$\partial \Omega$ is in $C^\infty$).

In case $\Omega$ is a ball, $\Sigma$ is just one point, its center. 
If we perturb the boundary of the ball by many small (but $C^\infty$) 
perturbations as in Fig. 1, we see that the set
$\Sigma$ consists of segments coming from
the origin

\bigskip

\centerline{ 
\epsfbox{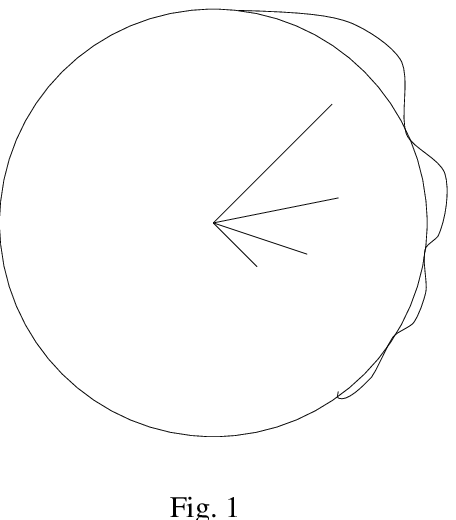}}

\bigskip

Another typical situation, with $\Omega$ not
simply connected is

\bigskip

\centerline{
\epsfbox{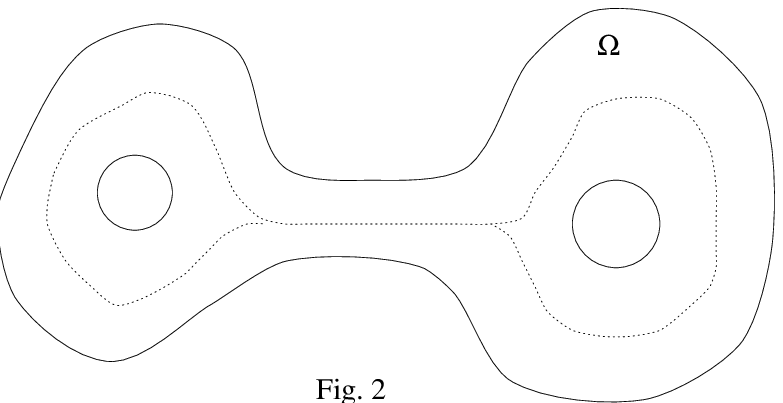}}

\bigskip

In this case $\Sigma$ is the dotted curve.

It is well known that $\Sigma$ is
always a connected set.  In Appendix C we will include 
a fairly short proof that it is
arcwise connected.

Concerning the set $\Sigma$ we proved that the
$(n-1)$-dimensional Hausdorff measure
$$
H^{n-1}(\Sigma)\ \mbox{is finite}.
$$
This is an immediate consequence of the following
result

\noindent{\bf Theorem A}\
{\it
 From every point $y$ on
$\partial \Omega$, move along
the inner normal until first hitting 
a point $m(y)$ on $\Sigma$.
The length $\bar s(y)$ of the resulting segment
is Lipschitz continuous in $y$.
}

\bigskip

\begin{rem}
The condition $C^{2,1}$ is sharp.  In Appendix
A, we present a convex domain $\Omega$ in the plane with $C^{2,\alpha}$
boundary, $0<\alpha<1$, for which
the conclusion of Theorem A does not
hold.
\label{rem1.1}
\end{rem}

If the domain $\Omega$ is unbounded the set $\Sigma$
may be empty, for example,
if $\Omega$ is the half-space $x_n>0$.  However,
the following form of Theorem A holds for general
$\Omega$ with $G$ and $\Sigma$ defined as before. 

\noindent {\bf Theorem A$'$}\
{\it For $y\in \partial\Omega$ let $\bar s(y)$ be
defined as in Theorem A (it may be infinite).  For
any $N>0$,
$\min(N, \bar s(y))$ is Lipschitz continuous in $y$ in
any compact subset of $\partial \Omega$.}

\medskip

After proving these theorems we extended
them to complete Riemannian manifold $(M,g)$.

\noindent {\bf Theorem A$''$}\
{\it
For any domain $\Omega$ in $M$, with $\partial\Omega$ 
locally in $C^{2,1}$, the 
conclusion of  Theorem A$'$
holds.  Here $\bar s(y)$ represents
the length of the geodesic going from
a point $y$ on $\partial \Omega$, normal to
 $\partial \Omega$, until it hits $\Sigma$.
}

\medskip

\begin{cor} For $\Omega$ as above in $(M^n,g)$,
$H^{n-1}(\Sigma\cap B)<\infty$ for any bounded
set $B$ in $M$.
\label{cor1}
\end{cor}

We then discovered that Theorem A$''$ had already 
been proved by J.I. Itoh and M. Tanaka
\cite{IT} in 2001.  In fact their domain
$\Omega$ may be the complement of a smooth submanifold
$X$ of $M$, of any dimension.  However the result
for such $X$ follows from the case $dim\ X=n-1$ by
taking for $\Omega$ the exterior of a tubular neighborhood
of $X$.

medskip

\noindent
{\bf Cut point}.\ In Theorem A$''$ we considered 
a geodesic from a point $y$ going into
$\Omega$ in the normal direction
until it first hits a point
$m(y)$ in $\Sigma$.  The point $m(y)$ is
called the cut point of $y$ on
$\partial\Omega$, meaning that if we go beyond $x$ on the geodesic, to any
point $x'$, then $x'$ has a closer point on $\partial
\Omega$ than $y$.
The collection
of these points $m(y)$ on $\Sigma$ for all $y$ (namely
$\Sigma$ itself) is called the  
{\it cut locus} of $\partial \Omega$.
That $\Sigma$ is the set of cut points is established in 
Section 4; see Corollary \ref{cor4-1}.

Recall the analogous notion of
{\it conjugate} point of $y$:  This is the first point
$\bar x$ on the normal geodesic such that any point
$x''$ on the geodesic beyond $\bar x$ has, in any
neighborhood of the normal geodesic, a point 
on  $\partial \Omega$ in the neighborhood which can
be connected to it by a path in the neighborhood
with length shorter than the arclength 
of the normal geodesic from  $y$ to it.
\begin{rem}
In case $\Omega$ is a domain in $\Bbb R^n$, the distance
from a point $y$ to the conjugate
point is the smallest of the principal radii 
of curvature of $\partial \Omega$ at $y$.
\label{rem1.2}
\end{rem}

In Corollary \ref{cor-new4}  we
give an analogous characterization for Finsler spaces.
It says that $m(y)$ is a conjugate point if and only if
 the (Finsler) sphere about  $m(y)$, of radius 
$s(y)$, has second order contact
                  with the boundary of $\Omega$  at
$y$ in some direction.
This result is not used
                  in the paper.

Remark \ref{rem1.2}  will be used in the construction given
in Appendix A.

Our proof of Theorem A$''$ is different from
that of \cite{IT}.  Some time ago Walter Craig
suggested that we might prove an analogue of
Corollary \ref{cor1} for viscosity
solutions of Hamilton-Jacobi equations
and we express our thanks to him.  The extension
is what we do in the paper.  As we learned,
to our surprise, for the problem (\ref{1.1prime}) and
(\ref{1.3}) it involved an extension
of Theorem A$''$ to Finsler 
geometry and we now proceed
to describe this.

\bigskip

\noindent{\bf 1.2.\ Hamilton-Jacobi equation}
Consider the problem
\begin{equation}
H(x,\nabla u)=1\quad\mbox{in}\ \Omega,
\label{1.4}
\end{equation}
\begin{equation}
u|_{\partial \Omega}=0.
\label{1.5}
\end{equation}
Here $H(x,p)\in C^\infty(\overline \Omega\times \Bbb R^n)$.
We assume that for every $x\in \overline \Omega$ the set
\begin{equation}
V_x=\{p\in \Bbb R^n\ |\ H(x,p)<1\}
\label{1.6}
\end{equation}
is a bounded convex surface containing $0$, with
smooth strictly convex boundary $S_x$
(i.e., having positive
principal curvatures).  For some
$r>0$ we assume that
\begin{equation}
B_r(0)\subset V_x\qquad
\forall\ x\in\overline\Omega.
\label{1.7}
\end{equation}

What is important are the sets $V_x$ rather than the particular function
$H(x,p)$.

Theorem 5.3 of \cite{L} gives an
explicit formula for the viscosity solution $u$
of (\ref{1.4}), (\ref{1.5}).
It involves, for each $x\in\overline \Omega$,
the support function $\varphi(x;\cdot)$ of $S_x$, i.e.
$$
\varphi(x;v)=\max\{v\cdot p\ |\ p\in
S_x\},\qquad v\in \Bbb R^n.
$$
The function $\varphi$ is in $C^\infty(\overline\Omega\times
(\Bbb R^n\setminus\{0\}))$, it is
positive homogeneous of degree $1$ in $v$,
is a convex function of $v$, in fact, for each $x\in \overline \Omega$,the set
$$
\{v\in \Bbb R^n\ |\ \varphi(x;v)=1\}
$$ is a smooth convex
hypersurface (with positive
principal curvatures) containing the origin in
its interior.
Furthermore, $\varphi$ satisfies
the triangle inequality in $v$.  Thus
for any curve $\xi(t)$, $0<t<T$, in
$\overline \Omega$
$$
\varphi(\xi(t);\dot \xi(t))dt
$$
is a Finsler metric. The length
of the curve, if $\dot\xi\in
L^1$, is
$$
\int_0^T \varphi(\xi(t);\dot \xi(t))dt.
$$
Because of the homogeneity it is independent of its
$t-$parameterization.

Note that the length of the curve depends on the
direction in which it is transversed, so we talk of
its length from $\xi(0)$ {\it to} $\xi(T)$.

For any $x,y\in\overline \Omega$ we denote
by $L(x,y)$ the infimum of length
of curves in $\overline\Omega$ going from
$y$ \underline{to} $x$,
\begin{eqnarray*}
L(x,y)&=&
\inf\bigg\{
\int_0^1 \varphi(\xi(t); \dot \xi(t))dt\ |\ 
\xi(t)\in \overline \Omega\ \mbox{for}\ 0\le t\le 1,\\
&&\qquad \dot\xi\in L^\infty(0,1)\ \mbox{and}\ \xi(0)=y, \xi(1)=x\bigg\}.
\end{eqnarray*}
Then for $x\in \overline \Omega$,
$$
u(x):=\inf_{y\in \partial \Omega}L(x,y)
$$
is the viscosity solution of (\ref{1.4}), (\ref{1.5}).
$u>0$ in $\Omega$ and $u\in W^{1,\infty}$.
See Theorem 5.3 in \cite{L}.

Thus the solution $u(x)$ is the distance from
$\partial \Omega$ to $x$ measured in the Finsler
metric.  What we do is to extend Theorem A$''$
to a general Finsler manifold.

\medskip

\noindent{\bf 1.3.}\
Consider an $n-$dimensional smooth manifold $M$ with
a complete, smooth Finsler metric.  Let
$\Omega$ be a domain in $M$ with
$$
\partial \Omega\in C^{2,1}_{loc}.
$$
Let $G$ be the largest open subset of $\Omega$ such that
for every $x$ in $G$ there is a unique 
closest point $y$
on $\partial \Omega$ to $x$; where we measure lengths of curves
in $\overline \Omega$ going from $\partial \Omega$
to $x$  in the 
Finsler metric.  It is easy to see that
the distance function from $\partial \Omega$
to $x$ is in $C^{1,1}(G\cup\partial\Omega)$.
Moreover $u$ belongs to $C^{k-1,\alpha}(G\cup\partial\Omega)$
if $\partial \Omega$ is $C^{k,\alpha}$
for $k\ge 3$ and $0<\alpha\le 1$.  But of course it never belongs to
$C^1$.

Set 
$$
\Sigma=\Omega\setminus G.
$$

As for Riemannian manifolds, $\Sigma$ is called
the cut locus of $\partial \Omega$.  The cut point
of $y$ on $\partial \Omega$ is defined
as in the Riemannian case, and the collection of
$m(y)$ for all $y\in \partial \Omega$ is
$\Sigma$ itself.
The cut point of $y$ on $\partial \Omega$ is usually defined 
 differently as follows.
We consider the geodesic from $y$ 
going into $\Omega$ in the ``normal'' direction
with unit speed, denoted
as $\xi(y, s)$.  The set of $s>0$ satisfying
$$
dist (\partial \Omega\add \xi(y,s))=s
$$
is either $(0, \infty)$ or $(0, \tilde s(y)]$
for some $0<  \tilde s(y)<\infty$.  In the latter case,
$\tilde m(y):=\xi(s, \tilde s(y))$
is the cut point of $y$ on $\partial \Omega$,
and the collection of $\tilde m(y)$ for all
$y\in \partial \Omega$, denoted as
$\tilde \Sigma$, is called the cut locus of $\partial \Omega$.
The two definitions are the same, i.e.
$\tilde m(y)=m(y)$ for all
$y\in \partial \Omega$, and $\tilde \Sigma
=\Sigma$.  This will be proved in Section 4.

The geodesic equations for the Finsler
metric $\varphi(\xi;v)$ are
\begin{equation}
\varphi_{\xi^i}(\xi(t);\dot\xi(t))=
\frac{d}{dt} \varphi_{v^i}(\xi(t);\dot\xi(t)),
\qquad i=1,\cdots,n.
\label{1.8}
\end{equation}
A $C^1$ solution, with nonvanishing
$\dot\xi$ is called a geodesic.  A geodesic locally minimizes
$$
\int_a^b \varphi(\xi(t);\dot\xi(t))dt.
$$

From any point $y$ on $\partial \Omega$ there is a
unique geodesic, in the metric, going into $\Omega$,
``normally'' at $\partial\Omega$.  This means
that for a point on the geodesic close to $y$, $y$
is the unique closest point on
$\partial \Omega$ to it.
This will be explained further below (see Lemma
\ref{lem2.2}).

\begin{thm} Let $\ell(y)$ denote the length
of the ``normal'' geodesic from $y$ until it first hits
a point $m(y)\in \Sigma$;
So $\Sigma=m(\partial \Omega)$.  
Then, for any $N>0$,
$$
\min(N, \ell(y))
$$
is Lipschitz continuous in $y$ on any
compact subset of $\partial\Omega$.
\label{thm1.1}
\end{thm}

\begin{cor} $H^{n-1}(\Sigma\cap B)<\infty$ for any
bounded set $B$.
\label{cor1.1}
\end{cor}

Returning to our viscosity solution of (\ref{1.4}), (\ref{1.5}), it means
that for its singular set $\Sigma$,
$$
H^{n-1}(\Sigma)<\infty.
$$

Some remarks on the general H-J equations
(\ref{1.1}) in a bounded $\Omega$.  Many authors have studied boundary
value problems
$$
u(x)=u_0(x)\quad\mbox{on}\ \partial \Omega.
$$
See for example papers below and references therein.
Usually it is considered that $H$ is convex in $p$.  Sometimes
it is also assumed that $H$ is convex in $(t,p)$.
And it is sometimes assumed 
that $H$ is nondecreasing in $t$; this
is usually used in proving uniqueness of the viscosity solution.
Adimurthi and Gowda (see \cite{AG1}, \cite{AG2} and references in it)
do not require $H$ nondecreasing in $t$.
In Theorem 5.5 of
\cite{L}, positive viscosity
solutions of (\ref{1.1}), (\ref{1.3}) are obtained
assuming $H$ is convex in $(t,p)$ and
nondecreasing in $t$ (and some additional conditions).

There are also a number of papers
which study the singular set of
solutions, which goes back at least to
\cite{T} by Ting.
A.C. Mennucci \cite{M}
studied the singular set for 
viscosity (and, what he calls ``minimal'') solutions
$u$ for the equation (\ref{1.1prime}) 
on a smooth $n-$dimensional manifold, with the value
of $u$ prescribed to be $u_0$ on a closed
subset $K$ of $M$.
$K$ and $u_0$ are usually assumed to be in $C^2$.
Among other things, he gives a
very fine characterization
of the set $A$ where the solution
$u$ is not differentiable, namely,
$A$ is the union of a
countable number of smooth
$(n-1)$-dimensional manifolds
with a set having zero $(n-1)$-dimensional 
Hausdorff measure.  Such sets are called ``rectifiable''.
This result does not contain ours,
since it does not show that
the total $(n-1)$-dimensional
measure is finite.
In an earlier paper \cite{MM}, he and
C. Montegazza
studied the distance function
to the boundary and showed
that the singular set
is ``rectifiable''
if $K$ is in $C^2$.
In addition, they presented an example
of a closed convex curve $K$ in $\Bbb R^2$,
$K$ of class $C^{1,1}$, such that
the singular set has positive Lebesgue
measure.  These papers contain
many more excellent results,
including some for the initial value problem,
as well  as many references to earlier work.

\bigskip

\noindent{\bf 1.4}.\  We wish to stress that what is 
important are the sets
\begin{equation}
V_x=\{(t,p)\in \Bbb R^{n+1}\ |\
H(x,t,p)<1\}\qquad\forall\
x\in \overline \Omega,
\label{1.9}
\end{equation}
and
$$
S_x=\partial V_x=\{(t,p)\ |\
H(x,t,p)=1\}\qquad\forall\
x\in \overline \Omega.
$$

For example, consider

\noindent{\bf Situation (*)}.\ Suppose $H$ is
smooth in a neighborhood of
$\displaystyle{\cup_x S_x}$ and that
$\forall\ x\in \overline \Omega$,
$V_x$ is convex and
$S_x$ is a smooth strictly convex
hypersurface with positive principal 
curvatures, and that
\begin{equation}
dist(0, S_x)\ge r_0>0\qquad
\forall\ x\in \overline \Omega.
\label{1.10}
\end{equation}
Suppose furthermore that each
$V_x$ lies in a fixed 
downward cone: 
for some $k, C_1>0$,
\begin{equation}
|p|\le k(C_1-t), \qquad t<C_1.
\label{1.11}
\end{equation}
Thus $t$ may be unbounded below in $V_x$.

Without loss of generality we may replace
the given $H$ by one that is
homogeneous in $(t,p)$ of
degree $1$.

\begin{rem} If $\tilde H(x,t,p)$ is another function satisfying
the condition above, with the same sets
$V_x$ as $H$, then 
a continuous viscosity solution of the problem
(\ref{1.1}), (\ref{1.3}) for
$H$ is also one for $\tilde H$----as is easily verified.
\label{rem1.3}
\end{rem}

For $H$ and $V_x$ as above, we take
$H$ to be homogeneous of degree one in
$(t,p)$, there is a viscosity solution.
See Claim \ref{claim10.1}.  However we do not
know if $H^{n-1}(\Sigma)<\infty$
for the singular set $\Sigma$.

In Section 10 we present a result, Theorem \ref{thm9.1}, with  
this picture, for which a
viscosity solution exists and its
singular set $\Sigma$ satisfies
$$
H^{n-1}(\Sigma)<\infty.
$$

Here are three special cases of that theorem.  In
the first two of these, $h(x,p)$ is
a function such that $\forall\ x\in\overline\Omega$,
$$
V(x)=\{p\ |\ h(x,p)<1\}
$$
is a bounded convex set with smooth boundary
$S_x$, strictly convex with
positive principal curvatures.  $h$ is assumed
to be smooth in a neighborhood of $\displaystyle{ \cup_x S_x }$.

\begin{prop}
There exists $\lambda_0>0$, depending on $h$ and on $\Omega$,
such that for any $0<\lambda<\lambda_0$, for the function
\begin{equation}
H(x,t,p)=\lambda t+h(x,p),
\label{1.12}
\end{equation}
problem (\ref{1.1}), (\ref{1.3}) has a positive viscosity
solution and its singular set $\Sigma$ satisfies
\begin{equation}
H^{n-1}(\Sigma)<\infty.
\label{1.13}
\end{equation}
\label{prop1.1}
\end{prop}

The existence of a positive viscosity solution
for any $\lambda>0$ is, of course,  part of
Theorem 5.4 in \cite{L}.  For large
$\lambda$ we have not succeeded in proving
(\ref{1.13}).

\begin{rem} One may ask what happens for $H$ given
in (\ref{1.12}) if $\lambda<0$.
Then there exists a negative viscosity solution, namely $u=-v$
where $v$ is the viscosity solution for
$$
\hat H=|\lambda| v+h(x,-\nabla v)=1
$$
as is easily verified.
\label{rem1}
\end{rem}

\begin{prop}
There exists $\epsilon_0>0$ depending on
$h$ and on $\Omega$, such that
$\forall\ 0<\epsilon<\epsilon_0$, for
$$
H=\epsilon t^2+h(x,p),
$$
problem (\ref{1.1}), (\ref{1.3}) has a viscosity solution
for which 
$$
H^{n-1}(\Sigma)<\infty.
$$
\label{prop1.2}
\end{prop}

\begin{prop}
Let $H(x,t,p)$, with corresponding $V_x$ and $S_x$,
satisfy the conditions of Situation (*) in a 
domain $\Omega$.  Then there exists a number $d_0>0$ 
depending on $H$ such that
if $\Omega'$ is any bounded subdomain of $\Omega$,
with $\partial \Omega'\in C^{2,1}$,
and such that the distance of any point $x$ in
$\Omega'$ to $\partial \Omega'$ is less than $d_0$
(i.e. $\Omega'$ is narrow)
then in $\Omega'$ the problem (\ref{1.1}), 
(\ref{1.3}) has a positive viscosity solution.
Furthermore, for its singular set $\Sigma$,
$$
H^{n-1}(\Sigma)<\infty.
$$
\label{prop1.3}
\end{prop}

The proofs of Proposition \ref{prop1.1}-\ref{prop1.3}
follow easily from Theorem \ref{thm9.1}
and will be presented in Section 10.

\medskip

We present one more proposition;
it is proved in Section 10.  Here we consider $H$ independent of
$x$,
$$
H=H(t,p)
$$
satisfying the conditions of Situation (*), in a bounded
domain $\Omega$.  Let $\bar t$ be the positive 
number satisfying $H(\bar t, 0)=1$ and
let
$$
\hat t=\max_{ H(0,p)=1}t;
$$
clearly $\bar t\le \hat t$.

\begin{prop}
Suppose $\bar t<\hat t$.  Then there is a positive
viscosity solution
of (\ref{1.1}), (\ref{1.3}) for this $H$, whose singular set $\Sigma$ satisfies
$$
H^{n-1}(\Sigma)<\infty.
$$
\label{prop1.4}
\end{prop}

In case $\bar t=\hat t$, we believe the same conclusion holds but,
as we explain in Section 10, our
method of proof cannot work.

\bigskip

\noindent{\bf 1.5.}\  Theorem \ref{thm9.1}, which concerns
general $H(x,t,p)$ is
derived from Theorem \ref{thm1.1}, where
$H$ does not involve $t$, by introducing
an extra independent variable
$\tau$ and by considering the function
\begin{equation}
z(\tau, x)=e^\tau u(x).
\label{1.13prime}
\end{equation}

We conclude the introduction by giving a brief
description of our proof of
Theorem \ref{thm1.1}.
For simplicity we assume $\overline \Omega$ is compact.

Consider a geodesic for the Finsler metric
$\varphi(\xi;v)$, starting at a point $y$ on $\partial \Omega$ and
going in the direction ``normal'' to $\partial \Omega$.
The geodesic is given by $\xi(t)$, with
$\xi(0)=y$ and satisfies
the geodesic equation
$$
\varphi_{\xi^i}(\xi(t);\dot\xi(t))
=\frac{d}{dt} \varphi_{v^i}(\xi(t);\dot\xi(t)).
$$
We may parameterize the geodesic using arclength
$s$, i.e.,
$$
\varphi(\xi(s);\dot\xi(s))\equiv 1.
$$
Denote the geodesic by 
$$\xi(y,s).$$ 
We have to explain the ``normal'' direction
$$
V(y)=\dot \xi(y,0).
$$
Let $\nu(y)$ be the unit inner 
normal to
$\partial \Omega$ at $y$.  Then
$V(y)$ is the unique
vector-valued function
on $\partial \Omega$ satisfying
\begin{equation}
\left\{
\begin{array}{l}
V(y)\cdot \nu(y)>0,\\
\varphi(y; V(y))=1,\\
\nabla_v\varphi(y, V(y))\
\mbox{is parallel to}\ \nu(y).
\end{array}
\right.
\label{1.14}
\end{equation}

From $y$ on $\partial \Omega$ we go along the geodesic
until we hit a point $m(y)$, set
$$
\bar s(y)=dist(y, m(y)).
$$

Without loss
of generality we may assume that $\bar s(\bar y)=1$,
i.e., $m(\bar y)=\xi(\bar y, 1)$.
We will show that 
there exist some large constant $K\ge 1$ and some small
constant $\delta>0$ such that
for all $y\in \partial \Omega$ satisfying
$0<|y-\bar y|\le \delta$, we can find 
$z=z(\bar y, y)\in \partial \Omega$ which satisfies
\begin{equation}
dist (z\add \xi(y, 1+K|y-\bar y|))<1+K|y-\bar y|=
 \bar s(\bar y)+K|y-\bar y|.
\label{co1}
\end{equation}

This implies that
$$
\bar s(y)\le \bar s(\bar y)+K|y-\bar y|,\qquad
\forall\  |y-\bar y|\le \delta.
$$
Since $K$ and $\delta$ are independent of 
$\bar y$ and $y$, we also have, by switching the roles
of $\bar y$ and $y$, that
$$
\bar s(\bar y)\le \bar s(y)+K|y-\bar y|,\qquad
\forall\ |y-\bar y|\le \delta.
$$
Thus
$$
|\bar s(y)-\bar s(z)|
\le K|y-z|,\qquad
\forall\ y,z\in \partial \Omega, |y-z|\le \delta.
$$
It follows, possibly for a larger $K$,  that
$$
|m(y)-m(z)|\le K|y-z|,\qquad
\forall\ y,z\in \partial \Omega, |y-z|\le \delta.
$$

To establish (\ref{co1}), 
we first use the triangle inequality
\begin{eqnarray*}
&&dist(z\add \xi(y, 1+K|y-\bar y|))
\\
&\le& dist(z\add  \xi(z, 1-K|y-\bar y|))
+ dist(\xi(z, 1-K|y-\bar y|)\add \xi(y, 1+K|y-\bar y|))\\
&\le& (1-K|y-\bar y|)+ dist(\xi(z, 1-K|y-\bar y|)\add \xi(y, 1+K|y-\bar y|)).
\end{eqnarray*}
We then construct a curve $\eta(t)$, $0\le t\le 1$,
satisfying
$$
\eta(0)=\xi(z, 1-K|y-\bar y|),\quad
\eta(1)= \xi(y, 1+K|y-\bar y|),
$$
and 
$$
\int_0^1\varphi(\eta(t); \dot \eta(t))dt<2K|y-\bar y|,
$$
from which we deduce
$$
dist(z\add \xi(y, 1+K|y-\bar y|))\le
 (1-K|y-\bar y|)+ \int_0^1\varphi(\eta(t); \dot \eta(t))dt<1+K|y-\bar y|.
$$

\bigskip

\centerline{
\epsfbox{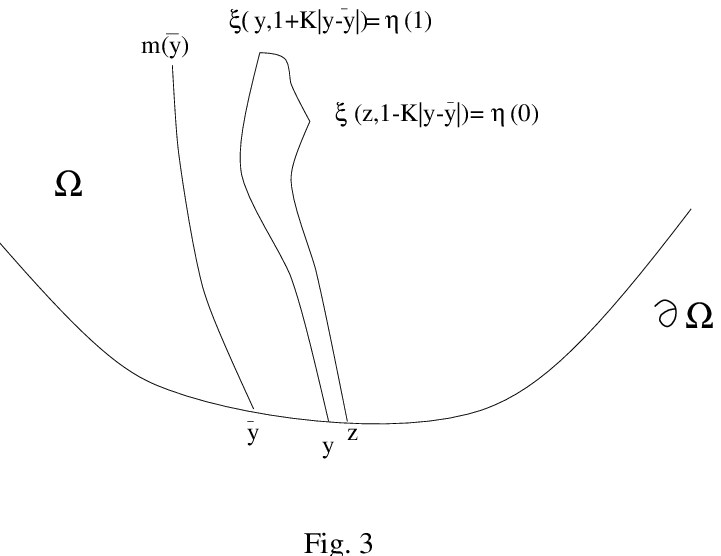}}

\bigskip

To construct the $\eta$, we make, for some small $\epsilon_0>0$,  
a diffeomorphism to map a neighborhood
of $\{\xi(\bar y, \tau)\}_{ -\epsilon_0\le \tau\le 1+\epsilon_0 }$
to  a neighborhood
of $\{\tau e_n\}_{ -\epsilon_0\le \tau\le 1+\epsilon_0 }$ so that 
in the new coordinates, $\{\tau e_n\}_{ -\epsilon_0\le \tau\le 1
+\epsilon_0 }$
is a geodesic for the new $\varphi$, and the new $\varphi$
has better properties.  Such new coordinates will be called  
special coordinates and they are produced in Section 3.
In the special coordinates, our $\eta$ is a straight segment 
connecting $\xi(z, 1-K|y-\bar y|)$ to  $\xi(y, 1+K|y-\bar y|)$.

\section{Preliminaries}
\setcounter{equation}{0}

\noindent{\bf 2.1.}\ It is convenient to
extend $\varphi$ so that
it satisfies 
\begin{equation}
\left\{
\begin{array}{l}
\varphi \in C^{2,1}\left(\Bbb R^n\times
(\Bbb R^n\setminus\{0\})\right),
\ \ \mbox{with derivatives smooth in}\ v\
\mbox{for}\ v\neq 0,
\\
\varphi(\xi; sv)\equiv s\varphi(\xi;v),
\qquad \forall\ s>0, \xi\in \Bbb R^n, v\in \Bbb R^n\setminus\{0\},
\\
\displaystyle{
0<\inf_{\xi\in \Bbb R^n, \|v\|=1}
 \varphi(\xi;v)\le \sup_{\xi\in \Bbb R^n, \|v\|=1}
 \varphi(\xi,v)<\infty,
}
\end{array}
\right.
\label{2.1}
\end{equation}
and
\begin{eqnarray}
0&<&\inf_{\xi\in \Bbb R^n, 
\|v\|=1, \|w\|=1} \frac{ \partial (\varphi^2) }
            { \partial v^i\partial v^j}(\xi;v)w^iw^j
\nonumber\\
&\le& \sup_{\xi\in \Bbb R^n, \|v\|=1, \|w\|=1}
 \frac{ \partial (\varphi^2) }
            { \partial v^i\partial v^j}(\xi;v)w^iw^j
<\infty.
\label{phi4}
\end{eqnarray}

Define,  for $x,y\in \Bbb R^n$,
$$
dist(y\ \mbox{to}\ x)=\inf\{ \int_0^1 \varphi(\xi(t), \dot \xi(t))dt\
|\ \xi(0)=y, \xi((1)=x, \dot \xi\in L^1(0,1)\}.
$$
Then $\Bbb R^n$, equipped with $dist(y\ \mbox{to}\ x)$, is a complete
(both forward and backward) Finsler manifold (see, e.g., \cite{BCS}).

Again, the geodesic equation for the Finsler metric is
$$
\varphi_{\xi^i}(\xi(t), \dot \xi(t))
=\frac {d }{d t} \varphi_{v^i}(\xi(t), \dot \xi(t)).
$$
We may always introduce a new $t$ variable so that
$$
\varphi(\xi; \dot \xi)\equiv 1,
$$
i.e. $t$ is arclength.

It is not difficult to see that 
$$
u(x)=\inf_{ y\in \partial \Omega }L(x,y)=
\inf_{ y\in \partial \Omega }dist(y\ \mbox{to}\ x),
\qquad
x\in \overline \Omega.
$$

Let
$$
\psi=\varphi^2.
$$
For $y\in \partial \Omega$, the vector $V(y)$ given in (\ref{1.14}) is simply 
$$
V(y)=\mu [\nabla_v\psi(y, \cdot)]^{-1}(\nu(y)),
$$
where
  $\mu>0$ is uniquely determined by
$$
\mu^2 \psi\left(y,  [\nabla_v\psi(y,  \cdot)]^{-1}(\nu(y))\right)=1.
$$

\bigskip

\centerline{
\epsfbox{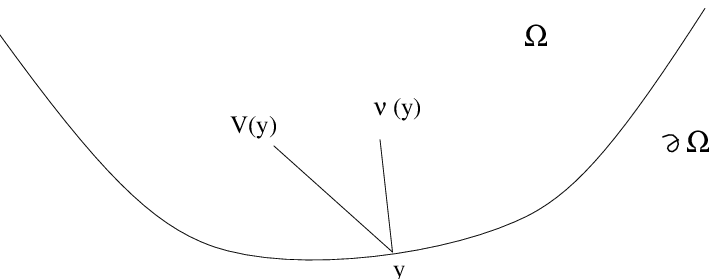}}

\bigskip

For $y\in \partial \Omega$, we consider 
 the following ODE:
$$
\psi_{\xi^i}(\xi(y,s); \dot \xi(y,s))
=\frac {\partial }{\partial s} \psi_{v^i}(\xi(y,s); \dot \xi(y,s)),
\qquad s\ge 0,
$$
\begin{equation}
\xi(y,0)=y, 
\label{ini1}
\end{equation}
and
\begin{equation}
\dot \xi(y,0)=V(y).
\label{ini2}
\end{equation}
Solutions $\xi(y,s)$ are geodesics 
starting from $y$ with unit speed, i.e.
$$\dot \xi(y,s)\neq 0,\quad 
\varphi(\xi(y,s); \dot \xi(y,s))\equiv 1,\qquad s\ge 0,
$$
and
$$
\varphi_{\xi^i}(\xi(y,s); \dot \xi(y,s))
=\frac {\partial }{\partial s} \varphi_{v^i}(\xi(y,s); \dot \xi(y,s)),
\qquad s\ge 0,
$$
with initial conditions (\ref{ini1}) and (\ref{ini2}).

For any $x,y\in \Bbb R^n$, 
 let
$$
X_1=\{\xi\in C([0,1],  \Bbb R^n)\ |\
\xi(0)=y, \xi(1)=x, \dot\xi\in L^1(0,1)\},
$$
$$
X_2=\{\xi \ \mbox{in}\ X_1 \ \mbox{with}\
 \dot\xi\in L^2(0,1)\},
$$
$$
I_1=\int_0^1\varphi(\xi(t); \dot\xi(t))dt, \qquad \xi\in X_1,
$$
and
$$
I_2=\int_0^1\varphi^2(\xi(t); \dot\xi(t))dt, \qquad \xi\in X_2.
$$

For any $\xi\in X_1$ and  any $t=t(\tau)\in C^1[0,1]$ 
satisfying  $t(0)=0$, $t(1)=1$ and $t'(\tau)>0, 
 0\le \tau\le 1$,
let $\eta(\tau)=\xi(t(\tau))$.  It  is easy to see that
$\eta\in X_1$ and
$$
I_1(\eta)=I_1(\xi).
$$

We list some elementary facts which can be found
in, e.g., \cite{BCS}.

\bigskip

\noindent {\it Fact 1.}\
 If $\bar \xi\in X_2$ is a critical
point of $I_2$, in the sense that
$$
\frac {d}{d\epsilon} I_2(\bar
\xi+\epsilon h)|_{\epsilon=0}=0, 
\qquad \forall\ h\in C_c^\infty((0,1),  \Bbb R^n).
$$
Then
 $\bar \xi$ belongs to $ C^\infty([0,1], \Bbb R^n)$,
$$
\dot {\bar \xi}(t)\neq 0,\qquad \forall\ 0\le t\le 1,
 $$
and $\bar \xi$
satisfies
$$
\psi_{\xi^i}(\bar \xi(t); \dot {\bar \xi}(t))
=\frac {d}{dt} \psi_{v^i}(\bar \xi(t); \dot {\bar \xi}(t)),
\qquad\mbox{on}\ [0,1],
$$
where $\psi=\varphi^2$.
Moreover,  if $\varphi$ is independent of
$\xi$, then 
$$
\bar
\xi(t)\equiv y+t(x-y).
$$

\bigskip

\noindent{\it Fact 2.}\
$$
I_1(\xi)\le \sqrt{  I_2(\xi) },\qquad
\forall\ \xi\in X_2.
$$

\bigskip

\noindent{\it Fact 3.}\ $\inf_{  X_1}I_1$ and $\inf_{X_2}  I_2$ are achieved,
and
$$
\inf_{  X_2}I_2=(\inf_{  X_1}I_1)^2.
$$
\bigskip

\noindent{\it Fact 4.}\
     Let $\bar \xi\in X_2$ be a minimum point of $I_2$,
i.e.
$$
I_2(\bar \xi)=\min_{  X_2 } I_2.
$$
Then $\bar \xi$ is also a minimum point of $I_1$, i.e.
$$
I_1(\bar \xi)=\min_{  X_1 } I_1.
$$
\bigskip

\noindent{\it Fact 5.}\
For $-\infty<a<b<\infty$, assume that 
$\xi\in C^2(a,b)$ satisfies
$$
\psi_{\xi^i}(\xi;\dot \xi)=\frac {d}{dt}
\psi_{v^i}(\xi;\dot \xi),\qquad \mbox{on}\ (a,b),
$$
where, as usual, $\psi=\varphi^2$.
Then 
$$
\frac {d}{dt}\psi(\xi;\dot \xi)\equiv 0,\qquad \mbox{on}\ (a,b),
$$
and, consequently,  $\xi$ satisfies the geodesic equation
$$
\varphi_{\xi^i}(\xi;\dot \xi)=\frac {d}{dt}
\varphi_{v^i}(\xi;\dot \xi),\qquad \mbox{on}\ (a,b).
$$
Moreover, either $\dot\xi\equiv 0$ on $(a,b)$ or
$
\dot\xi(t)\neq 0$ for all $t\in (a,b)$.

\bigskip

The following is a simple but useful lemma.
\begin{lem} Let $\xi(s,\sigma)$ be a $C^1$
family of geodesics with $s$ as arclength, depending on
some parameters $\sigma=(\sigma_1, \cdots, \sigma_k)$ and assume
that $\xi_{\sigma_\alpha}$ are twice
continuously differentiable in $s$.  Then
\begin{equation}
\frac {\partial }{ \partial s}\left(
\xi_{\sigma_\alpha}^i \varphi_{v^i}(\xi;\dot \xi)\right)\equiv 0.
\label{2.5}
\end{equation} 
Here $\dot{\ }=\partial_s$.
\label{lem2.1}
\end{lem}

\noindent{\bf Proof.}\ Differentiating
$$
\varphi(\xi; \dot \xi)=1
$$
with respect to $\sigma_\alpha$ we find
$$
\varphi_{ \xi^i}\xi^i_{\sigma_\alpha}
+\varphi_{v^i}\dot \xi^i_{\sigma_\alpha}=0.
$$
Identity (\ref{2.5}) then follows with the aid of the geodesic equations.

\vskip 5pt
\hfill $\Box$
\vskip 5pt

\noindent{\bf 2.2.}\ We now turn to a point
on $\partial \Omega$.  We may
assume it is the origin, and that $\Omega$ is given
by 
$$
x_n>f(x'),\qquad x'\in \Bbb R^{n-1}
$$
with $f$ a $C^{2,1}$ function defined
on $|x'|\le \epsilon_1$, with
$$
f(0')=0, \quad \nabla f(0')=0.
$$
Throughout, when we say that some constant
depends on $f$ we mean it depends on
the $C^{2,1}$ norm of $f$:
$$
\|f\|_{ C^{2,1} }=
\|f\|_{ C^2 }+\sup_{ x'\neq y'}
\frac {  |D^2f(x')-D^2f(y')|   }
{|x'-y'| }.
$$

We consider geodesics $\xi=\xi(x',s)$
which are $C^{1,1}$ functions of $x'$ and $s$,
with $\nabla_{x'}\xi$ smooth in $s$, with unit speed
starting at $z=(x', f(x'))$  i.e.,
$\xi$ satisfies
\begin{equation}
\varphi_{\xi^i}(\xi;\dot \xi)=\frac {\partial}{\partial s}
\varphi_{v^i}(\xi;\dot \xi),\qquad |x'|\le \epsilon_1,\
0\le s<a,
\label{2.6}
\end{equation}
\begin{equation}
\varphi(\xi; \dot\xi)\equiv 1,\qquad |x'|<\epsilon_1,
0\le s< a,
\label{2.7}
\end{equation}
and
$$
\xi(x',0)=z=(x', f(x')),\qquad |x'|<\epsilon_1,
$$
and entering $\Omega$,
$$
\dot \xi(x',0)\cdot (-\nabla f(0'), 1)>0.
$$

We have changed notation: before the
geodesic $\xi(x',s)$ was
denoted by $\xi( (x', f(x')), s)$.

\begin{lem} 
Suppose  that for some
fixed $w=(x', f(x'))$, and
$\bar s$ small, $w$ is the closest point on $\partial \Omega$
to $\xi(x', \bar s)$.  Then
\begin{equation}
\dot \xi(x', 0)=V(x'),
\label{2.8}
\end{equation}
where $V(x')$ is the vector satisfying (\ref{1.14}) i.e.
$$
V(x')\cdot (-\nabla f(x'), 1)>0,
$$
$$
\psi(w; V(x'))=1,
$$
$$
\nabla _v\psi(w; V(x'))\
\mbox{is parallel to}\
 (-\nabla f(x'), 1).
$$
The vector $V(x')$ is simply
$$
V(x')=\mu[\nabla_v\psi(w; \cdot)]^{-1}(-\nabla f(x'), 1)
$$
with $\mu$ determined by
$$
\psi(w; V(x'))=1.
$$
Here we have abused the notation a little since
by our earlier convention, $V(x')$ should be
denoted as $V(w)$.
\label{lem2.2}
\end{lem}

\noindent{\bf Proof.}\ For any
$0<s<\bar s$, $w$ is the closest point on
$\partial \Omega$ \underline{to} $\xi(x',s)$, so we may take $\bar s$
so small
that for every $y'$ close to $x'$ there is a minimal geodesic
$\eta(y', t)$,
$0\le t\le \bar s$, with
\begin{equation}
\eta(y', 0)=(y', f(y')),
\quad \eta(y', \bar s)=\xi(x', \bar s).
\label{2.9}
\end{equation}

Note that except for $\eta(x',t)$, $t$ may not be arc length
on the geodesics $\eta$.  By assumption,
$$
\int_0^{\bar s}\varphi(\eta(y',t); \dot \eta(y',t))dt
$$
has a minimum at $y'=x'$; so at $x'$,
for $\alpha<n$, its
$y_\alpha$-derivative is zero:
\begin{eqnarray}
0&=& \int_0^{\bar s}
\varphi_{\xi^i}(\eta;\dot\eta)\eta_{y_\alpha}^i+
\varphi_{v^i}(\eta;\dot\eta)\dot \eta_{y_\alpha}^idt\nonumber\\
&=&  \int_0^{\bar s}
\frac{\partial}{\partial t}[\varphi_{v^i}  \eta_{y_\alpha}^i]dt
=(\varphi_{v^i}  \eta_{y_\alpha}^i)(\bar s)
-(\varphi_{v^i}  \eta_{y_\alpha}^i)(0).
\label{2.10}
\end{eqnarray}
Here we have used the geodesic equations
satisfied by $\eta$.  By (\ref{2.9}),
$$
\eta_{y_\alpha}^i(y', \bar s)\equiv 0.
$$
Also, for $1\le \alpha, \beta\le n-1$,
\begin{equation}
\left\{
\begin{array}{rll}
\xi_{x_\alpha}^\beta(x',0)&=&\eta_{y_\alpha}^\beta(x',0)=\delta_\alpha^\beta,
\\
\xi_{x_\alpha}^n(x',0)&=&\eta_{y_\alpha}^n(x',0)=f_{x_\alpha}(x').
\end{array}
\right.
\label{2.11}
\end{equation}
Inserting these into (\ref{2.10}) we find, for $\alpha\le n-1$,
$$
\varphi_{v^\alpha}(\xi(x',0); \dot \xi(x',0))+
f_{x_\alpha}\varphi_{v^n}(\xi(x',0); \dot \xi(x',0))=0,
$$
i.e.,
\begin{equation}
\nabla_v\varphi(z; \dot \xi(x',0))\ \mbox{is
parallel to}\ (-\nabla f(x'), 1)
\label{2.12}
\end{equation}
so (\ref{2.8}) is proved.

\vskip 5pt
\hfill $\Box$
\vskip 5pt

Note that, from (\ref{2.11}),
\begin{equation}
\xi_{x_\alpha}^i(x',0)\varphi_{v^i}(\xi(x',0); \dot \xi(x',0))=0.
\label{2.13}
\end{equation}

In the following we continue to  use $\xi(x', s)$ to denote 
the solution of
$$
\psi_{\xi^i}(\xi(x',s); \dot \xi(x',s))
=\frac {\partial }{\partial s} \psi_{v^i}(\xi(x',s); \dot \xi(x',s)),
$$
$$
\xi(x',0)=(x', f(x')), 
$$
and
$$
\dot \xi(x',0)=V(x').
$$
By the choice of $V(x')$, $\psi(\xi(x',0);\dot \xi(x',0))=1$, so,
by Fact 5, $\psi(\xi(x',\cdot); \dot \xi(x',\cdot))\equiv 1$.
By the smooth dependence of solutions of ODEs
on initial datas, we have, for some smooth $\chi$,
that $\xi(x',s)=\chi( (x',f(x')), V(x'), s)$.
Since $f$ is in $C^{2,1}$, $V(x')$ is in $C^{1,1}$, and therefore,
 for some constant $E$, depending
only on $\varphi$, $f$
and $a$, we have, for all $1\le \alpha, \beta\le n-1,
|x'|\le \epsilon_1,$ and $-\epsilon_1\le s\le a,$  that
$$
\sum_{k=0}^3
(|\frac{\partial^k}{\partial s^k}\xi(x', s)|+
|\frac{\partial^k}{\partial s^k}
 \xi_{x_\alpha  }(x', s)|+
|\frac{\partial^k}{\partial s^k}
 \xi_{x_\alpha   x_\beta  }(x',s)|)
\le E,
$$
and
$$
\sum_{k=0}^3|\frac{\partial^k}{\partial s^k}
 \xi_{x_\alpha  }(x', s)
-\frac{\partial^k}{\partial s^k}
 \xi_{x_\alpha  }(0', s)|\le E|x'|, 
$$

The conditions of Lemma \ref{lem2.1} therefore
hold, and it follows from the lemma, and
(\ref{2.13}), that
\begin{equation}
\xi_{x_\alpha}^i(x',s)\varphi_{ v^i}
(\xi(x', s); \dot \xi(x', s))\equiv 0.
\label{2.14}
\end{equation}

We now show, in some sense, the converse of Lemma \ref{lem2.2}.
 
\begin{lem} Consider $|x'|\le \epsilon_1$.  
For some positive constant $\epsilon_2$, depending only on
$\varphi$ and  $f$,
 we have 
$$
dist(0\add \xi(0', s))
<  dist((x', f(x'))\add \xi(0', s)),\quad
\forall\ 0<s<\epsilon_2, \ 0<|x'|\le \epsilon_1,
$$
and
$$
dist(\xi(0', s)\add 0)<dist(\xi(0', s)\add (x', f(x'))),\quad
\forall\ -\epsilon_2<s<0, \ 0<|x'|\le \epsilon_1.
$$
\label{lem2.3}
\end{lem}

\noindent{\bf Proof.}\  For simplicity we assume $s>0$.
 There exists $\epsilon_2>0$,
depending only on
$f$ and $\varphi$, such that
$$
\psi_{\xi^i}(\xi(x',s); \dot \xi(x',s))
=\frac {\partial }{\partial s} \psi_{v^i}(\xi(x',s); \dot \xi(x',s)),
\qquad |x'|\le \epsilon_1/2, \ |s|\le 2\epsilon_2,
$$
$$
\xi(x',0)=(x', f(x')), \qquad |x'|\le \epsilon_1/2,
$$
and
$$
\dot \xi(x',0)=V(x'), \qquad |x'|\le \epsilon_1/2.
$$
has unique smooth solutions. Moreover,
for any $|x'|\le \epsilon_1/2$, $\xi(x'; s)$ is shortest geodesic
for $|s|\le \epsilon_2$.  From Lemma \ref{lem2.2} and
(\ref{2.13}) we see that for $|x'|<\epsilon_1$,
the Jacobian of the map $(x',s)\to \xi(x',s)$ is
positive at $s=0$.  Hence for $\epsilon_2$ small,  
the map $(x', s)\to \xi(x';s)$ is a diffeomorphism
for $|x'|\le \epsilon_1/2$ and $|s|\le \epsilon_2$,
and 
$$
s=dist(0\add \xi(0',s))<dist((x', f(x'))\add
 \xi(0',s)),\qquad\forall\ |x'|\ge
\epsilon_1/4,
$$
and
$$
dist(0\add  \xi(0',s))=\min_{|x'|\le \epsilon_1/4}
dist((x', f(x'))\add  \xi(0',s)).
$$

Let $\bar x'$ be a minimum point, i.e.,
$|\bar x'|\le  \epsilon_1/4$ and
$$
s=dist(0\add  \xi(0',s))=dist( (\bar x', f(\bar x'))\add  \xi(0',s)).
$$
By Lemma \ref{lem2.2},
$$
 \xi(0',s)=\xi(\bar x', s).
$$
Since the map $(x', s)\to \xi(x',s)$ is a diffeomorphism,
 we must have
$\bar x'=0'$.
Lemma \ref{lem2.3} is established.

\vskip 5pt
\hfill $\Box$
\vskip 5pt

\section{Special Coordinates}
\setcounter{equation}{0}

Let $\varphi(\xi;v)$ be as in Section 2, 
and let $\xi=\xi(t)$ be a geodesic with $\dot \xi(t)\ne 0$ and  
$$
\varphi_{\xi^i}(\xi(t); \dot \xi(t))=
\frac{ d}{dt} \varphi_{v^i}(\xi(t); \dot \xi(t)),
\qquad \forall\ 1\le i\le n.
$$
For a non-singular change of variables
$\xi=\xi(\eta)$ in $\Bbb R^n$, let
$$
\widetilde \varphi(\eta;w)=\varphi(\xi(\eta); \xi_{\eta} w)),
$$
where $\xi_{\eta}:= \left\{\frac {\partial \xi^i}{  \partial \eta^j}\right\}$.
Such a change of variables maps geodesics to geodesics.

With $\partial \Omega$ locally as in Section 2.2, so
that $\nu(0)=e_n=(0,\cdots,0,1)$, we
consider the geodesics $\xi(x',s)$
of that section.  In view of the above one may
make a smooth change of variables so that in the
new variables the geodesic
$\xi(0',s)$, with $s$ as arc length, runs on the
$x_n-$axis and such that we still have
$\nu(0)=e_n$.  We start with this situation.

Throughout, Greek letters, $\alpha, \beta$, run from
$1$ to $n-1$, while indices
$i,j,k$ etc. run from $1$ to $n$.

\begin{lem}  Let $\{t e_n\ |\ 0\le t\le 1\}$ 
be a geodesic for $\varphi(\xi;v)$ with unit speed, i.e.,
$$
\varphi_{\xi^i}(te_n; e_n)\equiv
\partial_t  \varphi_{v^i}(te_n; e_n),
\qquad \forall\  0\le t\le 1,  1\le i\le n,
$$
and
\begin{equation}    \varphi(te_n; e_n)\equiv 1,
\qquad  0\le t\le 1.
\label{g5}
\end{equation} 
 
    Then, in an open neighborhood of
the geodesic segment,  
 there exists some
non-singular change of variables $\xi=\xi(\eta)$ such that
\begin{equation}
\xi_\eta(0)=Id, \quad
 \xi(te_n)=te_n,\quad \xi_\eta(te_n)e_n=e_n\qquad  0\le t\le 1,
\label{g1}
\end{equation}  
and $\widetilde \varphi(\eta;w)=\varphi(\xi(\eta); \xi_{\eta} w))$
satisfies (\ref{g5}) and 
\begin{equation}  
\widetilde \varphi_{\eta^j}(te_n; e_n)=0, \qquad  1\le j\le n,\ 
 0\le t\le 1,
\label{g2}
\end{equation}
\begin{equation}
\widetilde \varphi_{w^\alpha}(te_n; e_n)=0, \qquad  1\le \alpha\le n-1,\
 0\le t\le 1,
\label{g6}
\end{equation}
and
\begin{equation}
\widetilde \varphi_{\eta^j w^k}(te_n; e_n)=0, \qquad  1\le j, k\le n,\
 0\le t\le 1.
\label{g3}
\end{equation} 
\label{lem2-1}
\end{lem} 
By the homogeneity, it then follows that
\begin{equation}
\widetilde \varphi_{w^n}(te_n; e_n)\equiv 1.
\label{3.6}
\end{equation}

The reader may choose to postpone reading
the long proof of the lemma and go on to the next section.

\noindent {\bf Proof.}\ 
By chain rule,
$$
\widetilde \varphi_{\eta^j}=
 \varphi_{\xi^i}\xi^i_j+  \varphi_{v^i}\xi^i_{lj}w^l,
$$
where we have used notations: $\xi^i_j:=
\frac {\partial \xi^i}{\partial\eta^j}$ and
$ \xi^i_{lj}:=  \frac {\partial^2  \xi^i}{\partial\eta^l
\partial\eta^j}$.

\noindent (i)\ Let
$$
b_\beta(t):=-\int_0^t \varphi_{\xi^\beta}(\tau e_n; e_n)d\tau.
$$
We take
$$
\xi=\xi(\eta):=
(\eta^1, \cdots, \eta^{n-1}, \eta^n+\sum_{\beta=1}^{n-1}
b_\beta(\eta^n)\eta^\beta).
$$

It is easy to check that
$$
\xi^\alpha_\beta(te_n)\equiv \delta^\alpha_\beta, \quad \xi^\alpha_n
(te_n)\equiv 0,
$$
$$
\xi^n_\beta(te_n)\equiv b_ \beta(t), \quad \xi^n_n(te_n)\equiv 1,
$$
$$
\xi^\alpha_{\beta\gamma}(te_n)\equiv  \xi^\alpha_{\beta n}
(te_n)\equiv \xi^\alpha_{n\beta}(te_n)\equiv 
 \xi^\alpha_{nn}(te_n)\equiv 0,
$$
$$
\xi^n_{\beta\gamma}(te_n)\equiv \xi^n_{nn}(te_n) 
\equiv 0, \quad \xi^n_{\beta n}
(te_n)\equiv \xi^n_{n\beta}(te_n)\equiv  b_\beta'(t).
$$
Identity (\ref{g1}) follows from the above.
Also, from the above,
$$
\det\left( \xi^i_j(te_n)\right)\equiv 1.
$$
Thus the change of variables is non-singular
near $\{t e_n\ |\ 0\le t\le 1\}$.

For  $1\le \beta\le n-1$, 
\begin{eqnarray*}
\widetilde \varphi_{\eta^\beta}(te_n; e_n)&=& 
 \varphi_{\xi^i}(te_n; e_n)\xi^i_\beta(te_n)+ \varphi_{v^i}(te_n; e_n)\xi^i_
{n\beta}(te_n)\\
&=&  \varphi_{\xi^\alpha}\xi^\alpha_\beta 
+\varphi_{\xi^n}\xi^n_\beta+\varphi_{v^n}\xi^n_{n\beta}
= \varphi_{\xi^\beta}+\varphi_{\xi^n}\xi^n_\beta+\varphi_{v^n}\xi^n_{n\beta}.
\end{eqnarray*}

Differentiating  (\ref{g5}) in $t$, we find 
\begin{equation}
 \varphi_{ \xi^n}(te_n; e_n)\equiv 0.
\label{c1}
\end{equation}
By (\ref{g5}) and the homogeneity of $\varphi$
in $v$, 
\begin{equation}
 \varphi_{ v^n}(te_n; e_n)\equiv  \varphi(te_n; e_n)\equiv 1.
\label{c2}
\end{equation}
Using (\ref{c1}) and (\ref{c2}), we have
$$
\widetilde \varphi_{\eta^\beta}(te_n; e_n)=
\varphi_{\xi^\beta}(te_n; e_n)+\xi^n_{n\beta}(te_n)
=\varphi_{\xi^\beta}(te_n; e_n)+b_\beta'(t)=0.
$$

Next,
by (\ref{c1}), 
$$
\widetilde \varphi_{\eta^n}(te_n; e_n)=
\varphi_{\xi^i}(te_n; e_n)\xi^i_n(te_n)+
\varphi_{v^i}\xi^i_{nn }(te_n)=
\varphi_{\xi^n}(te_n; e_n)=0.
$$

We have verified (\ref{g2}).

\medskip

\noindent (ii)\ 
Since we have verified (\ref{g2}) for $\widetilde \varphi$ and
the change of variables also preserve the hypotheses on $\varphi$,
we may assume without loss of generality that, to start, the $\varphi$
satisfies the additional hypothesis
\begin{equation}
\varphi_{\xi^j}(te_n; e_n)=0, \qquad  1\le j\le n,\
 0\le t\le 1.
\label{g2new}
\end{equation}

Now we try to make a change of variables such that
$\widetilde \varphi$
also satisfies (\ref{g1}), (\ref{g2}) and, in addition,
 (\ref{g6}).  Later we do another transformation to 
ensure also
(\ref{g3}).

Since $\{te_n\}$ is  a geodesic, we deduce from the
geodesic equations together with (\ref{g2new}) that
\begin{equation}
\varphi_{v^i}(te_n;e_n)\equiv \varphi_{v^i}(0;e_n),
\qquad \forall\ 1\le i\le n.
\label{c3}
\end{equation}

Let
$$
A=
\left(
\begin{array}{ccccc}
1&0 & \cdots & 0 &0 \\
0&1& \cdots & 0&0\\
\cdots	&\cdots&\cdots &\cdots &\cdots \\
0& 0	& \cdots& 1 &0 \\
-\varphi_{v^1}(0;e_n)	&-\varphi_{v^2}(0;e_n)  &\cdots  &
-\varphi_{v^{n-1}}(0;e_n)  &1
\end{array}
\right),
$$
and consider a linear change of variables
$$
\xi=\xi(\eta):=A\eta.
$$
Let
$$
\widetilde \varphi(\eta;w)=\varphi(\xi(\eta); \xi_{\eta} w))
=\varphi(A\eta; Aw).
$$

Clearly the change of variables satisfies (\ref{g1}).
By (\ref{c2}) and (\ref{c3}),
we have
\begin{eqnarray*}
\widetilde \varphi_{w^\alpha}(te_n;e_n)&=&
\varphi_{v^i}(te_n;e_n)A^i_\alpha=
\varphi_{v^\alpha}(te_n;e_n)+\varphi_{v^n}(te_n;e_n)
A^n_\alpha\\
&=& \varphi_{v^\alpha}(te_n;e_n)+
A^n_\alpha=\varphi_{v^\alpha}(te_n;e_n)-
\varphi_{v^\alpha}(0;e_n)=0.
\end{eqnarray*}
We have verified that $\widetilde \varphi$ satisfies
(\ref{g6}).  Clearly  $\widetilde \varphi$ satisfies 
(\ref{g2}),  since $\varphi$ satisfies 
(\ref{g2new}).

So from now on, we may assume without loss of generality that
$\varphi$ further satisfies (\ref{g2new}) and
\begin{equation}
\varphi_{v^\alpha}(te_n; e_n)=0, \qquad  1\le \alpha\le n-1,\
 0\le t\le 1,
\label{g6new}
\end{equation}

\noindent (iii)\
Let $\psi:=\varphi^2$.  For $ 1\le \alpha, \beta\le n-1$, we have,
 by (\ref{g5}) and  (\ref{g6new}), 
$$
\psi_{v^\alpha v^\beta}(te_n; e_n)=2\varphi(te_n, e_n)
\varphi_{v^\alpha v^\beta}(te_n; e_n)=2\varphi_{v^\alpha v^\beta}(te_n; e_n).
$$
So, by the positivity of 
$(\psi_{v^\alpha v^\beta})$,
 $A:=\left(\varphi_{v^\alpha v^\beta}(te_n; e_n)\right)$
is real symmetric and positive definite.

Let
$$
E(t):=\left(\varphi_{\xi^\alpha v^\beta}(te_n; e_n)\right).
$$
By
Lemma \ref{lema-2} in Appendix B, 
 the dimension of the space of  solutions of
$$
X^TA-AX=E^T-E
$$
is $\frac {(n-1)n}2$. 
 For fixed $t$, let $X(t)$ be the solution of of the above equation
with the least Euclidean norm.  Clearly $X(t)$ 
depends smoothly on $t$. 

Let $B(t)$ be the solution of
$$
\left\{
\begin{array}{rll}
\dot B(t):&=& \frac {d}{dt} B(t)=XB, \qquad 0\le t\le 1,\\
B(0)&=& I,
\end{array}
\right.
$$
---- clearly $\det(B(t))\neq 0, 0\le t\le 1$---- 
and let
$$
M(t):= B^T E^TB +B^TA\dot B.
$$
It is easy to see that $M$ is symmetric, i.e.
$$
M^T\equiv M.
$$

We introduce a final change of variables $\xi=\xi(\eta)$ by
$$
\left\{
\begin{array}{rll}
\xi^\alpha &=& \sum_{  1\le \beta\le n-1}
B^\alpha _\beta(\eta^n)\eta^\beta,
\\
\xi^n  &=&  \eta^n - \frac 12
\sum_{  1\le \gamma, \mu\le n-1} M_{\gamma\mu} (\eta^n)  \eta^\gamma
\eta^\mu.
\end{array}
\right.
$$
Then (\ref{g1}) holds,
$$
\xi_\eta(0)=ID, \qquad 
\xi(te_n)=te_n,
\qquad \xi_\eta(te_n) e_n=e_n,
$$
and
$$
\det( \xi_\eta(te_n) )= \det( B(t) )\ne 0,
$$
$$
\xi^\alpha_n(te_n)\equiv \xi^\alpha_{nn}(te_n)\equiv 
\xi^n_{jn}(te_n)\equiv \xi^n_{nj}(te_n)\equiv
\xi^n_\alpha (te_n)\equiv 0.
$$

Let $\widetilde \varphi(\eta;w)=\varphi(\xi(\eta); \xi_{\eta} w))$.
Using  (\ref{g2new}), (\ref{g6new}), and the
 above listed properties of the change of variables, we find 
$$
\widetilde \varphi_{\eta^j}(te_n;e_n)= \varphi_{v^l}\xi^l_{nj}=\varphi_{v^n}
\xi^n_{nj}=0, \qquad 1\le j\le n,
$$
$$
\widetilde \varphi_{w^\alpha}(te_n;e_n)=
\varphi_{v^i}\xi^i_\alpha= \varphi_{v^n}
\xi^n_\alpha =0,  \qquad 1\le \alpha\le n-1.
$$
We have verified that (\ref{g1}), (\ref{g2}) and (\ref{g6})
continue to hold in the new variables.

\noindent (iv)\ Finally, to verify (\ref{g3}), consider, at $(te_n,e_n)$, 
$$
\widetilde \varphi_{w^i\eta^j}=\varphi_{\xi^l v^m}\xi^m_i \xi^l_j +\xi ^n_{ij}+
 \varphi_{v^lv^m} \xi^m_i\xi^l_{nj}.
$$
By (\ref{g6new}), 
\begin{equation}
\varphi_{\xi^nv^\alpha}(te_n;e_n)\equiv 0,  \qquad 1\le \alpha\le n-1,
\label{1new}
\end{equation}
and, using also the homogeneity of $\varphi$ in $v$, 
\begin{equation}
\varphi_{v^\alpha v^n}(te_n;e_n)\equiv 0,  \qquad 1\le \alpha\le n-1.
\label{2new}
\end{equation}
By (\ref{c2}) and the homogeneity of $\varphi$ in $v$,
\begin{equation}
\varphi_{v^nv^n}(te_n;e_n)\equiv 0 .
\label{3new}
\end{equation}
By (\ref{g2new}) and  the homogeneity of $\varphi$ in $v$,
\begin{equation}
\varphi_{\xi^jv^n}(te_n;e_n)\equiv 0,  \qquad 1\le j\le n.
\label{4new}
\end{equation}

Simplifying the expression of $\widetilde \varphi_{w^i\eta^j}$ by using
(\ref{4new}), (\ref{2new}), (\ref{3new}) and (\ref{1new}), we have
\begin{eqnarray*}
\widetilde \varphi_{w^i\eta^j}(te_n;e_n)&=&
\varphi_{\xi^l v^\alpha} \xi^\alpha _i \xi^l_j + \xi^n_{ij}
+  \varphi_{v^\alpha v^\beta}  \xi^\beta_i  \xi^\alpha_{nj}\\
&=&  \varphi_{\xi^\beta v^\alpha} \xi^\alpha _i \xi^\beta_j +\xi^n_{ij}
+  \varphi_{v^\alpha v^\beta}  \xi^\beta_i  \xi^\alpha_{nj}.
\end{eqnarray*}

Since $\xi^\beta_n(te_n)\equiv \xi^\alpha_{nn}(te_n) \equiv
 \xi^n_{in}(te_n) \equiv 0$ for all
$1\le \alpha, \beta\le n-1$ and $1\le i\le n$, 
we have
$$
\widetilde \varphi_{w^i \eta^n}\equiv 0,
\qquad 1\le i \le n.
$$
Similarly,
$$
\widetilde \varphi_{w^n \eta^j }\equiv 0,
\qquad 1\le j \le n.
$$
Finally, for $1\le \gamma, \mu\le n-1$, 
as one may check,
\begin{eqnarray*}
\widetilde \varphi_{w^\gamma \eta^\mu}(te_n;e_n)&=&
\varphi_{\xi^\beta v^\alpha} \xi^\alpha_\gamma \xi^\beta_\mu+
\xi^n_{\gamma\mu} +
\varphi_{v^\alpha v^\beta}  \xi^\beta_\gamma \xi^\alpha_{n\mu}\\
&=& M_{\gamma\mu}+ \xi^n_{\gamma\mu}=0.
\end{eqnarray*}
In the above, we have used
$$
\dot B^\alpha_\mu(t)= \frac{d}{dt} B^\alpha_\mu(t)=  \frac{d}{dt} 
\xi^\alpha_\mu(te_n)=\xi^\alpha_{n\mu}(te_n).
$$
We have thus verified 
(\ref{g3}).  Lemma \ref{lem2-1} is established.

\vskip 5pt
\hfill $\Box$
\vskip 5pt

\section{}
\setcounter{equation}{0}
In this section we establish some properties 
of the cut points and conjugate points of $y$ on $\partial \Omega$.
In particular we first prove the continuity of the
map $m(y)$, defined on $\partial \Omega$, and
then  prove that $m(y)=\tilde m(y)$ for all
$y\in \partial \Omega$ and, consequently, $\Sigma=
\tilde \Sigma$.  

\noindent{\bf 4.1.}\ 
For $y\in \partial \Omega$,
without  loss of generality,
we may assume $\bar s(y)=\bar s(0)=1$. 
  Then we use
our special coordinates of Section 3; near
the origin  $\Omega$ is given by
$x_n>f(x')$ with
$$
f(0')=0, \qquad \nabla f(0')=0.
$$
Then $m(y)=m(0)=e_n$.  The ``normal'' geodesic
from $0$ lies along the $x_n-$axis.

For $\epsilon_0>0$, let $\Gamma
:=\{t e_n\ |\ -\epsilon_0\le t\le 1+\epsilon_0\}$
be the geodesic for $\varphi(\xi;v)$   satisfying,
for $-\epsilon_0\le t\le  1+\epsilon_0$, the conclusions
of Lemma \ref{lem2-1} and (\ref{3.6}):
\begin{equation}    \varphi(te_n; e_n)\equiv 1,
\label{ag5}
\end{equation}
\begin{equation}
 \varphi_{\xi^j}(te_n; e_n)=0, \qquad  1\le j\le n,
\label{g2n}
\end{equation}
\begin{equation}
\varphi_{v^\alpha}(te_n; e_n)=0, \qquad  1\le \alpha\le n-1,
\label{g6n}
\end{equation}
and
\begin{equation}
\varphi_{\xi^j v^k}(te_n; e_n)=0, \qquad  1\le j, k\le n.
\label{g3n}
\end{equation}
By (\ref{g6n}) and the homogeneity of $\varphi$ in $v$, we have
\begin{equation}
\varphi_{v^\alpha v^n}(te_n; e_n)\equiv 0,
\qquad 1\le \alpha\le n-1,  -\epsilon_0\le t\le 1+\epsilon_0.
\label{sss1}
\end{equation}
Differentiating (\ref{g2n}), we have
\begin{equation}
\varphi_{\xi^j\xi^n}(te_n; e_n)\equiv 0,
\qquad 1\le j\le n, -\epsilon_0\le t\le 1+\epsilon_0.
\label{aaa2}
\end{equation}

For $y\in \partial \Omega$, let $\xi=\xi(y,\tau)$ denote the
geodesic satisfying
$$
\varphi(\xi; \dot\xi)\equiv 1,
$$
$$
\xi(y,0)=y,
$$
and
$$
\dot\xi(y,0)=V(y),
$$
where $V(y)$ is as in (\ref{1.14}).

Recall that for $|x'|<\epsilon_1$, we  write $\xi((x',f(x')), \tau)$ as
  $\xi(x', \tau)$, i.e. $\xi=\xi(x', \tau)$ is
the geodesic satisfying
$$
\varphi(\xi; \dot\xi)\equiv 1,
$$
$$
\xi(x',0)=(x', f(x')),
$$
and
$$
\dot \xi(x',0)=V(x'),
$$
where $V(x')$ is the vector-valued function defined in
Section 2.

\bigskip

\centerline{
\epsfbox{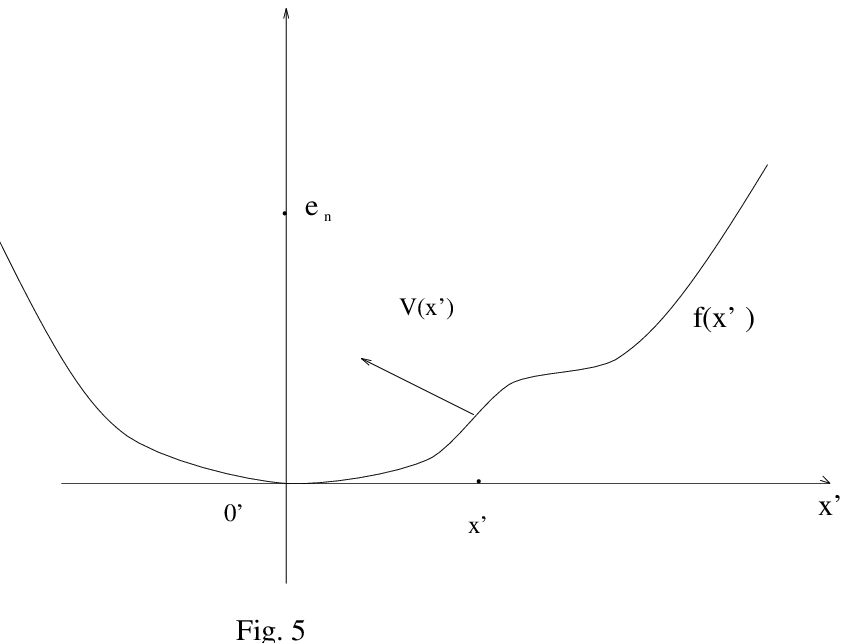}}

\bigskip

The following lemma  establishes the continuity of the map $m(y)$.

\begin{lem} Suppose, as above, 
 $m(0)=e_n$. Then
$\displaystyle{ \lim_{|x'|\to 0} m((x', f(x'))=e_n }$, i.e.,
$m$ is continuous at $0$.
\label{lem4.7}
\end{lem}

\noindent{\bf Proof.}\
We prove it by contradiction argument.
Suppose the contrary, there exist $x'_{ i }\to 0$
such that $m(( x'_{i}, f(x'_{ i }))=\xi(x'_{ i }, t_i)$
with $t_i\to \bar t\ne 1$.
We know that $\xi(x'_{ i }, t_i)\to \xi(0', \bar t)\in \Sigma$,
so we must have $\bar t\ge 1$.
On the other hand, if
 $\bar t>1$, then, by compactness,
there exists some $\delta>0$, independent of $i$,
such that the $\delta$-neighborhood of
$\{\xi(x'_{ i }, t)\ |\ 0\le t\le \frac {1+\bar t}2\}$
belongs to $G$, the complement of $\Sigma$, for large $i$.
Since $ \frac {1+\bar t}2>1$,
this  set would contain $e_n$ for large $i$, a contradiction.
Lemma \ref{lem4.7} is established.

\vskip 5pt
\hfill $\Box$
\vskip 5pt

\noindent{\bf 4.2.}\
We will prove that $m(0)= \tilde m(0)$.
We first show

\begin{lem} Suppose $m(0)=e_n$.  Then
$\tilde m(0)=\tilde t e_n$ for some $\tilde t\ge 1$.
\label{lem4-new1}
\end{lem}

\noindent{\bf Proof.}\  We argue by contradiction.
Suppose $0<\tilde t<1$, then, since
$\tilde t e_n\in G$ and $G$ is open, $G$ contains
a neighborhood of $\tilde te_n$.  For $X$
close to $\tilde te_n$, $X$ in $G$, there exists
a unique $z=z(X)\in \partial \Omega$
such that
$$
dist(z\add X)=dist(\partial \Omega\add X).
$$
Clearly, the map $X$ to $z(X)$ is continuous near
$\tilde te_n$.
Since $dist(0\add \tilde te_n)=
dist(\partial \Omega\add \tilde te_n)$ and
since $\tilde te_n\in G$, we find
$z(\tilde te_n)=0$ and, by the continuity 
of the map, $z(X)$ is close to $0$ for
$X$ close to $\tilde te_n$.  So we can write
$$
z(X)=(x'(X), f(x'(X))), 
$$
where $x'(X)$ is continuous near $\tilde te_n$
with $x'(\tilde te_n)=0'$.

For $X$ close to $\tilde te_n$, consider a geodesic,
with unit speed, joining 
$z(X)=(x'(X), f(x'(X)))$ to $X$
which realizes 
$$
\ell(X):=dist(z(X)\add X)=dist(\partial \Omega\add X).
$$
By Lemma \ref{lem2.2} and the fact that the geodesic
must enter $\Omega$ (otherwise it would not
realizes the distance of $\partial \Omega$ to
$X$ since it has to enter $\Omega$),  the geodesic is
$\xi(x'(X), s)$, and
\begin{equation}
X=\xi(x'(X), \ell(X)).
\label{B3-1}
\end{equation}

It is easy to see that $\ell(X)$ is a continuous function near
$\tilde te_n$.  Consider the following map defined
in a neighborhood of $\tilde te_n$:
$$
F(X):=(x'(X), \ell(X)).
$$
One verifies that $F$ is one-to-one near $\tilde te_n$.
A continuous one-to-one map is open, i.e.
it maps open sets to open sets.  So $F$ maps a neighborhood of
$\tilde te_n$ to a neighborhood of 
$(0', \tilde t)$.  For $t$ close to
$\tilde t$, let
$X_t=F^{-1}(0', t)$. Then
$(0', t)=F(X_t)$ and therefore
$x'(X_t)=0'$, $\ell(X_t)=t$.
By (\ref{B3-1}), $X_t=\xi(0', t)=te_n$, i.e.
$$
t=\ell(te_n)=dist(\partial\Omega\add te_n)\ \ 
\mbox{for}\ t\ \mbox{close to}\
\tilde t,
$$
violating $\tilde te_n=\tilde m(0)$.
Lemma \ref{lem4-new1} is established.

\vskip 5pt
\hfill $\Box$
\vskip 5pt

Sometimes, for convenience, we normalize 
so that $\tilde m(0)=e_n$ instead of $m(0)=e_n$.
We still have the same properties of our special
coordinates stated at the beginning of this section.

\bigskip

\noindent{\bf 4.3.}\
\begin{lem} Assume $\tilde m(0)=e_n$.  Then there 
exists some  $\mu>0$, and for all
$H\in C^1_0([0,1], \Bbb R^{n-1})$,
$$
\int_0^1 \left( \varphi_{\xi^\beta\xi^\gamma}(te_n;e_n)
H^\beta H^\gamma+ \varphi_{v^\beta v^\gamma}(te_n;e_n)
\dot H^\beta \dot H^\gamma\right)dt
\ge \mu\int_0^1 H^2dt.
$$
\label{lem-posi}
\end{lem}

An easy consequence is
\begin{cor} Under the same hypotheses of Lemma \ref{lem-posi},
there exists $\mu_1>0$ such that
 for all
$H\in C^1_0([0,1], \Bbb R^{n-1})$,
$$
\int_0^1 \left( \varphi_{\xi^\beta\xi^\gamma}(te_n;e_n)
H^\beta H^\gamma+ \varphi_{v^\beta v^\gamma}(te_n;e_n)
\dot H^\beta \dot H^\gamma\right)dt
\ge \mu_1\int_0^1 \dot H^2dt.
$$
\label{cor-posi}
\end{cor}

\begin{rem}
One sees from the proof that the conclusion of
Lemma \ref{lem-posi} and Corollary \ref{cor-posi}
holds when replacing $e_n=\tilde m(0)$ by $\hat t e_n$
for any $0<\hat t<1$.
\label{rem-posi}
\end{rem}

\noindent {\bf Proof of Lemma \ref{lem-posi}.}\
Let $\mu$ be the first eigenvalue of the quadratic form, i.e.,
$\mu$ is the largest number such that for all
$H\in C^1_0([0,1], \Bbb R^{n-1})$ we have
$$
\int_0^1 \left( \varphi_{\xi^\beta\xi^\gamma}(te_n;e_n)
H^\beta H^\gamma+ \varphi_{v^\beta v^\gamma}(te_n;e_n)
\dot H^\beta \dot H^\gamma\right)dt
\ge \mu\int_0^1 H^2dt.
$$
We only need to show that $\mu>0$.  If not, then
for $\epsilon>0$, there exists $\overline H\in C^1_0([-\epsilon,1],
 \Bbb R^{n-1})$
such that
$$
\int_{-\epsilon}^1 \left( \varphi_{\xi^\beta\xi^\gamma}(te_n;e_n)
\overline H^\beta\overline  H^\gamma+ \varphi_{v^\beta v^\gamma}(te_n;e_n)
\dot {\overline H}^\beta \dot {\overline H}^\gamma\right)dt< 0.
$$
We identify $\overline H(t)$ with $(\overline H(t), 0)$
in $\Bbb R^{n+1}$, and
perturb the geodesic
$te_n$ by considering
 $\zeta(\tau,t)=te_n+\tau \overline   H(t)$, $-\epsilon<t\le 1$.
Then at $\tau =0$, we have
$$
\frac {d}{d\tau}\int_{-\epsilon}^1
\varphi(\zeta; \dot \zeta)dt=0,
$$
and
$$
\frac {d^2}{d\tau ^2}\int_{-\epsilon}^1
\varphi(\zeta; \dot \zeta)dt<0.
$$
It follows that for $\tau>0$ small, we have
\begin{equation}
\int_{-\epsilon}^1
\varphi(\zeta; \dot \zeta)dt< 1+\epsilon.
\label{contra}
\end{equation}
On the other hand, let $\bar t=\bar t(\tau)>0$ be such that
$$
\zeta(\tau, \bar t)=(x', f(x'))
$$
for some $x'$.  Since $\tilde m(0)=e_n$, we find
$$
\int_{\bar t}^1
\varphi(\zeta; \dot \zeta)dt\ge 1,
$$
and, by Lemma \ref{lem2.3},
$$
\int_{-\epsilon}^1
\varphi(\zeta; \dot \zeta)dt\ge \epsilon
$$
for $\epsilon$ sufficiently small.
The above two estimates violate (\ref{contra}), a contradiction.

\vskip 5pt
\hfill $\Box$
\vskip 5pt

 \noindent{\bf 4.4.}\
 We still assume that $e_n=\tilde m(0)$,
and we  now consider geodesics ending at $e_n$.
For $\sigma'=(\sigma_1,\cdots, \sigma_{n-1})\in \Bbb R^{n-1}$ satisfying
$|\sigma'|\le 1/2$, let
$\tau=\tau(\sigma')$ be defined by
$$
\varphi\left(e_n; (\sigma', \tau)\right)=1,
$$
and
$$
\tau(0')=1.
$$
Since $\varphi_{v^n}(e_n;e_n)=1$, by
the Implicit Function Theorem, $\tau$ exists as
 a smooth function of 
$\sigma'$.

Let $\eta=\eta(\sigma', t)$ be the unique smooth solution of
$$
\psi_{\xi^i}(\eta; \dot \eta)= \frac {d}{dt}\psi_{v^i}
(\eta; \dot \eta), \qquad t\le 1,
$$
satisfying
$$
\eta(\sigma', 1) =e_n,
$$
$$
\dot \eta(\sigma', 1)= (\sigma',
\tau(\sigma')).
$$
The solution exists for all time until 
it hits the boundary $(x', f(x'))$ (in fact it
goes further since $\varphi$ has been extended to a fixed 
open neighborhood of the domain).

Clearly  $\eta(\sigma', t)$ is a geodesic and
(see Fact 5)
$$
\psi(\eta(\sigma', t); \dot \eta(\sigma', t))
\equiv \psi(\eta(\sigma', 1); \dot \eta(\sigma', 1))
= \psi(e_n; (\sigma', \tau(\sigma'))=1.
$$
Applying $\frac{\partial}{\partial \sigma_\alpha}$ to the
geodesic equations and setting $\sigma'=0$, we have,
by our special coordinates, 
$$
\varphi_{\xi^\beta\xi^\gamma}(te_n;e_n)
\eta^\gamma_{\sigma_\alpha}(0', t)
\equiv\frac {d}{dt} \left( \varphi_{v^\beta v^\gamma}(te_n;e_n)
\dot \eta^\gamma_{\sigma_\alpha}(0', t)\right),
\qquad 0\le t\le 1.
$$
We remark that  $(\varphi_{v^\beta v^\gamma}(te_n;e_n))=\frac 12
(\psi_{v^\beta v^\gamma}(te_n;e_n))$
is positive definite and 
$$
\eta^\gamma(0', 1)=0, \qquad
\dot \eta^\gamma_{\sigma_\alpha}(0',1)=\delta^\gamma_\alpha.
$$
With the aid of  Lemma \ref{lem-posi},
one sees that $\{\eta_{\sigma_1}(0',0),
\cdots, \eta_{\sigma_{n-1}}(0',0)\}$ 
are linearly independent.

\medskip

By compactness, for some positive number $\delta>0$, depending only
on $f$ and $\varphi$, we have
\begin{equation}
\det (\eta_{\sigma_1}(0',0),
\cdots, \eta_{\sigma_{n-1}}(0',0))\ge c>0.
\label{independent}
\end{equation}

Let
$$
x_\alpha=\eta^\alpha(\sigma_1, \cdots, \sigma_{n-1}, 0),
\qquad 1\le \alpha\le n-1.
$$
We know from the above, using the Implicit Function Theorem,
that the map $\sigma'$ to $x'$ is a diffeomorphism
in a fixed neighborhood of $0'$ (the size of the neighborhood
depends only on $f$ and $\varphi$).

Define
$$
\tilde f(x_1,\cdots, x_{n-1})=\eta^n(\sigma_1, \cdots, \sigma_{n-1}, 0).
$$
Then for some positive constants
$\tilde \epsilon_10$ and $C$, depending only on
$f$ and $\varphi$, we have
\begin{equation}
\| \tilde f\|_{ C^{2,1}(B_{{\tilde\epsilon}_1}) }\le
C.
\label{new1}
\end{equation}

In fact, the parameter sphere we have constructed is
a distance sphere near the origin, i.e. 
for possibly a smaller positive constant $\tilde \epsilon_1$,
still depending only on
$f$ and $\varphi$, we have
\begin{equation}
dist((x', \tilde f(x'))\add e_n)=1, \qquad |x'|<\tilde\epsilon_1.
\label{new2}
\end{equation}
Indeed, if the above does not hold for any $\tilde \epsilon_1$, 
then there exist $x_i'\to 0$ such that
$$
b_i:=dist ((x'_i, \tilde f(x'_i))\add e_n)<1.
$$
It may appear that the above statement is negating
(\ref{new2}) for $\tilde\epsilon_1$ which depends on
the initial base point we pick (the origin), but this can
be taken care by an easy compactness argument.

Let $\zeta_i$ be shortest geodesics, with unit speed, joining
$(x'_i, \tilde f(x'_i))$ to $e_n$.   
We know that
$e_n=\zeta_i(b_i)$.
After passing to a subsequence,
$b_i\to b\le 1$, $\zeta_i\to \zeta$ in $C^1$ norm.
Clearly $\zeta$ is a geodesic with unit
speed, $\zeta(0)=0$, $\zeta(b)=e_n$.
Since
$$
dist(0\add e_n)=1,
$$
we have $b\ge dist(0\add e_n)=1$.  Since we also know $b\le 1$, 
we find
$b=1$.  
Now we know that $dist(0\add \zeta(1))=dist(\partial\Omega\add \zeta(1))$,
we find, by Lemma \ref{lem2.2}, 
 $\zeta(t)$ is normal to $\partial\Omega$ 
at the origin.
Since  $\zeta$ must enter $\Omega$
(otherwise it would not realize the distance of $\partial \Omega$ to
$\zeta(1)$),  $\zeta(t)\equiv \xi(0', t)\equiv te_n$.
Thus $\dot \zeta_i(b_i)$ is, for large $i$,
 close to $e_n$, and the geodesics $\zeta_i$ comes from 
the spreading geodesics from $e_n$ we have constructed, i.e.,
for some $\sigma'_i\to 0'$,
$$
\zeta_i(t)\equiv \eta(\sigma'_i,   t+1-b_i).
$$
On the other hand, we know that $\zeta_i(0)$ is on the
graph of $\tilde f$, so
$\eta(\sigma'_i,   1-b_i)$ is on the graph
of $\tilde f$.  It follows that $b_i=1$, a contradiction.
(\ref{new2}) is established.

Summarizing the above, we have 
established the following
\begin{lem}  Under the hypotheses
stated at the beginning of Section 4, though assuming
$\tilde m(0)=e_n$ instead of $m(0)=e_n$, there
exists a smooth function $\widetilde f$ satisfying
(\ref{new1}) and (\ref{new2}) for some 
positive constants
 $\tilde\epsilon_1$ and $C$ depending only on $f$ and $\varphi$.
\label{lem4-4}
\end{lem}

\begin{rem} The distance sphere centered at $\tilde m(0)=e_n$ 
can be constructed the same way with center to be any point
before $\tilde m(0)$, i.e. with center $\hat t e_n$ for
any $0<\hat t<1$, though in this case, 
the $\tilde\epsilon_1$ depends also on 
the positive lower bound of $\hat t$.
\label{rem-new1}
\end{rem}

\bigskip

\centerline{
\epsfbox{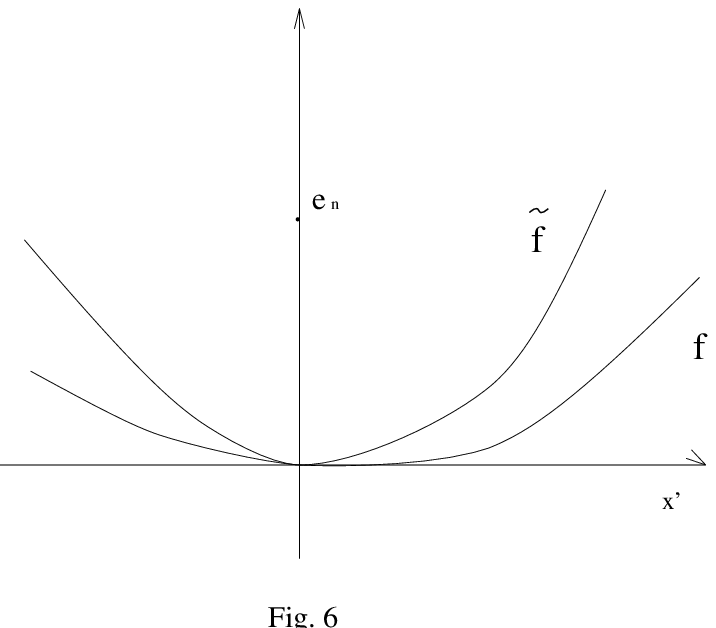}}

\bigskip

\begin{rem}  Clearly, under the assumption
of Lemma \ref{lem4-4},
$$
\widetilde f(x')-f(x')\ge 0, \qquad |x'|<\tilde\epsilon_1,
$$
$$
\widetilde f(0')=0,\qquad
\widetilde f_{x_\alpha  }(0')=0,\qquad 1\le \alpha\le n-1.
$$
\label{rem2-20}
\end{rem}

Let $\lambda$ denote the smallest eigenvalue of
$(\widetilde f_{x_\alpha   x_\beta  }(0')-f_{x_\alpha  x_\beta  }(0'))$;
we know that
$\lambda\ge 0$.

We may carry out the above for points $X$ near $e_n$ instead of
for $e_n$ only.  Indeed for $X$ close to $e_n$  and for small
$\sigma'$, let $\tau=\tau(\sigma', X)$  be defined by
$$
\varphi\left(X; (\sigma', \tau)\right)=1,
$$
and
$$
\tau(0',0)=1.
$$
$\tau$ is a smooth function of $\sigma'$ and $X$.

Let $\eta=\eta(\sigma', X, t)$ be the unique smooth solution of
$$
\psi_{\xi^i}(\eta; \dot \eta)= \frac {d}{dt}\psi_{v^i}
(\eta; \dot \eta), \qquad t\le 1,
$$
satisfying
$$
\eta(\sigma', X,1)=X,
$$
$$
\dot \eta(\sigma', X, 1)=(\sigma', \tau(\sigma', X)).
$$

Because of  (\ref{independent}), there exists some positive constant
$\epsilon$  such that
for every $|t|\le \epsilon$ and $|X-e_n|\le \epsilon$, 
$\{\eta(\cdot, X, t)\}$ is locally represented as
a graph, and the  
gradient and  Hessian of the function
representing the graph
converges to those of $\tilde f$ as $\epsilon$ tends to $0$.

Let us still assume that $\tilde m(0)=e_n$. Then
 $\tilde f$ is defined by
Lemma \ref{lem4-4}, with the nonnegative
least eigenvalue $\lambda$ of
$(\tilde f_{x_\alpha   x_\beta  }(0')-f_{x_\alpha  x_\beta  }(0'))$.
For $0<\epsilon<\frac 12$, an application
of Lemma \ref{lem4-4} together with Remark \ref{rem-new1}
yields a smooth function $\tilde f^{(\epsilon)}$
satisfying, for some constants
 $\delta, C>0$ depending only on $\varphi$ and $f$,
$$
dist((x', \tilde f^{(\epsilon)}(x'))\add (1-\epsilon)e_n)=1-\epsilon,
\qquad |x'|<\delta,
$$
$$
 \tilde f^{(\epsilon)}(0')=0, \quad
\| \tilde f^{(\epsilon})\|_{ C^{2,1}(B_\delta)}\le
C,
$$
and, by the triangle inequality for the Finsler metric,
$$
\tilde f^{ (\epsilon_2)}(x')\ge 
 \tilde f^{(\epsilon_1)}(x')\ge  \tilde f(x'), \qquad 
\forall\ |x'|<\delta, \ 0<\epsilon_1<\epsilon_2<\frac 12.
$$
Consequently,
$$
 \tilde f^{(\epsilon)}_{x_\alpha}(0')=0, \qquad
1\le \alpha\le n-1.
$$

Let $\lambda^{(\epsilon)}$ denote the least eigenvalue
of
$(\tilde f^{(\epsilon)}_{x_\alpha   x_\beta  }(0')
-f_{x_\alpha  x_\beta  }(0'))$,
and let $\gamma^{(\epsilon)}\ge 0$ be  the least eigenvalue
of
$(\tilde f^{(\epsilon)}_{x_\alpha   x_\beta  }(0')-
\tilde f_{x_\alpha  x_\beta  }(0'))$.
Clearly,
$$
 \lambda^{(\epsilon)}\ge \lambda+
\gamma^{(\epsilon)}.
$$

\begin{lem} Assuming $\tilde m(0)=e_n$. 
 For $0<\epsilon<\frac 12$, let $\gamma^{(\epsilon)}$,
$\lambda^{(\epsilon)}$ and $\lambda$ be as above.  Then
for some constant $c>0$, depending only on
$f$ and $\varphi$, such that
$$
\lambda^{(\epsilon)}-\lambda
\ge \gamma^{(\epsilon)}
\ge c\epsilon.
$$
\label{lem4.6}
\end{lem}

\noindent{\bf Proof.}\
Let, as usual,  $\tilde \xi(x', t)$
denote the geodesics,  with unit speed,
starting from $(x', \tilde f(x'))$ and ``normal''
to the graph of $\tilde f$.  By the property
of $\tilde f$, $\tilde \xi(x', 1)=e_n$.
Similarly, let $\tilde \xi^{(\epsilon)}(x', t)$
denote the geodesics for $\tilde f^{(\epsilon)}$ instead of for $\tilde f$.
 Let $\zeta^{(\epsilon)}$ be a unit eigenvector
of $ \left( \tilde f^{(\epsilon)}_{x_\alpha x_\beta}(0')-
\tilde f_{x_\alpha x_\beta}(0')\right) $ associated with
the least eigenvalue $\gamma^{(\epsilon)}$, and let 
 $x'$  be a multiple of $\zeta^{(\epsilon)}$,
 we find
\begin{equation}
|\tilde \xi(x', t)-\tilde \xi^{(\epsilon)}(x', t)|
\le C(\gamma^{(\epsilon)}|x'|+|x'|^2),\qquad \forall\ 0\le t\le 1.
\label{add1}
\end{equation}
For $t=1-\epsilon$,
$\tilde \xi^{(\epsilon)}(x', 1-\epsilon)=(1-\epsilon)e_n$, and
therefore
\begin{eqnarray}
&&\tilde \xi(x', 1-\epsilon)-\tilde \xi^{(\epsilon)}(x', 1-\epsilon)
\nonumber\\
&=&\tilde \xi(x', 1-\epsilon)-e_n+\epsilon e_n
=  \tilde \xi(x', 1)-
\dot  {\tilde \xi}(x', 1)\epsilon +
O(|x'|\epsilon^2) -e_n+\epsilon e_n\nonumber\\
&=& \epsilon (e_n- \dot {\tilde \xi}(x', 1))+
O(|x'|\epsilon^2).
\label{add2}
\end{eqnarray}
In the above, we have used, as usual, Taylor expansions and
the fact that $\ddot{ \tilde \xi}(0', t)\equiv 0$.

Since  $ \tilde \xi(x', \cdot)$
satisfies the geodesic equations, and since
$ \tilde \xi(0', 1)= \tilde \xi(x', 1)=\dot{\xi}(0', 1)=e_n$, we have,
for some positive constants $a$ and $b$, depending
only on $f$ and $\varphi$, such that
$$
|e_n- \dot{\tilde \xi}(x', 1)|=
|\dot{ \tilde \xi}(0', 1)-\dot{ \tilde \xi}(x', 1)|
\ge b|\tilde \xi(0', 0)-\tilde \xi(x', 0)|\ge a|x'|.
$$

This, together with (\ref{add1}) and (\ref{add2}),
yields
$$
a\epsilon |x'|\le C(|x'|\epsilon^2+\gamma^{(\epsilon)}|x'|+|x'|^2).
$$
Dividing the above by $|x'|$ and sending $|x'|$ to $0$, we find
$$
a\epsilon\le C\epsilon^2+C\gamma^{(\epsilon)}.
$$
The desired estimate follows if
$C\epsilon\le  \frac {a}2$.
If $C\epsilon>  \frac {a}2$, the desired estimate follows
from the estimate for $\epsilon=\frac {a}{2C}$ and the
monotonicity of
  $\gamma^{(\epsilon)} $ in $\epsilon$.

\vskip 5pt
\hfill $\Box$
\vskip 5pt

\noindent{\bf 4.5.}\
To establish $m=\tilde m$, we need, in addition
to Lemma \ref{lem4-new1}, the following

\begin{lem} Assuming $\tilde m(0)=e_n$, then
$te_n\in G$ for all $0<t<1$.
\label{lem4-new2}
\end{lem}

A consequence of Lemma \ref{lem4-new1} and Lemma \ref{lem4-new2} is
\begin{cor} $m(y)=\tilde m(y)$ for all $y\in \partial
\Omega$.  Consequently, $\Sigma =\tilde \Sigma$.
\label{cor4-1}
\end{cor}
 
\noindent{\bf Proof of Lemma \ref{lem4-new2}.}\
We argue by contradiction.  Suppose that
$m(0)=(1-\epsilon) e_n\in \Sigma$ for some
$0<\epsilon<1$.  Clearly $1-\epsilon>\bar \epsilon>0$ for some
$\bar \epsilon$ depending only on $f$ and $\varphi$.
Since $(1-\epsilon)e_n\in \Sigma$, there exist $X_i\to (1-\epsilon)e_n$,
$z_i, \hat z_i\in \partial \Omega$, $z_i\neq \hat z_i$,
such that
$$
b_i:= dist(\partial\Omega\add X_i)=
dist (z_i\add  X_i)
= dist (\hat z_i\add  X_i).
$$
After passing to a subsequence,
we may assume that $z_i\to z$, $\hat z_i\to \hat z$
and $b_i\to b$.
Clearly
$$
b=dist(\partial\Omega\add (1-\epsilon)e_n)
=1-\epsilon,
$$
and 
$$
dist(z\add (1-\epsilon)e_n)=dist(\hat z\add (1-\epsilon)e_n)=
dist(0\add (1-\epsilon)e_n)
=1-\epsilon.
$$
Since $\tilde m(0)=e_n$ and $1-\epsilon<1$, there can only be
one point on $\partial \Omega$ which realizes 
$dist(\partial\Omega\add (1-\epsilon)e_n)$.  So 
we must have $z=\hat z=0$.  Write
$$
z_i=(x'_i, f(x'_i)), \quad \hat z_i=(\hat x'_i, f(\hat x'_i)),
$$
and let $\zeta_i$ and
$\hat \zeta_i$ be  shortest geodesics, with unit speed,
joining respectively $z_i$ and $\hat z_i$ to $X_i$.
By Lemma \ref{lem2.2}, $\zeta_i\equiv \xi(x'_i, \cdot)$
and $\hat \zeta_i\equiv \xi(\hat x'_i, \cdot)$.
So, $\zeta_i\to \xi(0', \cdot)$ and $\hat \zeta_i\to 
\xi(0', \cdot)$ in $C^1$ norm.
It follows that $\dot \zeta_i(b_i)\to e_n$ and
$\dot {\hat \zeta}_i(b_i)\to e_n$.  Therefore, there exist
$\sigma_i', \hat \sigma_i'\to 0'$ such that
$$
\zeta_i(t)\equiv \eta(\sigma_i', X_i,  t+1-b_i),
\qquad
\hat \zeta_i(t)\equiv \eta(\hat  \sigma_i', X_i,  t+1-b_i),
$$
where $\eta(\sigma', X, t)$ are the spreading geodesics we have constructed.
In particular,
$$
 \eta(\sigma_i',X_i,  1-b_i)=\zeta_i(0)=(x_i', f(x_i')),
\qquad  \eta(\hat  \sigma_i', X_i,  1-b_i)=
\hat \zeta_i(0)=(\hat x_i', f(\hat x_i')).
$$
Let $\tilde f^i$ denote the function whose graph
is the parameter sphere given by $\eta(\cdot, X_i,  1-b_i)$,
then, by the previous arguments,
 $\tilde f^i$, $\nabla \tilde f^i$ and the Hessian
 converge to corresponding things  
of $\tilde f^{(\epsilon)}$ in a fixed neighborhood
of $0'$.  
Thus, by Lemma \ref{lem4.6},
for some $\delta'>0$ independent of $i$,  
$$
(\tilde f^i-f)(x')\ge 0,\quad 
\left( (\tilde f^i-f)_{x_\alpha x_\beta}(x')\right)>0,
\qquad \forall \ |x'|<\delta' ,
$$
for large $i$.
On the other hand,
$$
 (\tilde f^i-f)(x'_i)=  (\tilde f^i-f)(\hat x_i')=0,
$$
$$
x'_i\to 0', \ \hat x'_i\to 0', \ x'_i\neq \hat x'_i.
$$
This is impossible.  
Lemma \ref{lem4-new2} is established.

\vskip 5pt
\hfill $\Box$
\vskip 5pt

We assume that $m(0)=\tilde m(0)=e_n$.
Let $\tilde f$ be the one given by Lemma \ref{lem4-4}.
Recall that  $\lambda\ge  0$ is the  smallest eigenvalue of
$(\widetilde f_{x_\alpha   x_\beta  }(0')-f_{x_\alpha  x_\beta  }(0'))$.

\begin{lem} Suppose $m(0)=e_n$ and $\lambda>0$.  Then there is a point
$Q\ne 0$ on $\partial \Omega$ whose distance to $e_n=1$.
\label{lem4.5}
\end{lem}

\noindent{\bf Proof.}\ Since $m(0)=e_n$, there is a sequence of points
$X_i\to e_n$, and $Q_i, \hat Q_i\in \partial \Omega$,
$Q_i\ne \hat Q_i$,  such that
$$
b_i:=dist(\partial \Omega\add X_i)=dist(Q_i\add 
X_i)=dist(\hat Q_i\add X_i).
$$
Passing to a subsequence,
 $Q_i\to Q$, $\hat Q_i\to \hat Q$, $b_i\to dist(\partial \Omega\add e_n)=1$.
Clearly $dist(Q\add e_n)=dist(\hat Q\add e_n)=1$.
If either $Q$ or $\hat Q$ is not $0$, we are done. Otherwise,
$Q=\hat Q=0$, and we write
$$
Q_i=(x_i', f(x_i')), \qquad \hat Q_i=(\hat x_i', f(\hat x_i')),
$$
and let $\zeta_i$ and $\hat \zeta_i$ be shortest geodesics,
with unit speed, joining respectively $Q_i$ and $\hat Q_i$ to
$X_i$.  
By Lemma \ref{lem2.2}, $\zeta_i\equiv \xi(x_i', \cdot)$ and
$\hat \zeta_i\equiv \xi(\hat x_i', \cdot)$.
So $\zeta_i\to \xi(0', \cdot)$ and
$\hat \zeta_i\to \xi(0', \cdot)$ in $C^1$ norm.
It follows that $\dot\zeta(b_i)\to e_n$ and
$\dot {\hat \zeta}_i(b_i)\to e_n$.
Therefore, there exists $\sigma_i', \hat \sigma_i'\to 0$ such that
$$
\zeta_i(t)\equiv \eta(\sigma_i', X_i, t+1-b_i),
\qquad\hat\eta(\hat\sigma_i', X_i, t+1-b_i).
$$
In particular,
$$
\eta(\sigma_i', X_i, 1-b_i)=\zeta_i(0)=(x_i', f(x_i')),
\qquad 
\eta(\hat \sigma_i', X_i, 1-b_i)=\hat \zeta_i(0)=(\hat x_i', f(\hat x_i')).
$$
Let $\tilde f^i$ denote the function whose graph is the parameter sphere
given by $\eta(\cdot, X_i, 1-b_i)$,
then the Hessian of $\tilde f^i$ converges to the
Hessian of $\tilde f$ in a fixed neighborhood of $0'$.
Thus, since  $\lambda>0$, there exists some $\delta'>0$ independent 
of $i$, such that
$$
 (\tilde f^i-f)(x')\ge 0,\quad
\left(  (\tilde f^i-f)_{x_\alpha x_\beta}(x')\right)>0,
\qquad\forall\ |x'|<\delta',
$$
for large $i$.  On the other hand,
$$
(\tilde f^i-f)(x_i')=(\tilde f^i-f)(\hat x_i')=0,
$$
$$
x_i'\to 0, \ \ \hat x_i'\to 0, \ \ x_i'\ne \hat x_i'.
$$
This is impossible.
Lemma \ref{lem4.5} is established.

\vskip 5pt
\hfill $\Box$
\vskip 5pt

\noindent{\bf 4.6.}\  In this subsection we show that
$m(0)$ is a conjugate point iff $\lambda=0$.  Since we never
apply this result the reader may choose to skip it.
\begin{lem} Suppose $m(0)=e_n$, and suppose $\lambda>0$.
Then  
 $e_n$ is not a  conjugate point
of $0$,
along the normal geodesic
 $\{te_n\ |\ 0\le t\le 1\}$, as described in Section 1.1.
\label{lem4.8}
\end{lem}

\noindent {\bf Proof. }\
We first prove that 
\begin{equation}
0\ \mbox{is an isolated point in }\
\{y\in \partial \Omega\ |\ dist(y\add e_n)=dist(\partial\Omega\add e_n)\}.
\label{isolate1}
\end{equation}
We argue by contradiction.  Suppose that for some
$x_i'\to 0'$, $x_i'\ne 0'$, we have
$$
dist( (x_i', f(x_i'))\add e_n )=dist(\partial\Omega\add e_n)=1.
$$
Let $\zeta_i$ be a shortest geodesic, with unit speed,
joining $(x_i', f(x_i'))$ to $e_n$, then, by Lemma \ref{lem2.2},
$\zeta_i\equiv \xi(x_i', \cdot)$.
So $\zeta_i\to \xi(0', \cdot)$ in $C^1$ norm, and in particular,
$\dot \zeta\to e_n$ in $C^0$ norm. Since $\lambda>0$,
$\tilde f>f$ near $0'$, and therefore, 
for some $t_i>0$, $t_i\to 0$, we find
$\zeta_i(t_i)$ on the graph of $\tilde f$.  By 
Lemma \ref{lem4-4}, 
the graph of $\tilde f$ is the distance sphere near the
origin, so $dist(\zeta_i(t_i)\add e_n)=1$.
On the other hand, since $\zeta_i$ is a  shortest geodesic with unit speed,
$$
1= t_i+(1-t_i)=t_i+dist(\zeta_i(t_i)\add e_n).
$$
This leads to contradiction.  We have thus verified (\ref{isolate1}).

The property (\ref{isolate1}) implies that $e_n$ cannot be a
conjugate point.
Indeed, if  $e_n$ is a conjugate point,
then, by (\ref{isolate1}), we may enlarge $\Omega$, without
changing $\partial \Omega$ near the origin, so that
$dist(\partial\Omega\add e_n)$ is realized only at $0$.
For this larger $\Omega$, $e_n$ is still a conjugate point
and we still have $m(0)=e_n$ for the new $\Omega$.
In the following we still use $\Omega$ to denote the
new one.  Since $e_n$ does not belong to $G$, there exist
$X_i\to e_n$, $y_i\ne z_i$,
$y_i, z_i\in \partial \Omega$,
such that
$$
b_i:=dist(y_i\add X_i)=dist(z_i\add X_i)=dist(\partial\Omega\add X_i).
$$
Passing to a subsequence, $z_i\to z$, $y_i\to y$ and $b_i\to b=1$.
Since $0$ is the only point on
$\partial \Omega$ which realizes $dist(\partial\Omega\add e_n)$,
we must have
$y=z=0$.
Write
$$
y_i=(x_i', f(x_i')), \qquad z_i=(\hat x_i', f(\hat x_i')),
$$
then $x_i'\ne \hat x_i'$.
As usual, $\xi(x_i', \cdot)$ is a
shortest geodesic joining
$y_i$ to $X_i$, $\xi(x_i', 0)=y_i, \xi(x_i', b_i)=X_i$.
Similarly,  $\xi(\hat x_i', \cdot)$ is a 
shortest geodesic joining
$z_i$ to $X_i$, $\xi(\hat x_i', 0)=z_i, \xi(\hat x_i', b_i)=X_i$.
We also know, as usual, for some $\sigma_i'$ and $\hat \sigma_i'$,
$$
\xi(x_i', t)\equiv \eta(\sigma'_i, X_i, t+1-b_i),
\quad \xi(\hat x_i', t)\equiv \eta(\hat \sigma'_i, X_i, t+1-b_i).
$$

Let $\tilde f^i$ be the function whose graph
is $\eta(\cdot, X_i, 1-b_i)$.
We argue as before:  the Hessian of $\tilde f^i$ converges to
the Hessian of $\tilde f$ in a fixed neighborhood of $0'$.
Since $\lambda>0$, $(\tilde f^i-f)$ is strictly convex 
in a fixed neighborhood of $0'$, but
we know that $\tilde f^i-f\ge 0$,
 $(\tilde f^i-f)(x_i')=(\tilde f^i-f)(\hat x_i')=0$,
$x_i'\ne \hat x_i'$, $x_i'\to 0$,  and $\hat x_i'\to 0$.
This is a contradiction.
Lemma \ref{lem4.8} is established.

\vskip 5pt
\hfill $\Box$
\vskip 5pt

Next 
\begin{lem}
Suppose $m(0)=e_n$, and suppose $\lambda=0$.  Then $e_n$ is
a conjugate point.
\label{lem4.9}
\end{lem}

A consequence of Lemma \ref{lem4.8} and \ref{lem4.9} is
\begin{cor} Suppose $m(0)=e_n$.  Then
$e_n$ is a conjugate point if and only
if $\lambda=0$.
\label{cor-new4}
\end{cor}

We present two proofs of Lemma \ref{lem4.9},
 the second one is more traditional.

\noindent{\bf First proof of Lemma \ref{lem4.9}.}\ 
Let $\zeta$ be a unit eigenvector
of $(\tilde f_{x_\alpha x_\beta}-f_{x_\alpha x_\beta})(0')$
associated with the least eigenvalue $\lambda=0$, and let
$x'\ne 0'$ be a multiple of $\zeta$.  Then
$$
|(\tilde f-f)(x')|\le C|x'|^3,
$$
and therefore
$$
dist ( (x', f(x'))\add (x', \tilde f(x'))\le C|x'|^3.
$$

Let, for some $\delta>0$,
$$
s=\delta |x'|.
$$
We will  fix some small $\delta>0$, independent of $x'$,
and show, for small $|x'|>0$,  that
\begin{equation}
dist( (x', f(x'))\add (1+s)e_n)< 1+s.
\label{C2}
\end{equation}
In fact, we will produce a curve
joining $(x', f(x'))$ to
$(1+s)e_n$ in small neighborhood
of $\{te_n\ |\ 0\le t\le (1+s)\}$ which has length less than
$1+s$.
This means that $e_n$ is a conjugate point.

By the  triangle inequality, and using
$\tilde \xi(x', 0)=(x', \tilde f(x'))$, 
\begin{eqnarray}
&&dist ( (x', f(x'))\add (1+s)e_n)\nonumber\\
&\le & dist (  (x', f(x'))\add (x', \tilde f(x')))
+dist ( (x', \tilde f(x'))\add \tilde \xi( x', 1-s))\nonumber\\
&&
+dist ( \tilde \xi(x', 1-s)\add  (1+s)e_n)\nonumber\\
&\le & C|x'|^3+ (1-s)+
dist ( \tilde \xi(x', 1-s)\add  (1+s)e_n).
\label{C3}
\end{eqnarray}

Since $\tilde \xi(x', \cdot)$ and $\tilde \xi(0', \cdot)$
 satisfy the same geodesic equations, we have
\begin{eqnarray*}
|\tilde \xi(x', 0)-\tilde \xi(0', 0)|
&\le & C |\tilde \xi(x', 1)-\xi(0', 1)|
+C|\dot{\tilde \xi}(x', 1)- \dot{\tilde \xi}(0', 1)|\\
&=& C|\dot{\tilde \xi}(x', 1)-e_n|.
\end{eqnarray*}
It follows, for some $c>0$ depending only on
$f$ and $\varphi$,   that
\begin{equation}
|\bar e-e_n|\ge
c|x'|,
\label{C4}
\end{equation}
where
$$\bar e:= \dot{\tilde \xi}(x', 1).
$$
Now,
a crucial point:  
since $\varphi(e_n; \bar e) 
=\varphi(\tilde \xi(x', 1); \dot{\tilde \xi}(x', 1))
=\varphi(e_n; e_n)=1$, by the strict
convexity hypothesis on $\psi$, we have for some
$\hat c_0>0$ depending only on 
$\varphi$, that
\begin{equation}
\varphi(e_n; \frac {e_n+\bar e}2)\le 1- \hat c_0 |e_n-\bar e|.
\label{C5}
\end{equation} 

Let
$$
\eta(t)=(1-t)\tilde \xi(x', 1-s)+t(1+s)e_n, \qquad 0\le t\le 1,
$$
be the straight segment joining $\tilde \xi(x', 1-s)$ to
$(1+s)e_n$.  Then, since $\tilde \xi(x', 1)=e_n$, 
$$
\eta(t)=e_n+O(s).
$$
Here and below $O(s)$ denotes some quantity
which is bounded in absolute value by
$Cs$ for some constant $C$ independent of
$x'$ and $s$.

Using $\ddot{\tilde \xi}(0', \cdot)\equiv 0$,
\begin{eqnarray*}
\dot\eta(t)&=& (1+s)e_n-\tilde \xi(x', 1-s)\\
&=& (1+s)e_n-[ \tilde \xi(x', 1)
-\dot{\tilde \xi}(x', 1)s+O(|x'|s^2)\\
&=& s e_n+ \dot{\tilde \xi}(x', 1)s+O(|x'|s^2)
=s(e_n+\bar e)+O(|x'|s^2).
\end{eqnarray*}
It follows, using properties of our special coordinates and the homogeneity
of $\varphi$ in $v$,
and making Taylor expansions, that
\begin{eqnarray*}
&& dist(\tilde \xi(x', 1-s)\add (1+s)e_n)\\
&\le & \int_0^1\varphi(\eta; \dot \eta)dt
=s\int_0^1 \varphi(e_n+O(s);
e_n+\bar e+O(|x'|s))dt\\
&=& s\left(\varphi(e_n; e_n+\bar e)+O(s)\right), 
\end{eqnarray*}
and therefore, by (\ref{C4}) and (\ref{C5}),
$$
dist(\tilde \xi(x', 1-s)\add (1+s)e_n)
\le 2s(1-\hat c_0|x'|+O(s)),
$$
where $\hat c_0>0$ is some constant independent of
$x'$ and $s$.

Back to (\ref{C3}), we find 
\begin{eqnarray*}
dist ( (x', f(x'))\add (1+s)e_n)&\le&
C|x'|^3+ (1-s)+s(2-\hat c_0|x'|+O(s))\\
&=&1+s-s(2\hat c_0|x'|+O(s))+C|x'|^3.
\end{eqnarray*}
Now we fix some $\delta>0$ from the beginning so that
$2\hat c_0|x'|+O(s)\ge \hat c_0|x'|$, 
then for $|x'|>0$ small, we obtain
$$
dist ( (x', f(x'))\add (1+s)e_n)\le
1+s-\hat c_0\delta |x'|^2+C|x'|^3
<1+s,
$$
the estimate (\ref{C2}).  It is clear that we have actually
produced a curve in small neighborhood of
$\{te_n\ |\ 0\le t\le 1+s\}$ joining 
$ (x', f(x'))$ to $(1+s)e_n$ with length less than $1+s$. 
Lemma \ref{lem4.9} is established.

\vskip 5pt
\hfill $\Box$
\vskip 5pt

Now we present the 

\noindent{\bf Second proof of Lemma \ref{lem4.9}.}\ 
(i)\ Consider the spreading geodesics $\eta(\sigma', t)$.  Since
$\tilde f\ge f$, for small $\sigma'$, there exists
a unique $\bar t(\sigma')\le 0$ such that $\eta(\sigma', \bar t(\sigma'))$
lies on $\partial \Omega$, i.e.
\begin{equation}
\eta^n(\sigma', \bar t(\sigma'))=f(\eta'(\sigma', \bar t(\sigma'))).
\label{1}
\end{equation}
The function $\bar t(\sigma')$ is a $C^2$ function in $\sigma'$.
The curve $\{\eta(\sigma', t) \ | \ \bar t(\sigma')\le t\le 1\}$ has
length 
\begin{equation}
L(\sigma')=1-\bar t(\sigma').
\label{2}
\end{equation}
We also have
\begin{equation}
\eta^n(\sigma', 0)=\tilde f(\eta'(\sigma', 0)).
\label{3}
\end{equation}
Differentiating (\ref{1}) w.r.t. $\sigma_\alpha$ we find
$$
\eta^n_{\sigma_\alpha}+\dot \eta^n \bar t_{\sigma_\alpha}
=f_{ x_\gamma}(\eta^\gamma_{\sigma_\alpha}+
\dot \eta^\gamma \bar t_{\sigma_\alpha}).
$$
Differentiate w.r.t. $\sigma_\beta$, and set $\sigma'=0'$.
We get at $\sigma'=0'$,
\begin{equation}
\eta^n_{\sigma_\alpha\sigma_\beta}+
\bar t_{\sigma_\alpha\sigma_\beta}=f_{ x_\gamma x_\delta}(0')
\eta^\gamma_{\sigma_\alpha} \eta^\delta_{ \sigma_\beta}\qquad \mbox{at}\ (0',0),
\label{4}
\end{equation}
since, when $\sigma'=0'$,
$ \eta^n_{\sigma_\alpha}=0$ (following from $\tilde f_{x_\beta}(0')=0$),
and so, $\bar t_{\sigma_\alpha}(0')=0$.

Similarly, from (\ref{3}), we find
\begin{equation}
\eta^n_{ \sigma_\alpha\sigma_\beta}
=\tilde f_{ x_\gamma x_\delta}\eta^\gamma_{\sigma_\alpha}
\eta^\delta_{\sigma_\beta}\qquad \mbox{at}\
(0',0).
\label{5}
\end{equation}
So,
\begin{equation}
\bar t_{\sigma_\alpha\sigma_\beta}(0')=(f_{x_\gamma x_\delta}(0')
-\tilde f_{x_\gamma x_\delta}(0'))
\eta^\gamma_{\sigma_\alpha}(0',0)\eta^\delta_{\sigma_\beta}(0',0).
\label{6}
\end{equation}

Suppose, now, $\lambda=0$.  Then there is a unit vector $\hat \zeta=(\hat \zeta^1, \cdots, \hat \zeta^{n-1})$ such that
\begin{equation}
(f_{x_\gamma x_\delta}(0')
-\tilde f_{x_\gamma x_\delta}(0'))
\hat \zeta^\delta=0, \qquad 1\le \gamma\le n-1.
\label{7}
\end{equation}
The matrix $\{\eta^\gamma_{\sigma_\alpha}(0',0)\}$ is nonsingular.
Choose $a=(a_1, \cdots, a_{n-1})$ so that
\begin{equation}
a_\alpha \eta^\delta_\alpha(0', 0)=\hat \zeta^\delta.
\label{8}
\end{equation}
Inserting this in (\ref{7}), we find, by (\ref{6}),
\begin{equation}
a_\alpha a_\beta \bar t_{\sigma_\alpha \sigma_\beta}(0')=0.
\label{9}
\end{equation}

(ii)\ Now the second variation.  For $0\le t\le 1$,
$te_n$ is the shortest connection from $\partial \Omega$ to $e_n$.
For $\zeta(t)$ small, $0\le t\le 1$, we consider
the perturbation $te_n+\zeta(t)$.  Here $\zeta(t)=0$ and
$\zeta(0)\in \partial \Omega$, i.e.,
\begin{equation}
\zeta^n(0)=f(\zeta'(0))=\frac 12
f_{x_\gamma x_\delta}(0')\zeta^\gamma(0)\zeta^\delta(0)+O(|\zeta'(0)|^3).
\label{10}
\end{equation}
The length of the curve $te_n+\zeta$ is, by the properties
of our special coordinates, 
\begin{eqnarray*}
&&\int_0^1\varphi(te_n+\zeta; e_n+\dot \zeta)dt
\\
&=&1+ \int_0^1(\frac 12
\varphi_{ \eta^\alpha \eta^\beta}(te_n; e_n)
\zeta^\alpha\zeta^\beta+\frac 12\varphi_{ v^\alpha v^\beta}(te_n; e_n)
\dot \zeta^\alpha
\dot \zeta^\beta+\dot \zeta^n)dt
+\mbox{higher order}.
\end{eqnarray*}
So the second variation is, by (\ref{10}), 
\begin{eqnarray}
Q(\zeta')&:=&
\frac 12
\int_0^1  (\varphi_{ \eta^\alpha \eta^\beta}(te_n; e_n)
\zeta^\alpha\zeta^\beta+\varphi_{ v^\alpha v^\beta}(te_n; e_n)
\dot \zeta^\alpha
\dot \zeta^\beta)dt\nonumber\\
&&
- \frac 12
f_{x_\alpha x_\beta}(0')\zeta^\alpha(0)\zeta^\beta(0).
\label{11}
\end{eqnarray}

Now $\{\varphi_{v^\alpha v^\beta}(te_n; e_n)\}$
is positive definite, and the 
quadratic form $Q(\zeta')$ is
positive semidefinite.  If it vanishes
for some $\zeta'(t)$ not identically zero, then
$e_n$ is a conjugate point--by the
usual argument:  the second
variation of the curve
$te_n$ for $0\le t\le 1+\epsilon$,
for any $\epsilon>0$, is not positive semidefinite.

(iii)\  Suppose $\lambda =0$.

\noindent{\bf Claim:}\ For $\zeta^\alpha(t)=a_\gamma 
\eta^\alpha_{\sigma_\gamma}(0',t), \qquad Q(\zeta')=0$.

\medskip

This would then complete the proof of the lemma.

\noindent{\bf Proof of Claim.}\  From (\ref{2}), we have
\begin{equation}
L_{\sigma_\alpha\sigma_\beta}(0')=-\bar t_{\sigma_\alpha\sigma_\beta}(0').
\label{13}
\end{equation}
Now
$$
L(\sigma')=\int_{\bar t(\sigma') }^1
\varphi(\eta(\sigma', t); \dot\eta(\sigma', t))dt,
$$
and recall that $\varphi(\eta(\sigma', t); \dot\eta(\sigma', t))
\equiv 1$.
So
$$
L_{\sigma_\alpha}=-\bar t_{\sigma_\alpha}
+ \int_{\bar t }^1
(\varphi_{ \eta^i}\eta^i_{\sigma_\alpha}
+\varphi_{ v^i}\dot
\eta^i_{\sigma_\alpha})dt
$$
and, at $\sigma'=0'$, by properties of the special
coordinates,
$$
L_{\sigma_\alpha\sigma_\beta}(0')=
-\bar t_{\sigma_\alpha\sigma_\beta}+
\int_0^1(\varphi_{\eta^i\eta^j}\eta^i_{\sigma_\alpha}
\eta^j_{\sigma_\beta}+
\varphi_{v^iv^j}\dot \eta^i_{\sigma_\alpha}
\dot \eta^j_{\sigma_\beta}+\dot\eta^n_{ \sigma_\alpha\sigma_\beta})dt.
$$
By (\ref{13}), the last integral is zero.  
By properties of the special coordinates and
by the homogeneity of $\varphi$ in $v$, we have
$\varphi_{\eta^i\eta^n}(te_n; e_n)\equiv
\varphi_{v^iv^n}(te_n; e_n)\equiv 0$.  Therefore
we have 
$$
\int_0^1(\varphi_{\eta^\gamma\eta^\delta}(te_n;e_n)
\eta^\gamma_{\sigma_\alpha}(0',t)\eta^\delta_{\sigma_\beta}(0',t)
+
\varphi_{v^\gamma v^\delta}(te_n;e_n)
\dot \eta^\gamma_{\sigma_\alpha}(0',t)\dot \eta^\delta_{\sigma_\beta}(0',t)
-\eta^n_{ \sigma_\alpha\sigma_\beta}(0',0)=0.
$$
Multiplying the above by $a_\alpha a_\beta$ and summing, we find
\begin{equation}
\int_0^1(\varphi_{\eta^\gamma\eta^\delta}(te_n;e_n)
\zeta^\gamma\zeta^\delta+
\varphi_{v^\gamma v^\delta}(te_n;e_n)
\dot \zeta^\gamma \dot \zeta^\delta
-\eta^n_{ \sigma_\alpha\sigma_\beta}(0',0)a_\alpha a_\beta=0.
\label{14}
\end{equation}
From (\ref{4}) and (\ref{9}), we have
$$
-\eta^n_{ \sigma_\alpha\sigma_\beta}(0',0)a_\alpha a_\beta
=- f_{ x_\gamma x_\delta}(0')\zeta^\gamma
\zeta^\delta.
$$
Inserting this into (\ref{14}) we obtain the Claim.

\vskip 5pt
\hfill $\Box$
\vskip 5pt

\section{Main Estimates I}
\setcounter{equation}{0}

 We now start the argument described in
Section 1.5, with $y$ as the origin.  Without loss of generality,
we may assume $\bar s(y)=\bar s(0)=1$.  Then we use
our special coordinates of Section 3; near
the origin  $\Omega$ is given by
$x_n>f(x')$ with
$$
f(0')=0, \qquad \nabla f(0')=0.
$$
Then $m(y)=m(0)=e_n$.  The ``normal'' geodesic
from $0$ lies along the $x_n-$axis.

For $|x'|$ small, as in Section 2, $\xi(x',\tau)$
is the geodesic, with $\tau$ as arclength, starting from
$(x', f(x'))$ normal to $\partial \Omega$.
We wish to find a point $z$ on $\partial \Omega$
such that for $s=K|x'|$, with $K$ a fixed
large constant,
\begin{equation}
dist(z \add \xi(x', 1+s))<1+s.
\label{4.1}
\end{equation}
To prove (\ref{4.1}) we will follow the interior ``normal'' geodesic
from $z$ a distance $1-s$, then join its end
point by a straight line segment $\eta(t)$,
$0\le t\le 1$, to $\xi(x', 1+s)$,
and show that the Finsler length of $\eta$ is less
than $2s$.

To compute lengths we use expansions
in  $x', s$ etc.; the special coordinates
make the computations easier.  But things are
not very easy.

For $\epsilon_0>0$, let $\Gamma
:=\{t e_n\ |\ -\epsilon_0\le t\le 1+\epsilon_0\}$
be the geodesic for $\varphi(\xi;v)$   satisfying,
for $-\epsilon_0\le t\le  1+\epsilon_0$, (\ref{ag5}),
(\ref{g2n}), 
(\ref{g6n}), 
(\ref{g3n}),
(\ref{sss1}),
and (\ref{aaa2}).
We use notation $\xi(x', \tau)$ as in Section 4.

\begin{lem} Under the above hypotheses, 
\begin{equation}
\xi^n_{x_\alpha  }=0, \qquad \mbox{at}\
(0', t)\ \forall\ -\epsilon_0\le t\le 1+\epsilon_0, \
1\le\alpha\le n-1,
\label{43-1}
\end{equation}
and, for $|x'|\le \epsilon_1$,
$-\epsilon_0\le t\le 1+\epsilon_0$,
$1\le \alpha, \beta\le n-1$,
\begin{equation}
|\xi^n_{x_\alpha   x_\beta  }(x',t)+\varphi_{v^iv^j}(te_n;e_n)
\dot \xi^i_{x_\alpha  }(0',t)\xi^j_{x_\beta  }(0',t)|
\le C|x'|, 
\label{43-2}
\end{equation}
where $C$ depends only on
$f$ and $\varphi$.
\label{lem4-7}
\end{lem}

\noindent{\bf Proof of Lemma \ref{lem4-7}.}\
By (\ref{2.14}),
$$
\xi^i_{x_\alpha  }\varphi_{v^i}(\xi; \dot\xi)\equiv 0.
$$
The first equality in the lemma follows easily from the above
by the properties of the special coordinates.
Applying $\partial_{x_\beta  }$ to the above, we have
$$
\xi^n_{x_\alpha   x_\beta  }
\varphi_{v^n}= -\xi_{x_\alpha   x_\beta  }^\gamma\varphi_{v^\gamma}
-\varphi_{v^i\xi^j}  \xi^i_{x_\alpha  } 
 \xi^j_{x_\beta  }-\varphi_{v^iv^j}  \xi^i_{x_\alpha  } \dot \xi^j_{x_\beta  }.
$$
At
 $x'=0$,  using properties of the 
special coordinates, we have
$$
\varphi_{v^n}(te_n;e_n)=1,\ 
\varphi_{v^\gamma}(te_n;e_n)=0, \
\varphi_{v^i\xi^j}(te_n;e_n)=0,
$$
and estimate (\ref{43-2}) follows.

\vskip 5pt
\hfill $\Box$
\vskip 5pt

We assume that
$$
\tau=dist(0\add \tau e_n)=
\min_{y\in \partial \Omega}dist(y\add \tau e_n),
\qquad
\forall\ 0\le \tau\le 1,
$$
and
$$
\epsilon_0= dist( (0', -\epsilon_0)\add 0).
$$
In particular, 
$$
\tau=dist(0\add \tau e_n)\le dist((x',f(x'))\add \tau e_n),\qquad
\forall\ 0\le \tau\le 1, |x'|\le \epsilon_1,
$$

\bigskip

We 
now find $z$,
for our program,  in the simplest case
\begin{prop}
Assume that there exists
 $Q\in \partial \Omega$, $|Q|\ge\hat \epsilon>0$,
with
$$
dist(0\add e_n)=dist(Q\add e_n).
$$
Then, we take $z=Q$, i.e.,  there exist some large  
constant $K\ge 1$ and small constant $0<\hat \delta<\hat\epsilon$,
depending only on $\hat \epsilon$, $f$ and $\varphi$
 such that
for all $0<|x'|\le \hat \delta$ and $s=K|x'|$ we have
$$
dist(Q\add \xi(x', 1+s))<1+s.
$$
\label{prop4-0}
\end{prop}

\noindent{\bf Proof of Proposition \ref{prop4-0}.}\
Set $\bar e=\dot \xi(Q, 1)$.  Since $\xi(Q,1)=e_n$,
and the fact that $\xi(Q,s)$ satisfies the geodesic equations, it follows that
\begin{equation}
|\bar e_n-e_n|=
|\dot \xi(Q, 1)-e_n|\ge c_1 |Q|\ge c_1\hat \epsilon
\label{4.10}
\end{equation}
for some $c_1>0$ depending only on $\varphi$.
  We know that
$$
\xi(0',1)=e_n,\quad
\xi(Q,1)=e_n,\quad
\mbox{and}\quad
\dot \xi(0',1)=e_n.
$$
By Taylor expansion, since $|x'|=\frac sK\le s$, 
$$
\xi(x', 1+ s)=\xi(0',1)+O(s)=e_n+O(s),
$$
$$
\xi(Q, 1-s)=\xi(Q,1)-\dot \xi(Q,1)s+O(s^2)
=e_n-s \bar e+O(s^2).
$$
For the segment
$$
\eta(t):=(1-t)\xi(Q, 1-s)+t\xi(x',1+s)
=e_n+O(s),
$$
\begin{eqnarray*}
\dot \eta(t)&=&\xi(x',1+s)-\xi(Q, 1-s)
=\xi_{x_\alpha  }(0',1)x_\alpha  +
\dot\xi(0',1)s+\dot\xi(Q,1)s+O(s^2+|x'|^2)\\
&=&\xi_{x_\alpha  }(0',1)x_\alpha  +(e_n+\bar e)s+O(s^2+|x'|^2).
\end{eqnarray*}

Using homogeneity, it follows that
\begin{eqnarray*}
\int _0^1\varphi(\eta(t);\dot\eta(t))dt
&=&s\int _0^1\varphi\big(e_n+O(s);
(e_n+\bar e)+
 \xi_{x_\alpha  }(0',1)\frac{x_\alpha  }s+O(s)
+O(\frac{|x'|^2}s)\big)
dt\\
&=&s
 \varphi(e_n, e_n+\bar e) +O(\frac sK+s^2)=
2s \varphi(e_n, \frac { e_n+\bar e }2)+O(\frac sK+s^2). 
\end{eqnarray*}

Now, the crucial point as in the proof of Lemma \ref{lem4.9}.
Since $\varphi(e_n; e_n)=\varphi(e_n; \bar e_n)=1$, by the 
strict
convexity hypothesis on $\psi$, we have
 for some $\bar c_0$ depending only
on $\varphi$, that 
$$
 \varphi(e_n; \frac { e_n+\bar e }2)\le
1-\bar c_0|e_n-\bar e_n|\le
 1-c_0 
$$
with $c_0>0$  depending also on $\hat \epsilon$----by (\ref{4.10}), 
from which we deduce, for some large  $K$ and small $\hat \delta$
($K$ chosen first and then $\hat \delta$), that
$$
\int _0^1\varphi(\eta(t);\dot\eta(t))dt
\le 2s(1-c_0)+O(\frac sK+s^2)\le  2s(1-\frac {c_0}2)<2s.
$$
Consequently,
$$
dist(Q\add \xi(x', 1+s))\le
dist(Q\add \xi(Q, 1-s))+
dist(\xi(Q, 1-s)\add \xi(x', 1+s))<1+s.
$$
Proposition \ref{prop4-0} is established.

\vskip 5pt
\hfill $\Box$
\vskip 5pt

\section{Main Estimates II}
\setcounter{equation}{0}

  In the remaining cases we will take
$$
z=(x'+q, f(x'+q))
$$
for suitable choices of $q\in \Bbb R^{n-1}$,
$|q|<$small. In the following, the value of
$\epsilon_0$ is possibly smaller than the one appearing in Section 4.1.

We know that
$$
\xi(0',\tau)=\tau e_n, \qquad
-\epsilon_0\le \tau\le 1+\epsilon_0.
$$

For $x', x'+q\in \Bbb R^{n-1}$,
$|x'|, |x'+q|\le \epsilon_1$, let
$\eta(x', q, s; t)$, $0\le t\le 1$, denote
the straight segment going from
$\xi(x'+q, 1-s)$ to
$\xi(x', 1+s)$.  We consider its length,
$L(x', q, s)$,  as a function of $2(n-1)+1$
variables, the $x', q, s$ being free variables (with small
norms).  Thus 
\begin{equation}
\eta(x', q, s; t)=(1-t)
\xi(x'+q, 1-s)+t\xi(x', 1+s),
\qquad 0\le t\le 1,
\label{5.1}
\end{equation}
and
\begin{equation}
L(x', q, s)=\int_0^1 \varphi(\eta(x', q, s; t);
\dot \eta(x', q, s; t))dt
\label{5.2}
\end{equation}
where $\dot {} $ denotes $\partial_t$.

For suitable choice of $q$, and with
$s=K|x'|$, $K$ large, we wish to show
$$
L(x', q, s)<2s.
$$
The main term will be $L(0', q,s)$.
Proposition \ref{prop5.1} below presents a general
estimate for the difference.
This result is rather technical; it will be used for several cases.
We stress that $x', q,s$ are free variables.  The expression
$O(|q|)$ is used to denote quantities bounded 
in absolute value by $C|q|$, where $C$ depends only
on $f$ and $\varphi$.

The vector 
\begin{equation}
A=e_n-\xi(q,1)
\label{5.3}
\end{equation}
plays an important role.  Note that
\begin{equation}
A^j=\delta^j_n-\xi^j(q,1)=-\xi^j_{x_\alpha}(0',1)q_\alpha+O(|q|^2).
\label{5.4}
\end{equation}

\begin{prop}
There exist $\bar \epsilon_1\le \epsilon_1$
and $0<\epsilon_3$, depending only on $f$ and $\varphi$, such that
$\forall \ x',q,s$ satisfying
$|x'|, |q|, |x'+q|, s\le
\bar \epsilon_1$, $s>0$, and if
\begin{equation}
\frac {|A|}s<\frac 14\quad\mbox{and}\quad \frac {|x'|}s\le \epsilon_3,
\label{5.5}
\end{equation}
then we have
\begin{eqnarray}
J&:=&L(x', q, s)-L(0', q, s)\nonumber\\
&\le& C|x'|^2\left(|q|+s+\frac {|q|^2}s\right)
+C|x'|\left( |A|(1+\frac {|q|}s)
+|q|^2+s^2\right).
\label{5.6}
\end{eqnarray}
\label{prop5.1}
\end{prop}

\noindent{\bf Proof.}\ In formula (\ref{5.2}), $\eta$ is given
by (\ref{5.1}) and
$$
\dot \eta=\xi(x', 1+s)-\xi(x'+q, 1-s).
$$
Clearly
\begin{equation}
\eta=e_n+O(|q|+|x'|+s),
\label{5.7}
\end{equation}
while
\begin{eqnarray*}
\dot\eta&=& e_n(1+s)+O(|x'|)-\xi(q,1-s)\\
&=&  e_n(1+s)-\xi(q,1)+\dot\xi(q,1)s
+O(s^2)+O(|x'|).
\end{eqnarray*}
Thus
\begin{equation}
\dot \eta=2se_n+s\left(\frac As+B\right),
\label{5.8}
\end{equation}
where
$$
|B|\le C_1(|q|+s+\frac {|x'|}s)
$$
with $C_1$ depending only on
$f$ and $\varphi$.  We now make
$|B|\le 1/2$ by choosing
\begin{equation}
\bar \epsilon_1=\min(\epsilon_1, \frac 1{8C_1}),
\qquad \epsilon_3=\frac 1{4C_1}.
\label{5.9}
\end{equation}

In addition to (\ref{5.8}) we have
\begin{equation}
D_{x'}^k\dot \eta=
D_{x'}^k \xi(x', 1+s)-D_{x'}^k \xi(x'+q, 1-s)=O(|q|+s),
\qquad 0\le k\le 2.
\label{5.10}
\end{equation}

Using Taylor expansion in $x'$ about the origin,
we have
\begin{equation}
J=L(x', q,s)-L(0', q,s)=
L_{x_\alpha}(0,q,s)x_\alpha+\int_0^1\int_0^t
L_{x_\alpha x_\beta}(\tau x', q,s)x_\alpha
x_\beta d\tau dt.
\label{5.11}
\end{equation}

Now
\begin{equation}
L_{x_\alpha}(x',q,s)=\int_0^1(\varphi_{\xi^i}\eta_{x_\alpha}^i+
\varphi_{v^i}\dot \eta_{x_\alpha}^i)dt
\label{5.12}
\end{equation}
and
\begin{eqnarray*}
L_{x_\alpha x_\beta}(x',q,s)&=&
\int_0^1\bigg[
\varphi_{\xi^i}\eta^i_{x_\alpha x_\beta}+\varphi_{\xi^i\xi^j}
\eta_{x_\alpha}^i \eta_{x_\beta}^j+
\varphi_{\xi^i v^j}\eta_{x_\alpha}^i
\dot \eta_{x_\beta}^j\\
&&\quad +\varphi_{v^i}\dot \eta_{x_\alpha x_\beta}^i
+\varphi_{ v^i\xi^j} \dot \eta_{x_\alpha}^i  \eta_{x_\beta}^j
+\varphi_{v^iv^j} 
\dot \eta_{x_\alpha}^i  \dot  \eta_{x_\beta}^j\bigg]dt.
\end{eqnarray*}

By the properties of the special
coordinates, at $(e_n, e_n)$,
$$
\varphi_{v^\alpha}=\varphi_{\xi^i}=
\varphi_{ \xi^i v^j}=
\varphi_{v^n}-1=\varphi_{v^\alpha v^n}=\varphi_{\xi^j\xi^n}=0.
$$
Thus, by (\ref{5.7}), (\ref{5.8}), (\ref{5.10}),
and the homogeneity of $\varphi$ in $v$, if we set
$$
\{\qquad\}=
\frac {|A|}s+\frac {|x'|}s+|q|+s,
$$
we find, at $(\eta, \dot\eta)$, that
$$
|\varphi_{ \xi^i}|+|\varphi_{ \xi^i\xi^n}|\le Cs\{\qquad\},
$$
$$
|\varphi_{\xi^i v^j}|+
|\varphi_{v^\alpha}|+|\varphi_{v^n}-1|\le
C\{\qquad\},
$$
$$
|\varphi_{v^\alpha v^n}|\le \frac Cs\{\qquad\},
$$
$$
|\varphi_{\xi^i\xi^j}|+|\varphi_{\xi^i\xi^j\xi^k}|\le Cs,
$$
$$
|\varphi_{\xi^i\xi^j v^k}|\le C,
$$
$$
|\varphi_{ \xi^i v^j v^k}|+|\varphi_{v^i v^j}|\le \frac Cs,
$$
$$
|\varphi_{v^i v^j v^k}|\le\frac C{s^2}.
$$

We deduce from the above, since $|\{\qquad\}|$ is bounded,
that 
$$
|L_{x_\alpha x_\beta}(\tau x', q,s)|\le C(|q|+s)+\frac Cs(|q|+s)^2.
$$
Consequently
\begin{equation}
|\int_0^1\int_0^t L_{x_\alpha x_\beta}(\tau x', q,s)x_\alpha x_\beta d\tau dt|
\le C|x'|^2 (|q|+s+\frac {|q|^2}s).
\label{5.13}
\end{equation}

Next, we estimate $L_{x_\alpha}(0', q, s)x_\alpha$.
Here $x'=0'$ in $(\eta, \dot \eta)$.
By the estimates above,
\begin{equation}
|\int\varphi_{\xi^i}\eta_{x_\alpha}^i x_\alpha|
\le C( |A|+s|q|+s^2)|x'|,
\label{5.14}
\end{equation}
$$
|\int \varphi_{v^\beta}\dot \eta_{x_\alpha}^\beta x_\alpha|
\le C(\frac {|A|}s+|q|+s)(|q|+s)|x'|.
$$

Write
$$
\varphi_{v^n}\dot \eta_{x_\alpha}^n=\dot\eta_{ x_\alpha}^n+(\varphi_{v^n}-1)
\dot\eta_{ x_\alpha}^n.
$$
Then, using the estimates on $(\varphi_{v^n}-1)$ and
on $|\dot\eta_{ x_\alpha}^n|$, we find
\begin{equation}
|\int(\varphi_{v^n}-1) \dot\eta_{ x_\alpha}^n x_\alpha|\le
C(\frac {|A|}s+|q|+s)(|q|+s)|x'|.
\label{5.15}
\end{equation}

To complete the estimate of $L_{x_\alpha}(0',q,s)x_\alpha$ we
need to estimate $|\dot \eta_{x_\alpha}^n (0',q,s)x_\alpha|$.  
Using Taylor expansion, we find
\begin{eqnarray*}
\dot \eta_{x_\alpha}^n (0',q,s)&=& \xi_{x_\alpha}^n (0',1)
- \xi_{x_\alpha}^n (q,1)+s(\dot \xi_{x_\alpha}^n (0',1)
+\dot \xi_{x_\alpha}^n (q,1))+O(s^2)\\
&=& \xi_{x_\alpha}^n (0',1)
- \xi_{x_\alpha}^n (q,1)+O(s|q|+s^2),
\end{eqnarray*}
since $\dot \xi_{x_\alpha}^n(0',1)=0$,
which follows from differentiating (\ref{43-1}).
Writing
$$
\xi_{x_\alpha}^n(0',1)-\xi_{x_\alpha}^n(q,1)
=-\int_0^1 \xi_{x_\alpha x_\beta}^n(\tau q, 1)q_\beta d\tau,
$$
we find, using (\ref{43-2}), that
$$
\dot\eta_{x_\alpha}^n(0',q,s)
=\varphi_{v^i v^j}(e_n; e_n)
\dot \xi_{x_\alpha}^i(0',1)\xi_{x_\beta}^j(0',1)q_\beta
+O(|q|^2+s^2).
$$
With the aid of (\ref{5.4}), we see that
$$
\dot\eta_{x_\alpha}^n(0',q,s)=O\left(\frac {|A|}s(|q|+s)
+|q|^2+s^2\right)
$$
so that
\begin{equation}
|\int\varphi_{v^n}\dot \eta _{x_\alpha}^n x_\alpha|
\le C\left( |A|(\frac {|q|}s+1)+|q|^2+s^2\right)|x'|.
\label{5.16}
\end{equation}

Combining all the estimates (\ref{5.14}), (\ref{5.15}), (\ref{5.16})
and (\ref{5.13}), we obtain (\ref{5.6}).

\vskip 5pt
\hfill $\Box$
\vskip 5pt

\section{Main Estimates III}
\setcounter{equation}{0}

Recalling $\bar\epsilon_1$ of Proposition
\ref{prop5.1}, we now consider
the case there is a $\hat q$ satisfying the condition on $q$
of Proposition \ref{prop5.1} and in addition
$$
1=dist(0\add e_n)=dist( (\hat q, f(\hat q))\add e_n).
$$
In this case the vector $A$ of Section 6 is zero.
We take
$$
z=(\hat q+x', f(\hat q+x')).
$$

\begin{prop}
Under the conditions above, there exist
small positive constants $\bar \epsilon,
\bar\delta$ and a large constant $K>1$, depending only on
$\varphi$ and $f$ such that for
$s=K|x'|$ and $0<|x'|\le\min(\bar\delta, \bar\epsilon|\hat q|)$ we have
$$
L(x', \hat q, s)<2s.
$$
\label{prop6.1}
\end{prop}

\noindent{\bf Proof.}\  We will apply Proposition \ref{prop5.1}
with $q=\hat q$.  Since $A=0$, wee see that the conditions
are satisfied provided
$\frac 1K\le \epsilon_3$.  Then, from (\ref{5.6}) we find
\begin{equation}
J=L(x', \hat q,s)-L(0', \hat q,s)\le
Cs|\hat q|^2\left( \frac {\bar\epsilon}K+\bar\epsilon^2
+\frac 1{K^2}+\frac 1K+\bar\epsilon^2 K\right).
\label{6.1}
\end{equation}

We now consider the main term $L(0', \hat q, s)$.  The estimate is technical.  A crucial element, as
in the proof of Proposition \ref{prop4-0}, is the strict
convexity of $\{v\ |\ \varphi(e_n;v)=1\}$, and
the fact that
$$
1=\varphi(e_n;e_n)=\varphi(\xi(\hat q, 1); \dot\xi(\hat q,1))=\varphi(e_n; 
\dot \xi(\hat q,1)).
$$

By the strict convexity it follows that for some $c_1>0$,
depending only on $\varphi$,
\begin{equation}
\varphi(e_n; e_n+\dot\xi(\hat q, 1))\le 2-
2c_1|\dot\xi(\hat q,1)-e_n|^2.
\label{6.2}
\end{equation}
Since $\xi(\hat q,\cdot)$ satisfies the geodesic equations,
 $\xi(\hat q,1)=\xi(0',1)=e_n$,
and $|\xi(\hat q,0)-\xi(0',0)|=|\hat q|$, there are positive
constants $c_2, c_3$ so that
$$
c_2|\hat q|\le |\dot\xi(\hat q,1)-e_n|
=|\dot \xi(\hat q,1)-\dot \xi(0', 1)|\le c_3|\hat q|.
$$
Inserting this in (\ref{6.2}) we find, for some $c_0>0$ depending only
on $\varphi$,
\begin{equation}
\varphi(e_n; e_n+\dot\xi(\hat q, 1))\le
2-2c_0|\hat q|^2.
\label{6.3}
\end{equation}

\begin{lem} There exist positive constants  $c_0,
C$, depending only on $\varphi$ such that 
for all $0<s<\bar\epsilon_1$ and $0<|\hat q|<\bar\epsilon_1$
above,
\begin{equation}
L(0', \hat q, s)\le 2s(1-c_0|\hat q|^2)
+C(s^4+s^2|\hat q|^2).
\label{6.4}
\end{equation}
\label{lem6.1}
\end{lem}

\noindent {\bf Proof of Lemma \ref{lem6.1}.}\
Let 
$$
\eta(t)=\eta(s,t)=(1-t)\xi(\hat q,1-s)+t(1+s)e_n.
$$
Then
$$
L(0', \hat q, s)=
\int_0^1 \varphi(\eta(t);\dot{\eta}(t))dt.
$$
Since 
$$
\xi(\hat q, 1)=e_n \ \ \mbox{and}\
\ddot{\xi}(0',\tau)\equiv 
\frac{\partial ^3}{\partial \tau^3}\xi(0',\tau)\equiv 0,
$$
we have, by Taylor expansion, that
$$
\xi(\hat q, 1-s)=e_n-\dot{\xi}(\hat q,1)s+
\frac 12 \ddot{\xi}(\hat q,1)s^2+O(s^3|\hat q|),
$$

$$
\eta(t)=
e_n+s[te_n-(1-t) \dot{\xi}(\hat q,1)]+ 
O(s^2|\hat q|),
$$
$$
\dot{\eta}(t)=
s[e_n+ \dot{\xi}(\hat q,1)]- \frac 12 \ddot{\xi}(\hat q,1)s^2+O(s^3|\hat q|).
$$
It follows that
\begin{eqnarray*}
L(0', \hat q, s)
&=&  \int_0^1  \varphi\bigg(\eta(t);
se_n+s \dot{\xi}(\hat q,1)-\frac 12 \ddot{\xi}(\hat q,
1)s^2+O(s^3|\hat q|)\bigg)dt\\
&=&  s \int_0^1  \varphi\bigg(e_n+s[te_n-(1-t) \dot{\xi}(\hat q,1)];
  e_n+ \dot{\xi}(\hat q,1)
-\frac 12 \ddot{\xi}(\hat q,1)s\bigg)dt\\
&&
+O(s^3|\hat q|).
\end{eqnarray*}

Since
\begin{equation}
\dot \xi(\hat q,1)=\dot \xi(0',1)+O(|\hat q|)
=e_n+O(|\hat q|),
\label{6.5}
\end{equation}
and
$$
\ddot \xi(\hat q,\tau)
=\ddot \xi(0',\tau)+O(|\hat q|)=O(|\hat q|),\qquad\forall\ 0\le \tau
\le 1+\epsilon_0,
$$
we have, by properties of our special coordinates,
$$
 \varphi_{\xi^i}(e_n; e_n+ \dot{\xi}(\hat q,1))
= \varphi_{\xi^i}(e_n; 2e_n)+O(|\hat q|)=O(|\hat q|),
$$
$$
 \varphi_{\xi^iv^j}(e_n; e_n+ \dot{\xi}(\hat q,1))
=  \varphi_{\xi^iv^j}(e_n;2e_n)+O(|\hat q|)=O(|\hat q|).
$$
Making a Taylor expansion of
$\varphi$ about $(e_n, e_n+\dot{\xi}(\hat q,1))$, we have
\begin{eqnarray*}
&&L(0', \hat q, s)\\
&=& s  
 \varphi(e_n; e_n+ \dot{\xi}(\hat q,1))
-\frac 12 s^2
\varphi_{v^i}(e_n; e_n+ \dot{\xi}(\hat q,1))
\ddot{\xi^i}(\hat q,1)\\
&&+
s^2  \int_0^1  \varphi_{\xi^i}(e_n; e_n+ \dot{\xi}(\hat q,1))
[t\delta_n^i-(1-t)  \dot{\xi}^i(\hat q,1)]dt
\\
&&+ \frac 12 s^3  \int_0^1  \varphi_{\xi^i\xi^j}(e_n; e_n+\dot{\xi}(\hat q,1))
[t\delta_n^i-(1-t)  \dot{\xi}^i(\hat q,1)]
[t\delta_n^j-(1-t)  \dot{\xi}^j(\hat q,1)]dt\\
&&+O(s^3|\hat q|+s^4)\\
&=&I+II+III+IV+O(s^3|\hat q|+s^4).
\end{eqnarray*}

First
$$
III= \frac 12 s^2    \varphi_{\xi^i}(e_n; e_n+ \dot{\xi}(\hat q,1))
[\delta_n^i- \dot{\xi}^i(\hat q,1)]=O(s^2|\hat q|^2).
$$

Using (\ref{aaa2}) and (\ref{6.5}), we have
$$
\varphi_{\xi^i\xi^j}(e_n; e_n+ \dot{\xi}(\hat q,1))
[t\delta_n^i-(1-t)  \dot{\xi}^i(\hat q,1)]
[t\delta_n^j-(1-t)  \dot{\xi}^j(\hat q,1)]=O(|\hat q|)
$$
from which we deduce
$$
IV=O(s^3|\hat q|)=O(s^4+s^2|\hat q|^2).
$$

Differentiating $\varphi(\xi(\hat q,\tau), \dot{\xi}(\hat q,\tau))\equiv 1$
in $\tau$, 
we have, using $\xi(\hat q,1)=e_n$,
\begin{eqnarray}
&&\varphi_{v^i}(e_n;  \dot{\xi}(\hat q,1)) \ddot{\xi^i}(\hat q,1)\nonumber\\
&=& -\varphi_{\xi^i}(e_n;  \dot{\xi}(\hat q,1))  \dot{\xi^i}(\hat q,1)
= -\varphi_{\xi^i}(e_n;  e_n+[\dot{\xi}(\hat q,1)-e_n])  \dot{\xi^i}(\hat q,1)
\nonumber\\
&=&  -\varphi_{\xi^i}(e_n;  e_n) \dot{\xi^i}(\hat q,1)
-\varphi_{\xi^iv^j}(e_n;  e_n) \dot{\xi^i}(\hat q,1) [\dot{\xi^j}(\hat q,1)-
\delta^j_n]\nonumber\\
&&+O(|\dot{\xi}(\hat q,1)-e_n|^2
= O(|\hat q|^2). 
\nonumber
\end{eqnarray}
Since $\varphi_{v^i}(e_n; e_n+\dot
\xi(\hat q,1))=\varphi_{v^i}(e_n; \dot\xi(\hat q,1))+O(|\hat q|)$, we
conclude that
$$
II=O(s^2|\hat q|^2).
$$ 

Based on the above, we have
$$
L(0', \hat q,s)=s \varphi(e_n; e_n+\dot\xi(\hat q,1))
+O(s^4+s^2|\hat q|^2).
$$
Inserting (\ref{6.3}), we obtain (\ref{6.4}).

\vskip 5pt
\hfill $\Box$
\vskip 5pt

We now complete the proof of Proposition \ref{prop6.1}.
Combining (\ref{6.4}) and (\ref{6.1}) we obtain
$$
L(x',\hat q,s)\le 2s (1-c_0|\hat q|^2)+C(s^2|\hat q|^2+s^4)
+Cs|\hat q|^2(\frac {\bar\epsilon}K+\bar\epsilon^2 K+\frac 1K).
$$
Thus, by our conditions on $x'$, 
$$
L(x',\hat q,s)\le 2s (1-c_0|\hat q|^2)+ 
Cs|\hat q|^2
(\bar\epsilon^2 K^3\bar\delta
+K\bar\delta+\frac{\bar\epsilon}K+
\bar\epsilon^2 K+\frac 1K).
$$
Proposition \ref{prop6.1} follows, if we choose first $K$ large,
then $\bar\epsilon$ small, and, last, $\bar\delta$ small.

\section{Main Estimates IV}
\setcounter{equation}{0}

We now take up another case for which we,
again, choose $z$ of the form
$$
(x'+\bar q, f(x'+\bar q))
$$
with suitable $\bar q$.  The choice of $\bar q$ is made so
as to make
$|A|=|e_n-\xi(\bar q,1)|$ small.

Let $\zeta\in \Bbb R^{n-1}$ be a unit eigenvector of
\begin{equation}
(\tilde f_{x_\alpha x_\beta}(0')-f_{x_\alpha x_\beta}(0'))\zeta^\alpha=\lambda\zeta^\beta,
\quad 1\le\beta\le n-1,
\label{7.1}
\end{equation}
where we recall that $\lambda\ge 0$ is the smallest eigenvalue.  We set
\begin{equation}
\bar q=\rho |x'|\zeta
\label{7.2}
\end{equation}
with 
\begin{equation}
\rho\ge K^{ \frac 34}.
\label{7.3}
\end{equation}
As before $L(0', \bar q,s)$ is the
Finsler length of the segment joining $\xi(x'+\bar q, 1-s)$
to $\xi(x', 1+s)$.

\begin{prop}
For any given positive constant $\epsilon'>0$,
there exist some large
 constant $K\ge 1$ and some small constant $\delta'>0$, depending
only on $\epsilon'$, $f$ and $\varphi$
 such that
for all $\epsilon'\lambda\le |x'|\le \delta'$,
$s=K|x'|$, and $\qbar$ as above,
$$
L(x', \qbar,s)<2s.
$$
Consequently,
$$
dist\left( (x'+\bar q, f(x'+\bar q))
\add \xi(x', 1+s)\right)<1+s.
$$
\label{prop7.1}
\end{prop}

\begin{rem}  In proving the above proposition, 
$K$ will be chosen first and then $\delta'$.
\label{rem7.1}
\end{rem}

We first establish
\begin{lem} For 
some positive constants $c_0$, $C$, 
$K$ and $\delta'$,  depending only on $\varphi$ and $f$, 
we have, for $x'$, $\qbar$ and $s$ above, that
$$
L(0',\qbar,s)\le 2s(1-c_0|\qbar|^2).
$$
\label{lem7.1}
\end{lem}

\noindent {\bf Proof of Lemma \ref{lem7.1}.}\
Let $\widetilde \xi=\widetilde \xi(x', \tau)$, $\tau\ge 0$,
denote the geodesics satisfying
$$
\varphi(\widetilde \xi; \dot{\widetilde \xi})\equiv 1,
$$
$$
\widetilde \xi(x',0)=(x', \widetilde f(x')),
$$
$$
\dot {\widetilde \xi}(x';0)= \widetilde V(x'),
$$
where $ \widetilde V(x')$ is defined as $V(x')$ in Section 2,
but for $\widetilde f$ instead of for $f$.

By the property of $\widetilde f$, 
$$
\widetilde \xi(x', 1)=e_n, \qquad |x'|<\tilde \epsilon_1.
$$

For any $q\in \Bbb R^{n-1}$, $|q|$ small, let
\begin{eqnarray*}
\widetilde\eta(x',q,s;t)&=&
(1-t)\widetilde\xi(x'+q, 1-s)+t \widetilde\xi(x',1+s)\\
&=&
\widetilde\xi(x'+q, 1-s)+t[\widetilde\xi(x',1+s)-\widetilde\xi(x'+q, 1-s)],
\end{eqnarray*}
and 
$$
\widetilde L(x', q, s)
=\int_0^1 \varphi(\widetilde \eta(x',q,s;t);
\dot{\widetilde \eta}(x',q,s;t))dt.
$$

By Lemma \ref{lem6.1}, applied to $\widetilde f$ with $\hat q=\qbar$,
we have
\begin{equation}
\widetilde L(0', \qbar, s)\le 2s(1-c_0|\qbar|^2)
+ O(s^4+s^2|\bar q|^2).
\label{7.4}
\end{equation}

In the rest of the proof we mainly
estimate $|L(0', \qbar, s)-\widetilde L(0', \qbar, s)|$.

Clearly,  by our choice of the vector
$\zeta$, 
$$
|\widetilde\xi(\qbar, 0)-\xi(\qbar, 0)|
=\tilde f(\qbar)- f(\qbar)\le C(\lambda|\qbar|^2+|\qbar|^3),
$$
and
$$
|\dot{\widetilde\xi}(\qbar, 0)-\dot\xi(\qbar, 0)|\le C(\lambda|\qbar|
+|\qbar|^2).
$$
Since both $\xi(\qbar,\cdot)$ and $\widetilde \xi(\qbar,\cdot)$
satisfy the same ODE, we have
\begin{equation}
|\widetilde\xi(\qbar, t)-\xi(\qbar, t)|
+|\dot {\widetilde\xi}(\qbar, t)- \dot \xi(\qbar, t)|
\le C(\lambda|\qbar|
+|\qbar|^2),\qquad \forall\ t.
\label{7.5}
\end{equation}

Next, one verifies that
\begin{equation}
L(0', \qbar, s)=s\int_0^1\varphi\bigg( \widetilde \eta(t)
-(1-t)(\widetilde \xi-\xi)(\qbar, 1-s);
\frac 1s \dot{\widetilde \eta}(t)+\frac 1s 
(\widetilde \xi-\xi)(\qbar, 1-s)\bigg)dt,
\label{7.6}
\end{equation}
where
$$
\tilde \eta(t):=\tilde \eta(0', \qbar, s; t).
$$
By (\ref{7.5}), for $0\le t\le 1+\epsilon_0$,  
$$
(|\widetilde \xi-\xi|
+|\dot {\widetilde \xi}-\dot \xi|)
(\qbar, t)=O(\lambda|\qbar|+|\qbar|^2).
$$
We also have
$$
\widetilde \eta(t)
=[1+(2t-1)s]e_n+O(|\qbar|),
$$
$$
\frac 1s \dot{\widetilde \eta}(t) 
=2e_n+O(|\qbar|).
$$
The last equality above needs some explanation:  By Taylor expansion,
\begin{eqnarray*}
\dot{\widetilde \eta}(t) &=& \widetilde \xi(0', 1+s)-
 \widetilde \xi(\qbar, 1-s)\\
&=& (1+s)e_n-\widetilde \xi(\qbar, 1)+
\dot {\widetilde \xi}(\qbar, 1)s-
\frac 12 \ddot{\widetilde \xi}(\qbar, 1-\theta s)s^2,
\end{eqnarray*}
where $0\le \theta\le 1$.
Since $\widetilde \xi(\qbar, 1)=e_n$,
$\dot{\widetilde \xi}(0', 1)=e_n$ and
$\ddot {\widetilde \xi}(0', t)=0$ for all $0\le t\le 1+\epsilon_0$,
we have 
$\dot{\widetilde \xi}(\qbar, 1)=e_n+O(|\qbar|)$,
$ \ddot{\widetilde \xi}(\qbar, 1-\theta s)=O(|\qbar|)$, and therefore
$$
\dot{\tilde \eta}(t)=2se_n+O(s|\qbar|).
$$

It is clear that
$$
\xi^n(0', t)-t\equiv\xi^n_{x_\alpha  }(0', t)\equiv 0,
\quad  \widetilde \xi^n(0', t)-t\equiv
 \widetilde \xi^n_{x_\alpha  }(0', t)\equiv 0,
\quad 0\le t\le 1+\epsilon_0.
$$

It follows that
\begin{eqnarray*}
\widetilde \xi^n(\qbar, 1-s)&=&
\widetilde \xi^n(0', 1-s)+ \widetilde \xi^n_{x_\alpha  }(0', 1-s)\qbar_\alpha
+
\int_0^1 \int_0^t 
 \widetilde \xi^n_{x_\alpha  x_\beta  }(\tau\qbar, 1-s)
\qbar_\alpha\qbar_\beta d\tau dt
\\
&=& (1-s)
+\int_0^1  \int_0^t 
 \widetilde \xi^n_{x_\alpha  x_\beta  }(\tau\qbar, 1-s)
\qbar_\alpha\qbar_\beta d\tau dt.
\end{eqnarray*}
Similarly
$$
 \xi^n(\qbar, 1-s)=
(1-s) +\int_0^1  \int_0^t \xi^n_{x_\alpha  x_\beta  }(\tau\qbar, 1-s)
\qbar_\alpha\qbar_\beta d\tau dt.
$$

By (\ref{43-2}), applied to both $\xi$ and $\tilde \xi$,
 we deduce from the above that
\begin{eqnarray*}
&&
(\widetilde  \xi^n- \xi^n)(\qbar, 1-s)
\\
&=&\frac  12\varphi_{v^iv^j}((1-s)e_n; e_n)
(\dot \xi_{x_\alpha  }^i \xi_{x_\beta  }^j
-\dot {\widetilde \xi^i}_{x_\alpha  } \widetilde \xi_{x_\beta  }^j)
(0', 1-s)\qbar_\alpha\qbar_\beta+O(|\qbar|^3).
\end{eqnarray*}

Thus, by  (\ref{7.5}), we have
\begin{equation}
|(\widetilde \xi^n-\xi^n)(\qbar, 1-s)|\le C(\lambda
|\qbar|+|\qbar|^2) |\qbar|^2+C|\qbar|^3\le C|\qbar|^3.
\label{7.7}
\end{equation}

Estimate (\ref{7.7})
 will be used  below.

By Taylor expansion in (\ref{7.6}), we have, using (\ref{7.5}), 
\begin{eqnarray*}
&&L(0', \qbar, s)\\
&=& s
\int_0^1\varphi(\widetilde \eta, \frac 1s \dot{\widetilde \eta})dt
-s\int_0^1\varphi_{\xi^i}(\widetilde \eta, \frac 1s
\dot{\widetilde \eta})(1-t)(\widetilde \xi-\xi)^i(\qbar, 1-s)dt\\
&&+\int_0^1\varphi_{v^i}(\widetilde \eta, \frac 1s\dot{\widetilde \eta})
(\widetilde \xi-\xi)^i(\qbar,1-s)]dt\\
&&+O(\lambda|\qbar|^2+|\qbar|^3)
+O\left(  \frac{  (\lambda|\qbar|+|\qbar|^2)^2  } s \right).
\end{eqnarray*}
By the properties of special coordinates and the
expressions of $\tilde \eta$ and $\frac 1s \dot{\tilde \eta}$, 
$$
(|\varphi_{\xi^i}|+|\varphi_{v^\alpha}|)(\widetilde \eta; \frac 1s
\dot{\widetilde \eta})
\le C|\qbar|,\quad \forall\
1\le \alpha\le n-1,
$$
$$
|\varphi_{ v^n}(\widetilde \eta; \frac 1s
\dot{\widetilde \eta})|\le C.
$$
It follows that
$$
L(0', \qbar, s)
= \widetilde L(0',\qbar, s)+ 
O\left(  (\lambda|\qbar|^2+|\qbar|^3)(1+ \frac {\lambda+|\qbar|}s)\right).
$$
Combining this with (\ref{7.7}) and (\ref{7.4}), we find
\begin{eqnarray*}
L(0', \bar q,s)&\le&
2s(1-c_0|\bar q|^2)+
C(s^4+s^2|\bar q|^2)+
C(\lambda|\bar q|^2+|\bar q|^3)(1+ \frac {\lambda+|\qbar|}s)\\
&\le & 2s(1-c_0|\bar q|^2)+Cs|\bar q|^2(\frac{s^3+s|\qbar|^2}{ |\qbar|^2})
+ Cs|\bar q|^2( \frac {\lambda+|\qbar|}s )
(1+ \frac {\lambda+|\qbar|}s ).
\end{eqnarray*}
Since
$$
\frac{s^3+s|\qbar|^2}{ |\qbar|^2}=
\frac {K^3}{\rho^2}|x'|+K|x'|\le(K^{\frac 32}+K)\delta',
$$
and
$$
\frac  {\lambda+|\qbar|}s\le
\frac 1{K\epsilon'}+\frac \rho K=
\frac 1{K\epsilon'}+\frac 1{ K^{\frac 14} },
$$
we obtain the desired estimate by choosing  first $K$ large and then $\delta'$ small
(recall that
 we also want $K\delta'<\bar \epsilon_1$ in Proposition
\ref{prop5.1}).   Lemma \ref{lem7.1} is proved.

\vskip 5pt
\hfill $\Box$
\vskip 5pt

\noindent{\bf Proof of Proposition \ref{7.1}.}\
We make use of Proposition \ref{prop5.1}, and for this we need an estimate of
$$
A=e_n-\xi(\bar q,1).
$$
In fact, since $\tilde\xi(\bar q,1)=e_n$ we have from
(\ref{7.5}),
$$
|A|\le C(\lambda|\bar q|+|\bar q|^2)\le
C|\bar q|(\frac { |x'| }{ \epsilon'}+|\bar q|).
$$
We have to verify (\ref{5.5}).
Well, $|\bar q|<1$ since
$K\delta'\le \bar \epsilon_1<1$, so
\begin{eqnarray*}
\frac { |A|+|x'| }s &\le &
C \frac {|x'|}s(1+\frac { |\bar q|}{\epsilon'})+C\frac {|\bar q|}s
\le C(\frac 1K+\frac 1{ \epsilon' K }+\frac 1{ K^{\frac 14}})\\
&<&\epsilon_3\qquad\mbox{of Proposition \ref{prop5.1}}
\end{eqnarray*}
provided we increase $K$ still further, which means decreasing
$\delta'$. 

We may thus apply Proposition \ref{prop5.1}
and conclude that
$$
L(x', q,s)-L(0',q,s) 
\le Cs|\bar q|^2(\frac 1{K\rho}+\frac 1{\rho^2}
+\frac \rho{K^2}+\frac 1{K\rho\epsilon'}+\frac 1K+
\frac 1{K^2\epsilon'}+\frac K{\rho^2}).
$$
Recalling that $\rho=K^{\frac 34}$
 and combining
the above with Lemma \ref{lem7.1} we 
obtain the desired result again, if necessary,
by increasing $K$ and decreasing $\delta'$.

\vskip 5pt
\hfill $\Box$
\vskip 5pt

\section{ Proof of Theorem \ref{thm1.1}}
\setcounter{equation}{0}

We consider $\Omega$ bounded.  The proof for
unbounded $\Omega$ goes the same.
Following the notations in the introduction, we need
to prove that $\bar s(y)$ is a Lipschitz function
 on $\partial \Omega$.
Namely, we need to show that there  exist some positive  constants $K$ 
and $\delta$ such that for any $\bar y\in 
\partial \Omega$, 
\begin{equation}
\bar s(y)\le \bar s(\bar y)+K|y-\bar y|,
\qquad\forall\ y\in \partial\Omega,\
|y-\bar y|\le \delta.
\label{want1}
\end{equation}

As before, by making
 a change of variables, we may assume without loss of generality
that $\bar y=0\in \partial \Omega$, $\bar s(\bar y)=1$,
and, for some $\epsilon_0>0$, $\xi(\bar y, t)= te_n$
for all $ -\epsilon_0\le t\le 1+\epsilon_0$.
By our result in Section 3 on the existence of
special coordinates, we may also assume, for all
$-\epsilon_0\le t\le 1+\epsilon_0$,
(\ref{ag5})-(\ref{aaa2}) hold.

We may assume that for some $\epsilon_1>0$,
 $f(x')$ is a  $C^{2,1}$ function  defined
in $|x'|<\epsilon_1$, $x'\in \Bbb R^{n-1}$, $f(0')=0$, $\nabla f(0')=0$,
  and
$\{(x', f(x'))\ |\ |x'|<\epsilon_1\}$ is a local representation of
$\partial \Omega$.  In the following, as before, we  use
$\xi(x', t)$ to denote $\xi((x', f(x')), t)$.  With this notation, we have
$$
\xi(x', 0)=(x', f(x')),
$$ 
and, by 
 Lemma \ref{lem2.2},
$$
\dot \xi(x', 0)=V(x'),
$$
where $V(x')$ is the vector field given in Section 2.

To prove (\ref{want1}), we only need to show that for
some constants $K$ and $\delta$, depending only on 
$\partial \Omega$ and $\varphi$,
we have
\begin{equation}
dist(\partial \Omega, 
\xi(x', 1+K|x'|))
 <1+K|x'|,\qquad\forall\
|x'|<\delta.
\label{want2}
\end{equation}

We put together the results of Sections 4-7.

\noindent{\bf Proof of Theorem \ref{thm1.1}.}\
We distinguish two cases.

\noindent{\bf Case 1.}\ There exists some 
$Q\in \partial \Omega\setminus\{0\}$
such that
$$
dist(0\add e_n)=dist(Q\add e_n).
$$

\medskip

\noindent{\bf Case 2.}\ For all $y\in \partial \Omega\setminus\{0\}$,
 we have
$$
dist(0\add e_n)<dist(y\add e_n).
$$

\medskip

In Case 1, we may assume, because of Proposition
\ref{prop4-0},  that
$Q=(\hat q, f(\hat q))$ for some $\hat q\in \Bbb R^{n-1}$
satisfying $|\hat q|\le \bar \epsilon_1/9$, where
$\bar\epsilon_1$ is that of Proposition \ref{prop5.1}.

Since
$$
dist(0\add e_n)=dist((x',\tilde f(x'))\add e_n),
\qquad\forall\ |x'|\le \tilde\epsilon_1,
$$
we have
$$
\widetilde f(\hat q)=f(\hat q),
$$
$$
\widetilde f(x')\ge f(x'),
\qquad\forall\ |x'|\le \tilde\epsilon_1,
$$
and
$$
(\widetilde f-f)(x')=\frac 12 (\widetilde f_{x_\alpha   x_\beta  }-
f_{x_\alpha   x_\beta  })(0')x_\alpha   x_\beta  +O(|x'|^3).
$$
Recall that $\lambda\ge 0$ is the least eigenvalue of
$\big( (\widetilde f_{x_\alpha   x_\beta  }-
f_{x_\alpha    x_\beta   })(0') \big)$. 

We thus have
$$
0=(\widetilde f-f)(\hat q)
\ge \frac 12\lambda |\hat q|^2+O(|\hat q|^3),
$$
and therefore
$$
\lambda\le C|\hat q|,
$$
where $C$ depends only on the $C^{2,1}$ norm of
$f$ and $\widetilde f$.

Let $\bar \epsilon$, $\bar \delta$ and $K$
 be the positive constants in Proposition \ref{prop6.1}, 
then, by Proposition \ref{prop6.1},
$$
dist\left((x'+\hat q, f(x'+\hat q))\add
\xi(x', 1+K|x'|)\right)<1+K|x'|,\qquad
\forall\ |x'|\le \min\{\bar\delta, \bar\epsilon|\hat q|\}.
$$

For $|x'|\ge \bar \epsilon |\hat q|$, we have, by the above,
$$
|x'|\ge \frac{\bar\epsilon}{C}\lambda.
$$
Let $\epsilon'=\frac{\bar\epsilon}{C}$, we have, by
Proposition \ref{prop7.1}, for some $K$ and $\delta'$ depending on $\bar \epsilon$,
$$
dist(\xi(x', 1+K|x'|)\add \qbar)<1+K|x'|,\qquad
\forall\ \epsilon'\lambda\le |x'|\le \delta'.
$$
Thus we have established (\ref{want2}) for
 $\delta=\min\{\bar\delta, \delta'\}$ and some positive constant $K$.

\medskip

In Case 2,
we have, by Lemma \ref{lem4.5},  $\lambda=0$.
The desired estimate (\ref{want2}) 
then follows from Proposition \ref{prop7.1}.

\vskip 5pt
\hfill $\Box$
\vskip 5pt

\section{}  
\setcounter{equation}{0}

In this section we consider the general
case
\begin{equation} 
H(x,u,\nabla u)=1\qquad\mbox{in}\ \Omega, 
\label{0.1} 
\end{equation} 
a bounded domain in $\Bbb R^n$,with $\partial \Omega$ in 
$C^{2,1}$.   In 
Theorem \ref{thm9.1}, 
under certain strict  conditions, we find a
positive viscosity solution $u$  
 satisfying 
\begin{equation} 
u=0\qquad \mbox{on}\ \partial \Omega 
\label{0.2} 
\end{equation} 
and show that for its singular set $\Sigma$,  
\begin{equation} 
H^{n-1}(\Sigma)<\infty .
\label{0.3} 
\end{equation} 
We then derive Propositions \ref{prop1.1}-\ref{prop1.3} of 
Section 1.4.

We shall make two conditions.
The first is Situation (*) of Section 1.4
which we repeat here as
 
\noindent{\bf Assumption I.}\ {\it
The function $H(x,t,p)$, $t\in \Bbb R$,
$p\in \Bbb R^n$, 
is assumed to satisfy:
For every $x$ in $\overline \Omega$ the set 
$$
V_x=\{(t,p)\ |\ 
H(x,t,p)<1\} 
 $$
is a convex set in $\Bbb R\times \Bbb R^{n+1}$ lying in
a fixed downward cone
\begin{equation}
|p|\le k(C_1-t),
\qquad t<C_1\ \mbox{with}\ k, C_1>0.
\label{9.1}
\end{equation}
Thus $t$ may be unbounded below in $V_x$.  The boundary of $V_x$, 
$$ 
S_x=\{(t,p)\ |\ 
H(x,t,p)=1\} 
$$
is assumed to be a smooth strictly convex
hypersurface in $(t,p)$ space
with positive principal
curvature for $t$ in the region
\begin{equation}
-1\le t\le C_1,
\label{6.6}
\end{equation}
uniformly for $x$ in $\overline \Omega$.
Furthermore, the origin
in $\Bbb R\times \Bbb R^n$ lies in $V_x$
and is bounded away from $S_x$ by some number
$$
r_0>0.            
$$
In addition, $H(x,t,p)$ is smooth in a neighborhood
of $\displaystyle{ \cup _x S_x }$.
}

\medskip

See Fig. 7  for example.

\bigskip

\centerline{
\epsfbox{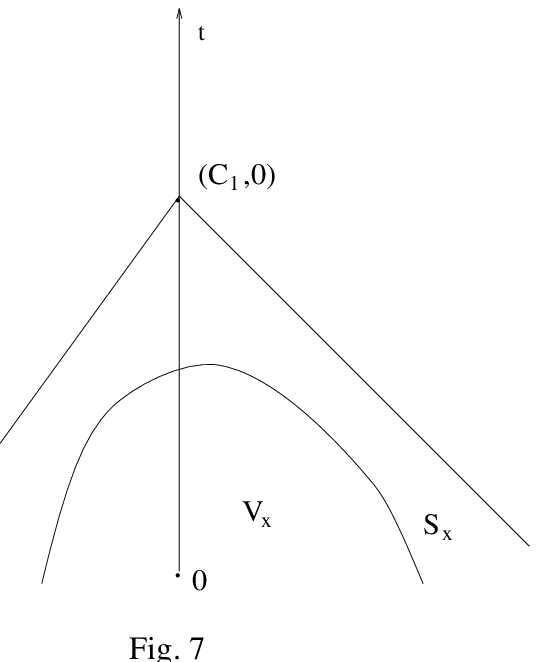}}

\bigskip

Thus a common assumption that $H$ is monotone in $t$ does not 
necessarily hold here.  Under another
assumption on $H$ we construct a  
viscosity solution; it need not, however, be unique. 
As we have said in Section 1.2, 
the function $H(x,t,p)$
 is not so important, the important things are the sets $V_x$. 
 
Because of Remark \ref{rem1.3}
 we may take 
$H$ to be homogeneous of degree $1$ in 
$(t,p)$.  It is thus completely determined by the $S_x$.  From 
now on we assume this homogeneity. 
 
Our way  of studying the problem is to set up a 
related problem in one higher dimension--in $\Bbb R\times \overline \Omega$: 
 
For $\tau\in \Bbb R$, given a function $u$ in 
$\overline \Omega$, we 
define $z$ in $\Bbb R\times \overline \Omega$ by 
$$ 
z(\tau, x)=e^\tau u(x). 
$$ 
Multiplying the equation (\ref{0.1}) by $e^\tau$--recall the homogeneity 
condition--we obtain 
$$ 
H(x, e^\tau u, e^\tau \nabla u)=e^\tau, 
$$ 
which we rewrite as 
\begin{equation} 
e^{-\tau}H(x, z_\tau, \nabla_x z)=1. 
\label{0.4} 
\end{equation} 
 
This is a Hamilton-Jacobi equation in  
$\Bbb R\times \Omega$, for $z$, and we solve 
it under the boundary condition 
\begin{equation} 
z=0\qquad\mbox{on}\ \Bbb R\times \partial \Omega. 
\label{0.5} 
\end{equation} 
 
As in Section 1 the solution involves 
the support function 
\begin{eqnarray*} 
\tilde \varphi(x,\tau; s,v)&=& \sup_{ e^{-\tau}H(x,t,p)=1 } 
(st+v\cdot p) 
= e^{\tau} \sup_{H(x,t,p)=1 } 
(st+v\cdot p)\\ 
&=& e^\tau \varphi(x; s,v) 
\end{eqnarray*} 
where $\varphi$ is the support function of $S_x$: 
$$
\varphi(x; s,v)=\sup_{ H(x,t,p)=1 }(st+v\cdot p). 
$$
 
According to Section 1 which uses formula 
(55)$'$ on page 132 of \cite{L}, the viscosity solution 
of (\ref{0.4}), (\ref{0.5}) is obtained using curves 
$(w(t), \xi(t))$, $0\le t\le T$ lying 
in $\Bbb R\times \overline \Omega$ with 
\begin{equation} 
w(0)=\tau, \ \xi(0)=x;\ \ 
w(T)=\mu, \ \xi(T)=y\in \partial \Omega. 
\label{0.7} 
\end{equation} 
 
The solution is given by 
\begin{equation} 
z(\tau,x)=\inf_{ \mu\in \Bbb R, y\in \partial \Omega} 
\ \ \inf_{  (w,\xi) } 
\int_0^T e^{ w(t) }\varphi(\xi(t); -\dot w, -\dot\xi)dt. 
\label{0.8} 
\end{equation} 
 
Here $\displaystyle{ \inf_{w, \xi} }$ means infimum over curves satisfying 
(\ref{0.7}). By Remark 5.5 in \cite{L},
$z$ is a viscosity solution even though
$\Bbb R\times \Omega$ is unbounded. 
 
Note that 
\begin{equation} 
z(\tau,x)=e^\tau z(0,x). 
\label{0.9} 
\end{equation} 
This follows from 
\begin{rem} If $(w(t), \xi(t))$ is an eligible curve in 
(\ref{0.8}) for $z(\tau,x)$ then
\newline  
$(w(t)-\tau, \xi(t))$ is one for $z(0,x)$. 
\label{remX}
\end{rem} 
 
We are really only interested in 
\begin{equation} 
u(x):=z(0,x), 
\label{0.10} 
\end{equation} 
because of 
\begin{clm} Since $z(\tau,x)$ is a viscosity solution of 
(\ref{0.4}), (\ref{0.5}),  $u(x)$, given 
by (\ref{0.10}), is a viscosity solution of 
(\ref{0.1}), (\ref{0.2}). 
\label{claim10.1}
\end{clm} 
This is easily seen.  For instance, to check that $u$ is a viscosity subsolution we have to 
show that for any $\varphi\in C^1(\Omega)$ such that 
$u-\varphi$ has a local maximum$=0$
 at some point $x_0\in \Omega$, necessarily, 
\begin{equation} 
H(x_0, u(x_0), \nabla \varphi(x_0))\le 1. 
\label{0.11} 
\end{equation} 
To see this for such a $\varphi$, consider 
$$ 
\tilde\varphi(\tau,x)=e^\tau \varphi(x). 
$$ 
Because of (\ref{0.9}), $z-\tilde\varphi$ has a local maximum$=0$  
at $(0,x_0)$ and since $z$ is a viscosity subsolution 
of (\ref{0.4}), 
$$ 
H(x_0, \tilde \varphi_\tau(0,x_0), \nabla_x\tilde\varphi(0,x_0))\le 1, 
$$ 
i.e. (\ref{0.11}) holds. 
 
Turning now to the singular sets of 
$z$ and $u$, we see from (\ref{0.9}) that the singular
set $\widetilde \Sigma$ of $z$ is a straight cylinder with 
generators parallel to the $\tau-$axis lying over
the singular set $\Sigma$ of $u$.
Thus if we know that
\begin{equation}
H^n(\widetilde \Sigma)<\infty,
\label{0.12}
\end{equation}
it follows that
$$
H^{n-1}(\Sigma)<\infty
$$
----our desired conclusion (\ref{0.3}).

Indeed our main result,  Theorem \ref{thm1.1} of
Section 1,  yields exactly (\ref{0.12}).

$\underline{
\mbox{\it Wrong}}$.
  We have to be more careful: The domain
$\Bbb R\times \Omega$ is not bounded and  we cannot
 apply our Lipschitz continuity result of Theorem \ref{thm1.1} in
Section 1; it holds for {\it compact
subsets } of the boundary.

 We are thus 
led to add a further restriction on the sets
$V_x$ relative to the domain $\Omega$:

Set
$$
\overline C=\sup_{x\in \Omega}\  \inf_{ y\in\partial \Omega}
\ \  \inf_{ \xi, \xi(0)=x, \xi(T)=y}
\int_0^T \varphi(\xi(t); 0, -\dot \xi(t))dt.
$$
This is the shortest distance from $x$ in 
$\Omega$ to $\partial \Omega$ in the restricted Finsler metric
$\varphi(\xi(t); 0, -\dot\xi(t))dt$.

Next, consider the support function
$\varphi(x; s,v)$.  From its definition, we have
\begin{equation}
c_0(|s|+|v|)\le \varphi(x; s,v)\le C_0(|s|+|v|)
\label{0.14}
\end{equation}
for suitable positive constants $c_0$ and $C_0$.

Set
$$
\sigma:=\sup_{ \varphi(x; s,v)=1, x\in \overline \Omega} s.
$$
The additional condition we impose is

\medskip

\noindent{\bf Assumption II.}\ $\sigma\overline C<1$.

\medskip

Assumption II  may be expressed more directly in terms
of the sets $S_x$:  For any $x\in \overline
\Omega$ denote by $\bar t=\bar t(x)$ 
the point $(\bar t,0)$ on $S_x$
with $\bar t>0$.
Since for every $x$, $H(x,t,p)$ is the support 
function of the convex hypersurface
$$
\hat S_x=\{(s,v)\ |\
\varphi(x; s,v)\equiv 1\},
$$
it follows that
$$
1=H(x,\bar t,0)=\sup_{ \varphi(x;s,v)=1} s\bar t
$$
is achieved at a point where $s=\bar s$,  the maximum value
of $s$ on $\hat S_x$.  Thus $\bar t=1/{\bar s}$ and
so
$$
\frac 1\sigma =\min_{  x\in \overline \Omega} \bar t(x).
$$
Hence Assumption II is equivalent to the condition
$$
\min_{ x\in \overline \Omega}\bar t(x)>\overline C.
$$
We now state the main result of this section.
\begin{thm} Under Assumption I
and Assumption II, the problem 
(\ref{0.1}), (\ref{0.2}) possesses a
positive viscosity solution and its singular set
$\Sigma$ satisfies
$$
H^{n-1}(\Sigma)<\infty.
$$
\label{thm9.1}
\end{thm}

The proof of Theorem \ref{thm9.1} is based on the following
lemma--we assume  Assumption I
and Assumption II.

For $x\in \Omega$ fixed and $0<\epsilon$ fixed, so that
\begin{equation}
\sigma(\overline C+\epsilon)<\frac {1+\sigma\overline C}2,
\label{0.17}
\end{equation}
consider a competing curve $(w(t), \xi(t))$,
$0\le t\le T$, such that 
$w(0)=0$, $\xi(0)=x$,
$w(T)=\mu$, $\xi(T)=y\in \partial \Omega$ and such that
\begin{equation}
\int_0^T  e^{ w(t) }\varphi(\xi(t);
-\dot w(t), - \dot \xi(t))dt<\overline C+\epsilon.
\label{0.18}
\end{equation}
By our definition of $\overline C$, such 
a curve exists.  Let us normalize the parameter $t$ so that
\begin{equation}
\varphi(\xi(t); -\dot w(t),  -\dot \xi(t))\equiv 1.
\label{0.19}
\end{equation}
$T$ is of course unknown.

\begin{lem} In the situation above,
\begin{equation}
T, |w(t)|\le C(\overline C, \sigma).
\label{0.20}
\end{equation}
\label{lem-XX}
\end{lem}

\noindent {\bf Proof.}\  By the definition of $\sigma$, because
of (\ref{0.19}), 
$$
-\dot w(t)\le \sigma.
$$
Thus
$$
w(t)\ge -\sigma t
$$
and inserting this in (\ref{0.18}), we obtain
$$
\overline C+\epsilon>\int_0^T e^{ -\sigma t}dt=\frac 1\sigma (1-e^{-\sigma T}).
$$
Using (\ref{0.17}), we find
$$
e^{ -\sigma T}> \frac {1-\sigma\overline C}2,
$$
from which a bound for $T$, as in
(\ref{0.20}), follows.

By (\ref{0.14})
$$
|\dot w(t)|\le \frac 1{ c_0}
$$
and thus
$$
|w(t)|\le \frac T{c_0},
$$
completing (\ref{0.20}).

\vskip 5pt
\hfill $\Box$
\vskip 5pt

We have proved that if we consider
competing curves for $z(0,x)$ with
``lengths'' close to $z(0,x)$ then, on them,
\begin{equation}
|w|\le C_1, \qquad \mbox{uniform in}\ x.
\label{0..21}
\end{equation}

By Remark \ref{remX} it follows that for $|\tau|\le 1$
if we consider competing curves for $z(\tau,x)$
in (\ref{0.8}) with lengths
sufficiently close to $z(\tau,x)$ then
on these
\begin{equation}
|w|\le C_2, \qquad \mbox{uniform in}\ x.
\label{0.22}
\end{equation}

We are now in a position to give the 

\noindent {\bf Proof of Theorem \ref{thm9.1}.}\
We change $S_x$ to $\tilde S_x$
by making it bounded from below.
We can do so with Assumption II unchanged--for $\tilde S_x$.
This can be done by taking $\delta>0$ very small, and changing
 $S_x$ smoothly
(also in $x$) so that it is unchanged for $t\ge -\delta$
but does not extend below $t=-2\delta$.  It is clear that
Assumption II still holds if we take
$\delta>0$ small and change $S_x$ properly.
By Remark \ref{rem1.3} 
we may assume that $S_x$ satisfies this additional property.

Finally, we construct a bounded domain $D$ in $\Bbb R^{n+1}$, over
$\Omega$, with $C^{2,1}$ boundary, which agrees with the cylinder when 
$|\tau|<2C_2$, as pictured.

\bigskip

\centerline{
\epsfbox{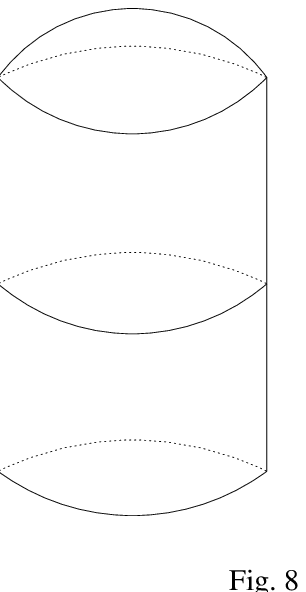}}

\bigskip

In D we solve (\ref{0.4}), (\ref{0.5})
by the formula (\ref{0.8}) where the curves $(w(t), \xi(t))$ go from
$(\tau,x)$ to the boundary of $D$--obtaining function
$z$.

  As we indicated
previously, $u(x)=z(0,x)$
is then a viscosity solution of
(\ref{0.1}), (\ref{0.2}).
Applying our main result, Theorem \ref{thm1.1}, to
$z$ in $D$, we see that the singular set
$\tilde \Sigma$ of $z$ has
$$
H^n(\tilde \Sigma)<\infty.
$$
Now (\ref{0.9}) holds for $|\tau|\le 1$ and
hence, for $|\tau|\le 1$, the singular set of $z$ is a finite cylinder over
the singular set $\Sigma$ of $u=z(0,x)$.  Consequently,
$$
H^{n-1}(\Sigma)<\infty,
$$
and we are through.

\begin{conj} Theorem \ref{thm9.1} holds merely under
Assumption I.  
\end{conj}

\noindent{\bf Proofs of Proposition \ref{prop1.1} 
and \ref{prop1.2}.}\ From the conditions in these propositions it is clear that Assumption II is satisfied.  Thus Theorem \ref{thm9.1} applies, proving
the propositions.

\vskip 5pt
\hfill $\Box$
\vskip 5pt

\noindent{\bf Proof of Proposition \ref{prop1.3}.}\
For $d_0$ small we verify Assumption II by showing that $\bar C$
is small.  Namely, from any point $x\in\Omega'$ we
join it to $y$ on $\partial \Omega'$ minimizing
$|y-x|$ by a straight segment
$$
\xi(t)=x+t(y-x)\qquad 0\le t\le 1.
$$
Then its Finsler length from $y$ to $x$ is
$$
\int_0^1 \varphi(\xi(t);0, x-y)dt\le C_0d_0
$$
by (\ref{0.14}).  Hence
$$
\bar C\le C_0d_0;
$$
it follows that for $d_0$ small, depending only on
$H$, Assumption II holds, and Theorem \ref{thm9.1}
applies.

\vskip 5pt
\hfill $\Box$
\vskip 5pt

\noindent{\bf Proof of Proposition \ref{prop1.4}.}\
As usual, we may suppose that the set
$$
V=\{(t,p)\ |\ H(t,p)=1\}
$$
is bounded and satisfies Assumption I as in the proof of 
Theorem \ref{thm9.1}, and that $H$
is positive homogeneous of degree one.
As in the proof of Theorem \ref{thm9.1} we consider the H-J equation
(\ref{0.4}) involving the extra
variable $\tau$:
$$
e^{-\tau}H(z_\tau, \nabla_x z)=1,
$$
and
consider the solution given by (\ref{0.8}).

First we obtain a bound on 
$$
u(x)=z(0,x).
$$
To this end we consider a competing
curve of the form
$$
w(t)=-\lambda t,\quad
\xi(t)=x-\lambda t V, \qquad 0\le t\le T
$$
where $V$ is a constant vector in $\Bbb R^n$ and
$$
T=\frac {d_V(x) }{ \lambda |V|}.
$$
Here $d_V(x)$ is the length of the segment
from $x$ in the direction $V$ until it hits $\partial\Omega$.
The curve is an eligible one and its length
$$
L=\int_0^T e^{ -\lambda t}
\varphi(\lambda, \lambda V)=\varphi(1, V)
(1-e^{ -\lambda T}).
$$
We now choose $V$ so as to minimize
$\varphi(1, V)$.

Letting
$$
\sigma=\max_{ \varphi(s,v)=1}s,
$$
it's clear that
$$
\varphi(\sigma, V)\ge 1\qquad \forall\ V
$$
and
$$
\min_{ V}\varphi(\sigma, V)=1.
$$
So
$$
\min_{ V}\varphi(1, V)=\frac 1\sigma.
$$
Now fix $V$ so that
$$
\varphi(1, V)=\frac 1\sigma.
$$
Since $\bar t<\hat t$, $V\ne 0$.  Recall that $\sigma=\displaystyle{
\frac 1{\bar t} }$.  Thus
$$
L=\bar t(1-e^{ -\frac{ d_V(x)}{ |V| } })=\bar t-a,\qquad
a>0,
$$
and hence
$$
u(x)=z(0,x)\le \bar t-a.
$$

We now follow the proof of Theorem \ref{thm9.1}.
Consider a competing curve $(w(t), \xi(t))$,
$0\le t\le T$, satisfying $w(0)=0$,
$\xi(0)=x$, $w(T)=\mu$,
$\xi(T)=y\in\partial \Omega$, and such that
\begin{equation}
\int_0^T e^{ w(t) }\varphi(-\dot w; -\dot\xi)dt\le
\bar t-\frac a2.
\label{9.21}
\end{equation}
As usual, we normalize the parameter $t$ so that
$$
\varphi(-\dot w; -\dot\xi)\equiv 1.
$$

\begin{lem} In the situation above,
\begin{equation}
T, |w(t)|\le C\ \mbox{independent of}\ x.
\label{9.22}
\end{equation}
\label{lem9.1prime}
\end{lem}

\noindent{\bf Proof.}\ It is the same as that of Lemma
\ref{lem-XX}.  Namely, we have
$$
-\dot w\le \sigma=\frac 1{\bar t}.
$$
Thus
$$
w\ge -\sigma t.
$$
Inserting this in (\ref{9.21}) we find
$$
\bar t-\frac a2\ge \int_0^T e^{w(t)}dt 
\ge  \int_0^T e^{-\sigma t}dt=\frac 1\sigma (1-e^{ -\sigma T})
$$
i.e.
$$
e^{ -\sigma T}\ge \frac {\bar t a}2.
$$
The bound for $T$ in (\ref{9.22})
follows.
Then, as before, we have
$ |\dot w(t)|\le \frac 1{c_0}$, so
$|w(t)|\le \frac T{c_0}$.
Lemma \ref{lem9.1prime} is proved.

The proof of Proposition \ref{prop1.4} then proceeds
as in the proof of Theorem \ref{thm9.1}.

\vskip 5pt
\hfill $\Box$
\vskip 5pt

The assumption $\bar t<\hat t$ in Proposition \ref{prop1.4}
seems strange.  However, in case
$\bar t=\hat t$, our method of proof
must fail.  Indeed, if we take
\begin{equation}
H(t,p)=(t^2+|p|^2)^{ \frac 12}
\label{9.23}
\end{equation}
the corresponding Finsler metric is
$$
e^w\varphi(-\dot w; -\dot \xi)=e^w(\dot w^2+|\dot \xi|^2)^{ \frac 12}
$$
in $\Bbb R\times \Omega$ and is, in fact, an incomplete Riemannian metric.
In Case $n=1$ and $\Omega=(-R, R)$ then,
for $R>\pi$, \underline{there is no
geodesic $(w(t), \xi(t))$
starting at $(0,0)$
going to}
\underline{ the boundary of the strip $\Bbb R\times \Omega$}.
Nonetheless, for a bounded domain
$\Omega$ in $\Bbb R^n$, and
for $H$ of (\ref{9.23}), the function
$$
u(x)=
\left\{
\begin{array}{ll}
1&\mbox{if}\ d(x)\ge \frac \pi 2\\
\sin(d(x))& \mbox{if}\ d(x)\le \frac \pi 2,
\end{array}
\right.
$$
where $d(x)$ is the Euclidean distance from
$x$ to $\partial \Omega$, is a
viscosity solution of (\ref{1.1}), (\ref{1.3}).
In addition for its singular set $\Sigma$,
\begin{equation}
H^{n-1}(\Sigma)<\infty.
\label{9.24}
\end{equation}
Indeed,
$$
\Sigma=\Sigma_1\cup \Sigma_2
$$
where $\Sigma_1=\{x\in \Omega\ |\ d(x)=\frac \pi 2\}$
and $\Sigma_2=$singular set of the
distance function to $\partial \Omega$.  Since
$\Sigma_1$ is contained in the set
of the points in $\Omega$ of all straight segments going normal to the boundary
and having length $\frac \pi 2$,
$$
H^{n-1}(\Sigma_1)<\infty.
$$
And by Theorem A in the Introduction, rather, Corollary \ref{cor1},
$$
H^{n-1}(\Sigma_2)<\infty.
$$

We plan to take up the general case $\bar t=\hat t$ in a later work.

\section{Appendix A}
\setcounter{equation}{0}

\noindent{\bf About Remark \ref{rem1.1}}:
Examples with $C^{2,\alpha}$ boundary.
Now we present the examples.
We start with $n=2$. Essentially the same examples work for $n\ge 3$. 
For $0<\alpha<\alpha+3\epsilon\le 1$, let
$$
f(x)=1-\sqrt{ 1-x^2 }-g(x), \qquad x\in \Bbb R,
$$
where
$$
g(x)=|x|^{ 2+\alpha+3\epsilon}\left( 2+ \sin(|x|^{-\epsilon})\right).
$$
Clearly $f$ is smooth in $(-1,0)\cup(0,1)$ and
$$
f'(x)=x-g'(x)+O(|x|^3)=O(|x|),
$$
$$
f''(x)=1-g''(x)+O(x^2)=O(1),
$$
$$
g'(x)=O(|x|^{1+\alpha+2\epsilon}),
$$
$$
g''(x)=
 -\epsilon^2 |x|^{\alpha+\epsilon} \sin(|x|^{-\epsilon})
+O(|x|^{\alpha+2\epsilon}),
$$
and
$$
g'''(x)=O(|x|^{ \alpha-1}).
$$
It follows from the above that for any $0<x<y\le \frac 12$,
$$
|g''(x)-g''(y)|\le \int_0^1|g'''(t)|dt
\le C\int_0^1 t^{ \alpha-1}dt\le C|x-y|^\alpha.
$$
So  $g''\in C^\alpha(-\frac 12, \frac 12)$ 
and $f\in C^{2,\alpha}(-\frac 12, \frac 12)$.  We also see that
for
 some $0<\delta<\frac 12$,
$$
f''(x)\ge \frac 12,\qquad \forall\ |x|<\delta.
$$

Since $1-\sqrt{ 1-x^2}$ is a part of the graph of the 
unit circle centered at $(0,1)$ and $f(x)\le
1-\sqrt{ 1-x^2}$ with equality holds only at $x=0$, we can construct a 
strictly convex $C^{2,\alpha}$ domain $\Omega$ which
has $\{ (x, f(x))\ |\ |x|<\delta\}$
as a part of its boundary $\partial
\Omega$, and
$$
dist ( (0,1), Q)>1,
\qquad \forall\ Q\in \partial \Omega \setminus\{ (0,0) \}.
$$
See Fig. 9 below

\bigskip

\centerline{
\epsfbox{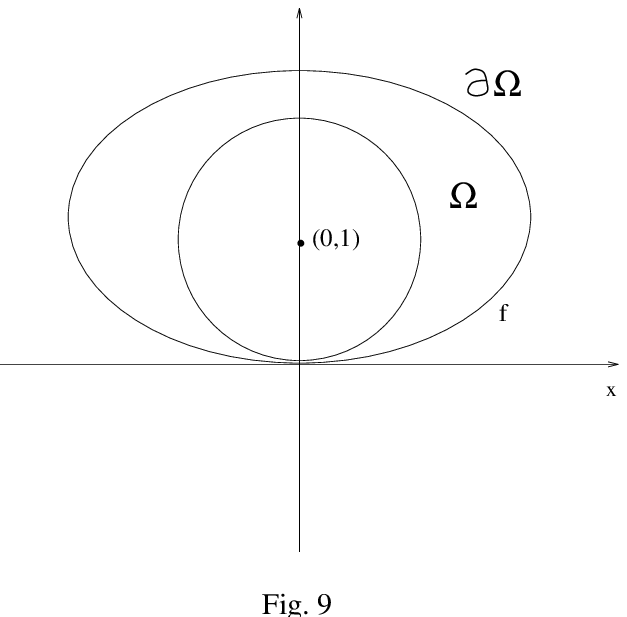}}

\bigskip

Clearly, $m(0,0)=(0,1)$.  We will show that
there exists some positive constant $c>0$ such that for
any $0<x<\delta$ satisfying $\cos(x^{-\epsilon})=0$ and
$\sin(x^{-\epsilon})=1$, we have
\begin{equation}
|m (x,f(x))-(x, f(x))|\le
1-c|x|^{\alpha+\epsilon}.
\label{A1}
\end{equation}
This implies that $m$ is not in $C^\beta$ for
any $\beta>\alpha+\epsilon$.
Indeed, for $x_k=(2k\pi +\frac \pi 2)^{ -1/\epsilon}\to 0$ as $k\to
\infty$, we have, for large $k$,
\begin{eqnarray*}
&&
|m (x_k, f(x_k))-m(0,0)|\\
&\ge &
|m(0,0)|- |m(x_k, f(x_k)) - (x_k, f(x_k))|
-|(x_k, f(x_k))|\\
&\ge & 1-( 1-c |x_k|^{ \alpha+\epsilon})-C|x_k|=
c |x_k|^{ \alpha+\epsilon}-C|x_k|
\ge  \frac c2 |x_k|^{ \alpha+\epsilon} 
\ge \frac c4 |(x_k, f(x_k))|^{  \alpha+\epsilon}.
\end{eqnarray*}

In the following we establish (\ref{A1}).
The curvature of the graph of $f$ is given by
$$
k(x)=\frac {f''(x) }{  \sqrt{ 1+f'(x)^2 } }.
$$
Thus
$$
k(x)= f''(x)+O(x^2)=1-g''(x)+O(x^2).
$$

Since $\cos( x_k^{-\epsilon} )=0$ and $\sin( x_k^{-\epsilon} )=1$, we have
$$
k(x)=1-g''(x)
+O(x^2)=1+\epsilon^2 x^{ \alpha +\epsilon} +O( x^{ \alpha+2\epsilon}).
$$
This implies that
$$
|m (x, f(x))- (x, f(x))|
\le 1- \epsilon^2 x^{ \alpha +\epsilon} +O( x^{ \alpha+2\epsilon}),
$$
from which (\ref{A1}) follows.

For $n\ge 3$, 
$$
f(x)=1-\sqrt{ 1-|x|^2 }-g(x), \qquad x\in \Bbb R^{n-1},
$$
where
$$
g(x)=|x|^{ 2+\alpha+3\epsilon}\left( 2+ \sin(|x|^{-\epsilon})\right).
$$
We still have $f\in C^{2,\alpha}$, and we can still 
construct $\Omega$ essentially the same way. 
For $x=(x_1, x_2, \cdots, x_n)$, considering the
curve, $(  (x_1, 0, \cdots, 0), f (x_1, 0, \cdots, 0))$, we
already know that for $x_1>0$, $\cos(x_1^\epsilon)=1$ and $\sin (x_1^\epsilon)
=0$,
the curvature of the curve
is  $\ge  1+c |x_1|^{ \alpha+\epsilon}$ for some 
 constant $c>0$, and therefore, for such $x_1$, 
$|m( (x_1, 0, \cdots, 0), f (x_1, 0, \cdots, 0))
-( (x_1, 0, \cdots, 0), f (x_1, 0, \cdots, 0)) |\ge
\frac c5 |x_1|^{ \alpha+\epsilon}$.
So $m$ is not in $C^\beta$
for any $\beta>\alpha +\epsilon$.

\section{Appendix B}
\setcounter{equation}{0}

\begin{lem}
Let ${\cal X}$ be the set of $k\times k$ real matrices.  For
$A\in {\cal X}$, $A$  positive definite, consider the following linear equations for $X\in {\cal X}$
$$
AX=X^TA.
$$
The dimension of the space of solutions 
is $ \frac {  k(k+1)}2$.
\label{lema-1}
\end{lem} 

\noindent {\bf Proof.}\  Let $Y=AX$.  Then the equation
takes the form 
$Y^T=Y$, i.e., $Y$ is symmetric.  The dimension of
the space of real symmetric matrices is  $ \frac {  k(k+1)}2$.

\vskip 5pt
\hfill $\Box$
\vskip 5pt

\begin{lem} Let  $A$ be a  $k\times k$
real symmetric  positive definite
matrix, and let   $D$ be  a  $k\times k$
real anti-symmetric matrix, i.e.,
$D^T=-D$.
Then the dimension of the space of solutions
 to the following
linear equations
$$
X^TA-AX=D, \qquad X\in {\cal X}
$$
is $ \frac {  k(k+1)}2$. 
\label{lema-2}
\end{lem}       

\noindent {\bf Proof.}\  Both sides of equations are anti-symmetric, so 
the number of equations is: $\frac { k(k-1)}2
$.
By Lemma \ref{lema-1}, the  dimension of the kernel is
$ \frac {  k(k+1)}2$.  The lemma follows since $
dim\ {\cal X}=k^2=
\frac { k(k-1)}2 + \frac {  k(k+1)}2$.

\vskip 5pt
\hfill $\Box$
\vskip 5pt

\section{Appendix C}
\setcounter{equation}{0}

In this appendix we give a proof of the path-connectedness
of the singular set $\Sigma$, as mentioned in the introduction.

\noindent{\bf Proof.}\
The proof is based on
Lemma \ref{lem4.7}, the continuity of the map $y\to m(y)$
for $y\in \partial \Omega$.  
Suppose $X$ and $Y$ are points in $\Sigma$. Connect
them by a smooth curve lying in $\Omega$.  It suffices
to show that if we have a smooth arc 
$x(t)$ lying in $G$ except for its end points,
$X_0, X_1$, which lie in
$\Sigma$, then $X_0$ can be joined
to $X_1$ by a continuous
arc lying in $\Sigma$.

Consider the smooth arc $x(t)$,
$0\le t\le 1$, with $X(0)=X_0$, $X(1)=X_1$.
For every $t$ in $(0,1)$ there
is a unique point $y(t)$ on $\partial \Omega$
which is the closest point on $\partial \Omega$
to $x(t)$.  Clearly $y(t)$
is a continuous curve for $0<t<1$.

As $t\to 0$, $y(t)$ need
not have a unique limit.  Choose a sequence $t_i\to 0$,
$t_{i+1}<t_i$, 
so that $y(t_i)$ converge to some point
$y_0$.  We have
$$
m(y(t_i))=x(t_i)\to X_0.
$$
For $i\ge k$, large, replace
the curve $y(t)$, for
$t_{i+1}\le t\le t_i$ by the shortest
arc on $\partial \Omega$
from $y(t_i)$ to $y(t_{i+1})$.
Continuing this for all $i\ge k$
we get a new curve $\bar y(t)
$ tending to $y_0$ as $t\to 0$,
and $m(\bar y(t))\to X_0$ as
$t\to 0$ by the continuity of the map
$y\to m(y)$.  Doing the same near
the other end point, for $t\to 1$, we
obtained the desired arc $m(\bar y(t))$ in
$\Sigma$ connecting $X_0$ to $X_1$.

\vskip 5pt
\hfill $\Box$

\end{document}